\documentclass[a4paper,11pt]{amsart}
\usepackage{amssymb}
\usepackage{amscd}
\usepackage[all]{xypic}

\title[Three-fold divisorial contractions]
{Three-fold divisorial contractions to singularities of higher indices}
\author{Masayuki Kawakita}
\address{Graduate School of Mathematical Sciences, University of Tokyo,
3-8-1 Komaba, Meguro, Tokyo 153-8914, Japan}
\email{kawakita@ms.u-tokyo.ac.jp}

\newtheorem{theorem}{Theorem}[section]
\newtheorem{lemma}[theorem]{Lemma}
\newtheorem{corollary}[theorem]{Corollary}
\newtheorem*{addendum}{Addendum}

\theoremstyle{definition}
\newtheorem{example}[theorem]{Example}
\newtheorem*{example*}{Example}

\theoremstyle{remark}
\newtheorem{step}{Step}

\numberwithin{equation}{section}

\newcommand{\bC}{{\mathbb C}}
\newcommand{\bP}{{\mathbb P}}
\newcommand{\bQ}{{\mathbb Q}}
\newcommand{\bZ}{{\mathbb Z}}
\newcommand{\cF}{{\mathcal F}}

\newcommand{\cO}{{\mathcal O}}
\newcommand{\cQ}{{\mathcal Q}}
\newcommand{\calR}{{\mathcal R}}

\newcommand{\axis}{\textrm{-}\mathrm{axis}}

\newcommand{\Img}{\operatorname{Im}}
\newcommand{\mult}{\operatorname{mult}}
\newcommand{\ord}{\operatorname{ord}}
\newcommand{\others}{\mathrm{others}}
\newcommand{\red}{\mathrm{red}}
\newcommand{\Sing}{\operatorname{Sing}}
\newcommand{\something}{\mathrm{something}}
\newcommand{\torsion}{\mathrm{torsion}}

\newcommand{\wt}{\operatorname{wt}}
\newcommand{\word}{\operatorname{w\textrm{-}ord}}

\newcommand{\rd}[1]{\lfloor{#1}\rfloor}
\newcommand{\system}[1]{\left|{#1}\right|}

\begin{document}
\maketitle

\begin{abstract}
We complete the explicit study of a three-fold divisorial contraction
whose exceptional divisor contracts to a point,
by treating the case where
the point downstairs is a singularity of index $n \ge 2$.
We prove that if this singularity is of type c$A/n$
then any such contraction is a suitable weighted blow-up;
and that if otherwise then the discrepancy is $1/n$ with a few exceptions.
Every such exception has an example.
Some exceptions allow the discrepancy to be arbitrarily large,
but any contraction in this case is described as
a weighted blow-up of a singularity of type c$D/2$
embedded into a cyclic quotient of a smooth five-fold.
The erratum to the previous paper \cite{Ka03} is attached.
\end{abstract}

\tableofcontents

\section{Introduction}\label{sec:introduction}
The minimal model program has been formulated to generalise the theory of
minimal models of surfaces to higher dimensional varieties,
see \cite{KMM87} and \cite{KM98}.
For a given variety with mild singularities,
it produces a good variety
after a finite sequence of elementary transformations
called divisorial contractions and flips.
Since Mori completed this program in dimension three
by proving the existence of three-fold flips in \cite{Mo88},
we have hoped the explicit study of three-folds
to strengthen our grasp of the theory of minimal models in higher dimension
as well as three-folds themselves.
In this paper,
we complete the explicit study of a three-fold divisorial contraction
whose exceptional divisor contracts to a point,
after my previous works \cite{Ka01}, \cite{Ka02} and \cite{Ka03}
based on the assumption that the target space is Gorenstein.
We work in the complex analytic category.

The minimal model program has to work
in a category of varieties with mild singularities
to avoid problems caused by the existence of small contractions.
Reid has introduced the notion of terminal singularities
to meet this requirement in \cite{Re83}.
We say that a normal variety $X$ has at worst \textit{terminal singularities}
if $X$ is $\bQ$-Gorenstein and
if every exceptional divisor has positive coefficient
in the discrepancy divisor $K_Y-f^*K_X$
for a resolution of singularities $f \colon Y \to X$.
Divisorial contractions are defined
in the category of varieties with terminal singularities.
Let $f \colon Y \to X$ be a proper morphism with connected fibres
between normal varieties with at worst terminal singularities.
We say that $f$ is an \textit{elementary contraction} if $-K_Y$ is $f$-ample,
and a \textit{divisorial contraction}
if moreover the exceptional locus of $f$ is a prime divisor.

We let $f \colon Y \to X$ be a three-fold divisorial contraction.
There exist two ways to study $f$,
one starting from $Y$ and the other from $X$.
The former approach is standard in the context of the minimal model program.
Mori classified them in the case where $Y$ is smooth in \cite{Mo82},
and Cutkosky extended this result
to the case where $Y$ is Gorenstein in \cite{Cu88}.
The work \cite{Mo82} of Mori has become the foundation of
the modern study of three-folds.
He proved the cone theorem in original form
and applied it to constructing an elementary contraction
from a smooth three-fold whose canonical divisor is not nef.
On the other hand,
the latter approach has turned out to be effective when $Y$ is fairly singular.
It is natural to regard $f$ as an extraction from $X$
in view of Sarkisov program,
a program which factorises a birational map between Mori fibre spaces,
because Sarkisov links of types I and II in this program start from
the converse of divisorial contractions,
see \cite{Co95}, \cite{CM00} and \cite{CPR00}.
Inspired by Sarkisov program,
we have again recognised the importance of
the explicit study of three-fold divisorial contractions.
In the case where the exceptional divisor $E$ contracts to a curve,
$f$ is if exists uniquely determined by the structure of $f(E) \subset Y$,
and we can also apply the results \cite{KM92} and \cite{Mo88}
on flipping contractions.
Note the results \cite{Tzi02a} and \cite{Tzi02p} of Tziolas in this direction.
Our object is a divisorial contraction $f$ for which $E$ contracts to a point.

From now on,
we let
\begin{align*}
f \colon (Y \supset E) \to (X \ni P)
\end{align*}
be a germ of a three-fold divisorial contraction
whose exceptional divisor $E$ contracts to a point $P$ of index $n$.
The index $n$ is defined as the smallest positive integer
such that $nK_X$ is Cartier.
We let $a/n$ denote the discrepancy of $f$,
which is defined by the relation $K_Y=f^*K_X+(a/n)E$.
Kawamata proved that $f$ is a certain weighted blow-up determined uniquely
when $P$ is a terminal quotient singularity in \cite{Km96},
and Corti proved that $f$ is the blow-up
when $P$ is an ordinary double point in \cite{Co00}.
I systematically studied $f$
in the case where $X$ is Gorenstein in \cite{Ka03},
after the classification results \cite{Ka01} and \cite{Ka02}.
I described them when $P$ is a c$A$ point,
and bounded the discrepancy by three when $P$ is a c$E$ point.
Unfortunately one case has been omitted in \cite{Ka03}
when $P$ is a c$D$ point,
the correction to which is supplied at the end of this paper.
In this paper
we treat the case where $X$ is not Gorenstein
to complete the study of $f$,
and deduce the following classification theorem.

\begin{theorem}\label{thm:MT}
Let $f\colon (Y \supset E) \to (X \ni P)$
be a germ of a three-fold divisorial contraction
whose exceptional divisor $E$ contracts to a point $P$ of index $n \ge 2$.
\begin{enumerate}
\item\label{itm:MTord}
If $P$ is of type c$A/n$
\textup{(}see Section \textup{\ref{sec:numerical}} for definition\textup{)},
then there exists a suitable identification
\begin{align*}
P \in X \cong o \in (x_1x_2+g(x_3^n,x_4)=0) \subset
\bC^4_{x_1x_2x_3x_4} /\frac{1}{n}(1,-1,b,0)
\end{align*}
such that $f$ is the weighted blow-up with weights
$\wt(x_1,x_2,x_3,x_4)=(r_1/n,r_2/n,a/n,1)$
which satisfies the following conditions.
\begin{enumerate}
\item
$a \equiv br_1$ modulo $n$ and $r_1+r_2 \equiv 0$ modulo $an$.
\item
$(a-br_1)/n$ is co-prime to $r_1$.
\item
$g$ has weighted order $(r_1+r_2)/n$
with weights $\wt(x_3,x_4)=(a/n,1)$.
\item
The monomial $x_3^{(r_1+r_2)/a}$ appears in $g$
with non-zero coefficient.
\end{enumerate}
Moreover, any such $f$ is a divisorial contraction.
\item\label{itm:MTexc}
If $P$ is of type other than c$A/n$,
then the discrepancy $a/n$ is $1/n$ except for the cases
listed in Table \textup{\ref{tbl:MT}}.
Moreover,
there exists an example of $f$ in every case in Table \textup{\ref{tbl:MT}},
see Examples \textup{\ref{exl:cE/2}}, \textup{\ref{exl:cD/2_14-1}},
\textup{\ref{exl:cD/2_15}}, \textup{\ref{exl:cD/2_14-2}} and
\textup{\ref{exl:cD/2_15g}}.
\begin{table}[ht]
\caption{}\label{tbl:MT}
\begin{tabular}{l|c|l}
\hline
\multicolumn{1}{c|}{$P$} & $a/n$ &
\multicolumn{1}{c}{\textup{singularities on $Y$
at which $E$ is not Cartier}} \\
\hline
\textup{c$D/2$}  & $1$ &
\textup{one point $\frac{1}{2r'}(1,-1,r'+2)$,
$r'$ odd, $r'\neq1$} \\[1.5pt]
\textup{c$D/2$}  & $1$ &
\textup{one point deforming to $2 \times \frac{1}{2r'}(1,-1,r'+1)$,
$r'$ even} \\[1.5pt]
\textup{c$D/2$}  & $2$ &
\textup{one point $\frac{1}{2r'}(1,-1,r'+4)$,
$r'\equiv\pm1$ modulo $8$, $r'\neq1$} \\[1.5pt]
\textup{c$D/2$}  & $a/2$ &
\textup{two points $\frac{1}{r}(1,-1,\frac{a+r}{2})$,
$\frac{1}{r+2}(1,-1,\frac{a+r+2}{2})$,
$a,r$ odd} \\[1.5pt]
\textup{c$D/2$}  & $a/2$ &
\textup{two points $\frac{1}{r}(1,-1,\frac{a+r}{2})$,
$\frac{1}{r+4}(1,-1,\frac{a+r+4}{2})$,
$a+r$ even} \\[1.5pt]
\textup{c$E/2$}  & $1$ &
\textup{two points $\frac{1}{6}(1,5,5)$, $\frac{1}{2}(1,1,1)$} \\[1.5pt]
\hline
\end{tabular}
\end{table}
\end{enumerate}
\end{theorem}

We can consult the works \cite{Ha99} and \cite{Ha00} of Hayakawa
when $f$ has discrepancy $1/n$
in the exceptional case (\ref{itm:MTexc}) of Theorem \ref{thm:MT}.
He studied the divisors with minimal discrepancy
over a three-fold non-Gorenstein terminal singularity in detail.
On the other hand,
the discrepancy $a/n$ is arbitrarily large
in the cases in the fourth and fifth lines of Table \ref{tbl:MT}.
It is desirable to provide the explicit description of $f$ in these cases.

\begin{theorem}\label{thm:cD/2}
Let $f$ be a germ of a three-fold divisorial contraction
to a point $P$ of type c$D/2$ with discrepancy $a/2>1/2$.
\begin{enumerate}
\item\label{itm:cD/2_2}
If $f$ belongs to the case in the fourth line of Table \textup{\ref{tbl:MT}},
then there exists a suitable identification of $P \in X$ with
\begin{align*}
o \in (\phi:=x_1^2+x_1x_3q(x_3^2,x_4)+x_2^2x_4+\lambda x_2x_3^{2\alpha-1}
+p(x_3^2,x_4)=0)
\end{align*}
in the space $\bC^4_{x_1x_2x_3x_4}/\frac{1}{2}(1,1,1,0)$,
such that $f$ is the weighted blow-up with weights
$\wt(x_1,x_2,x_3,x_4)=((r+2)/2,r/2,a/2,1)$
which satisfies the following conditions.
\begin{enumerate}
\item
$a$, $r$ are odd and $a \mid r+1$.
\item
$\phi$ has weighted order $r+1$ with the weights distributed above.
Its weighted homogeneous part of weight $r+1$ is irreducible.
In fact $x_1x_3q$ has weighted order $r+1$ unless $q=0$.
\end{enumerate}
\item\label{itm:cD/2_4}
If $f$ belongs to the case in the fifth line of Table \textup{\ref{tbl:MT}},
then there exists a suitable identification of $P \in X$ with
\begin{align*}
o \in \bigg(
\begin{array}{c}
\phi_1:=x_1^2+x_2x_5+p(x_3^2,x_4)=0 \\
\phi_2:=x_2x_4+x_3^{(r+2)/a}+q(x_3^2,x_4)x_3x_4+x_5=0
\end{array}
\bigg)
\end{align*}
in the space $\bC^5_{x_1x_2x_3x_4x_5}/\frac{1}{2}(1,1,1,0,1)$,
such that $f$ is the weighted blow-up with weights
$\wt(x_1,x_2,x_3,x_4,x_5)=((r+2)/2,r/2,a/2,1,(r+4)/2)$
which satisfies the following conditions.
\begin{enumerate}
\item
$a \mid r+2$ with $a\neq r+2$ and $(r+2)/a$ is odd.
\item
$\phi_1$ has weighted order $r+2$ with the weights distributed above.
Similarly $qx_3x_4$ has weighted order $(r+2)/2$ unless $q=0$.
\end{enumerate}
\end{enumerate}
\end{theorem}

Although Theorems \ref{thm:MT} and \ref{thm:cD/2}
are the strongest statements,
we emphasise that they are intrinsically obtained as corollaries to
the conceptual Theorems \ref{thm:co-prime} and \ref{thm:GE} explained below.
In contrast to the Gorenstein case,
the discrepancy $a/n$ is no longer an integer.
The annoying problem is that
$a$ and $n$ may have a common divisor $g$,
and consequently there may exist
a numerically trivial but linearly non-trivial
$\bQ$-Cartier divisor $(n/g)K_Y-(a/g)E$.
This phenomenon really happens as indicated in Theorem \ref{thm:MT},
but we can still expect the co-primeness of $a$ and $n$
by observing the weights of semi-invariant coordinates
of the index-one cover of a three-fold terminal singularity.
We prove this co-primeness with the exceptions listed in Table \ref{tbl:MT}
when $f$ is of exceptional type.

\begin{theorem}\label{thm:co-prime}
Let $f \colon Y \to X$
be a three-fold divisorial contraction
whose exceptional divisor contracts to a point of index $n$,
and $a/n$ the discrepancy of $f$.
If $f$ is of exceptional type
\textup{(}see Section \textup{\ref{sec:numerical}} for definition\textup{)},
then $a$ is co-prime to $n$
except for the cases in Table \textup{\ref{tbl:MT}}.
\end{theorem}

The main theorem in \cite{Ka03} is
the affirmation of the \textit{general elephant conjecture}
for a three-fold divisorial contraction to a Gorenstein point,
and the classification is deduced from it.
This conjecture has been proposed by Reid in \cite{Re87}.
From the observation of the classification of
three-fold terminal singularities
by himself in \cite{Re83} and Mori in \cite{Mo85},
he pointed out that
a general element in the anti-canonical system,
dubbed a \textit{general elephant}, of a three-fold
should have at worst Du Val singularities
in appropriate situations involving an elementary contraction.
The remarkable point is that
this approach has settled the problem of the existence of three-fold flips
in \cite{Mo88}.
We prove that this conjecture holds also for our divisorial contractions.
Combining the results in \cite{KM92} and \cite{Mo88}
for the flipping/divisor-to-curve contraction case,
we finally obtain the following theorem.

\begin{theorem}\label{thm:GE}
Let $f \colon Y \to X$
be a germ of a three-fold elementary birational contraction
whose central fibre is irreducible.
Then a general elephant of $Y$ has at worst Du Val singularities.
\end{theorem}

We have positive answers \cite{Re83p}, \cite{Sh79} and \cite{Ta02}
also for Fano three-folds with mild singularities,
whereas there exist counter-examples in \cite{ABR03}
when the anti-canonical system has small dimension.
Note the result \cite{Al94} of Alexeev in this direction.

We shall explain how to prove these theorems.
Our investigation of $f$ is, as in the case where $X$ is Gorenstein,
founded on a classification theorem in the numerical sense,
Theorem \ref{thm:classification},
which is obtained by
the singular Riemann--Roch formula \cite[Theorem 10.2]{Re87} of Reid and
a relative vanishing theorem \cite[Theorem 1-2-5]{KMM87} of Kawamata--Viehweg.
Theorem \ref{thm:classification} classifies $f$ into types according to
the non-Gorenstein singularities on $Y$,
and we use the dichotomy that
$f$ is either of exceptional type or of general type,
following this theorem.
It is worth while mentioning that
Theorem \ref{thm:classification} does not encode
any information on the singularities on $Y$ at which $E$ is Cartier.
This is the result of the fact that,
because of the locality of $f$,
we need to take the difference along $E$ of Euler--Poincar\'e characteristics
in applying the singular Riemann--Roch formula.
Although this lack of information generates
larger number of types for $f$ than in the Gorenstein case,
it is not difficult to handle most of the exceptional cases.

The method of approaching the co-primeness of $a$ and $n$ is simple.
If $a$ and $n$ have a common divisor $g$,
then the index-$(n/g)$ cover $(X' \ni P') \to (X \ni P)$,
the covering defined by the $\bQ$-Cartier divisor $(n/g)K_X$,
is lifted to the crepant covering morphism $(Y' \supset E') \to (Y \supset E)$.
The induced morphism $f' \colon (Y' \supset E') \to (X' \ni P')$
satisfies the conditions to be a divisorial contraction
except for the irreducibility of $E'$,
and we can apply the numerical classification Theorem \ref{thm:classification}
to $f'$,
which narrows down the possibility that $a$ is not co-prime to $n$.
Then, in most cases,
Theorem \ref{thm:co-prime} is deduced from the numerical analysis
which uses the numerically but not linearly trivial divisor $(n/g)K_Y-(a/g)E$.
As seen in Corollary \ref{cor:a=1} and Theorem \ref{thm:easy_GE},
it is fairly obvious that for $f$ of exceptional type
the co-primeness implies the equality $a=1$ and the general elephant theorem.

On the other hand,
most of the study of $f$ of general type is devoted to
the proof of the general elephant Theorem \ref{thm:hard_GE} in strong form.
In the case $a/n>1$,
we prove it by analysing the scheme-theoretic intersection $C$
of $E$ and the birational transform of a general hyperplane section on $X$,
as demonstrated in \cite{Ka03}.
This $C$ is a (possibly non-reduced) tree of $\bP^1$,
and the way of the embedding of $C$ into $Y$ encodes
the behaviour of the general elephants of $Y$.
It should be emphasised that we need to keep paying attention to
the hidden non-Gorentstein points on $Y$,
the non-Gorenstein points at which $E$ is Cartier,
throughout the geometric investigation of $Y$.
In the case $a/n\le1$,
the general elephant theorem itself is obvious,
but the complete proof of Theorem \ref{thm:hard_GE} involves
very complicated division into cases and
delicate geometric investigation.

We obtain Theorems \ref{thm:MT}(\ref{itm:MTord}) and \ref{thm:cD/2}
on the explicit description of $f$
by specifying $E$ as a valuation on the function field.
Starting with the explicit description of the germ $P \in X$
provided by Mori in \cite{Mo85},
we find good coordinates of this germ in the sense
that $f$ must be a certain weighted blow-up
with weights distributed to these coordinates,
thanks to the general elephant theorem.
These coordinates are constructed to meet the structure of the graded ring
$\bigoplus_{i\ge0} f_*\cO_Y(-iE)/f_*\cO_Y(-(i+1)E)$
in the part of lower degrees,
determined by the numerical results from the singular Riemann--Roch formula.
We are convinced that once a three-fold singularity is explicitly given,
we can classify all the three-fold divisorial contractions to this singularity
by combining this method with
the singular Riemann--Roch technique and the general elephant theorem.
The method of explicit description used for proving Theorem \ref{thm:cD/2}
is applied to whichever case classified in Theorem \ref{thm:classification}
by formal modification.
It is enough to perform such a remaining classification,
where the discrepancy is very small,
in the course of the application to the global study of three-folds.

The analysis of the cases in Table \ref{tbl:MT} yields
interesting examples of $f$.
In particular for a divisorial contraction $f$ to a point $P$ of type c$D/2$,
it is not always true that $f$ is described
as a weighted blow-up of a hypersurface
in the quotient space $\bC^4/\frac{1}{2}(1,1,1,0)$.
The description of $f$ involves the identification of $P \in X$
with a cyclic quotient of a complete intersection three-fold in $\bC^5$.
Such a complication in the c$D/2$ case is recognised
also in the study \cite{Ha00} of divisors with minimal discrepancy by Hayakawa.

\begin{figure}[b]
\caption{}
\begin{align*}
\xymatrix@!0@=15.8pt{
&&&&&&&&&&
*+[F]\txt{Numerical classification Theorem \ref{thm:classification}}
\ar@{-}[dd]
&&&&&&&&\\
&&&&&&&&&&&&&&&&\\
*=<0pt>{}\ar@{-}[dddddddd]^*\txt{some \\ cases}
&&&&&
*=<0pt>{}\ar[dd]\ar@{-}[lllll]
&&&&&
*=<0pt>{}\ar@{-}[lllll]_>*\txt{exceptional type \quad}
\ar@{-}[rrrrr]^>*\txt{general type}
&&&&&
*=<0pt>{}\ar@{-}[dd]
&&&\\
&&&&&&&&&&&&&&&&&&\\
&&&&&
*+[F]\txt{Co-primeness by \\ Theorems
\ref{thm:co-prime_almost}, \ref{thm:co-prime_8}, \ref{thm:14/15}}
\ar[ddd]
&&&&&&
*=<0pt>{}\ar[dd]
&&&&
*=<0pt>{}\ar@{-}[llll]_>*{a/n\le1}\ar@{-}[rrr]^>*{a/n>1}
&&&
*=<0pt>{}\ar[dddd]^*\txt{Theorems \\{}\ref{thm:hard_GE16}, \ref{thm:hard_GE15}}
\\
&&&&&&&&&&&&&&&&&&\\
&&&&&&&&&&&
*+[F]\txt{Easy general elephant \\ Theorem \ref{thm:easy_GE}}\ar@{-}[dd]
&&&&&&&\\
&&&&&
*+[F]\txt{$a=1$ by \\ Corollary \ref{cor:a=1}}\ar[ddd]
&&&&&&&&&&&&&\\
&&&&&&&&&&&
*=<0pt>{}\ar[rrrrrrr]^<<<<<<<<*{a>1}_<<<<<<<<<*\txt{Theorem\\
{}\ref{thm:15/16small}}
\ar[ddd]_<<<<<*{a=1}
&&&&&&&
*+[F]\txt{Hard general elephant \\ Theorem \ref{thm:hard_GE}}
\ar@{-}[ddd]
\\
&&&&&&&&&&&&&&&&&&\\
*=<0pt>{}\ar[rrrrr]\ar@{-}[ddd]
&&&&&
*+[F]\txt{Easy general elephant \\ Theorem \ref{thm:easy_GE}}
&&&&&&&&&&&&&\\
&&&&&&&&&&&
*+[F]{a=1}
&&
*=<0pt>{}\ar@{-}[dd]
&&&&&
*=<0pt>{}
\ar@{-}[lllll]_>>>>>>>>>>>*\txt{Theorem \ref{thm:hard_GE}(\ref{itm:hard_A})}
\ar@{-}[dd]^<<*\txt{Theorem\\
{}\ref{thm:hard_GE}(\ref{itm:hard_D})}
\\
&&&&&&&&&&&&&&&&&&\\
*=<0pt>{}\ar[rrrrrr]_<<<<<<<*\txt{Theorems\\
{}\ref{thm:co-prime_8}, \ref{thm:14/15}}
&&&&&&
*+[F]\txt{$P$ type c$D/2$, c$E/2$ \\ $(a,n)=(2,2), (4,2)$}\ar[ddd]
&&&&&&&
*+[F]\txt{$P$ type c$A/n$}\ar[ddd]
&&&&&
*+[F]\txt{$P$ type c$D/2$\\{}by Corollary \ref{cor:cD/2}}\ar[ddd]
\\
&&&&&&&&&&&&&&&&&&\\
&&&&&&&&&&&&&&&&&&\\
&&&&&&
*+[F]\txt{Examples\\{}\ref{exl:cE/2}, \ref{exl:cD/2_14-1},
\ref{exl:cD/2_15}, \ref{exl:cD/2_14-2}
}
&&&&&&&
*+[F]\txt{Classification\\{}by Section \ref{sec:cA/n},\\$a=1$ included}
&&&&&
*+[F]\txt{Classification\\{}by Section \ref{sec:cD/2}}\ar[dd]
\\
&&&&&&&&&&&&&&&&&&\\
&&&&&&&&&&&&&&&&&&
*+[F]\txt{Example \ref{exl:cD/2_15g}}
}\
\end{align*}
\end{figure}

The structure of this paper is as follows.
In Section \ref{sec:numerical},
we deduce the numerical classification Theorem \ref{thm:classification}
and some auxiliary results concerned with this theorem.
In Section \ref{sec:co-primeness},
we discuss the co-primeness and prove Theorem \ref{thm:co-prime} except for
a few cases to be treated in Section \ref{sec:small}.
These two sections are devoted to the numerical analysis,
while the rest are devoted to the geometric analysis.
The general elephant theorem is proved in Section \ref{sec:GE}.
Most of this section is allocated to the proof of Theorem \ref{thm:hard_GE}
in the case $a/n>1$.
To complete Theorems \ref{thm:co-prime} and \ref{thm:hard_GE},
we treat some remaining cases quite explicitly in Section \ref{sec:small}.
The last two sections are devoted to the classification of $f$.
In Section \ref{sec:cA/n},
we present the explicit description of divisorial contractions
to points of type c$A/n$.
In Section \ref{sec:cD/2},
we complete Theorems \ref{thm:MT}(\ref{itm:MTexc}) and \ref{thm:cD/2}
by studying divisorial contractions of general type to points of type c$D/2$.
Examples for the cases in Table \ref{tbl:MT}
are provided in Sections \ref{sec:small} and \ref{sec:cD/2}.
Finally,
the erratum to the previous paper \cite{Ka03} is attached
at the end of this paper.
All the corrections are made by the methods used in this paper.
The framework of the proof of
Theorems \ref{thm:MT}, \ref{thm:cD/2}, \ref{thm:co-prime} and \ref{thm:GE}
is summarised in the flowchart.

I would like to thank Professor Yujiro Kawamata
for his stimulating encouragement.
Financial support was provided by
the Japan Society for the Promotion of Science.

\section{Numerical results}\label{sec:numerical}
Our study of three-fold divisorial contractions relies on
the classification of three-fold terminal singularities.
We let $P \in X$ be a three-fold germ.
We say that $P$ is a \textit{cDV} (\textit{compound Du Val}) \textit{point}
if a general hyperplane section has at worst a Du Val singularity at $P$.
The singularity $P$ is said to be
\textit{c$A_n$, c$D_n$, c$E_n$} (\textit{compound $A_n, D_n, E_n$})
according to the type of the Du Val singularity on
a general hyperplane section.
The characterisation of three-fold Gorenstein terminal singularities
is due to Reid.

\begin{theorem}[{\cite{Re83}}]\label{thm:Re83}
Let $P \in X$ be a three-fold germ.
$P$ is a Gorenstein terminal singularity if and only if
$P$ is an isolated cDV point.
\end{theorem}

We then consider a three-fold non-Gorenstein terminal singularity $P \in X$.
Let $n$ denote the index of $P \in X$,
that is, the smallest positive integer such that $nK_X$ is Cartier at $P$.
Take the index-one cover $\pi \colon (X^\sharp \ni P^\sharp) \to (X \ni P)$,
which is a cyclic $\mu_n$-cover.
Mori gave a precise description of this covering,
and the classification was completed by
Koll\'ar and Shepherd-Barron in \cite[Theorem 6.4]{KS88}.

\begin{theorem}[{\cite{Mo85}}]\label{thm:Mo85}
There exists a $\mu_n$-equivariant identification
\begin{align*}
P^\sharp \in X^\sharp \cong o \in (\phi=0) \subset \bC^4_{x_1x_2x_3x_4},
\end{align*}
where $x_1, x_2, x_3, x_4$ and $\phi$ are $\mu_n$-semi-invariant.
$P$ belongs to one of the types in Table \textup{\ref{tbl:Mo85}}.

\begin{table}[ht]
\caption{}\label{tbl:Mo85}
\begin{tabular}{l|c|l|l}
\hline
\multicolumn{1}{c|}{\textup{type}} & $n$ & \multicolumn{1}{c|}{$\phi$} &
\multicolumn{1}{c}{$\wt(x_1,x_2,x_3,x_4; \phi)$} \\
\hline
\textup{c$A/n$}  & $n$ &
$x_1x_2+g(x_3^n,x_4)$ & $(1,-1,b,0;0)$ \textup{$b$, $n$ co-prime} \\
\textup{c$Ax/2$} & $2$ &
$x_1^2+x_4^2+\phi_{\ge 4}(x_2,x_3)$   & $(1,1,1,0;0)$ \\
\textup{c$Ax/4$} & $4$ &
$x_1^2+x_2^2+\phi_{\ge 2}(x_3^2,x_4)$ & $(1,3,1,2;2)$ \\
\textup{c$D/2$}  & $2$ &
$x_1^2+\phi_{\ge 3}(x_2,x_3,x_4)$     & $(1,1,1,0;0)$ \\
\textup{c$D/3$}  & $3$ &
$x_4^2+\phi_{\ge 3}(x_1,x_2,x_3)$     & $(1,2,1,0;0)$ \\
\textup{c$E/2$}  & $2$ &
$x_1^2+x_4^3+ \phi_{\ge 4}(x_2,x_3)$  & $(1,1,1,0;0)$ \\
\hline
\end{tabular}
\end{table}
\end{theorem}

Reid made notable observations on
the geometric aspects of a three-fold terminal singularity $P \in X$
in \cite{Re87}.
He pointed out that a general member of $\system{-K_X}$
on a germ of any such singularity has at worst a Du Val singularity at $P$,
and then proposed the general elephant conjecture
that a general elephant of a three-fold
should have at worst Du Val singularities
in appropriate situations involving an elementary contraction.
We repeat that a \textit{general elephant} is by definition
a general element in the anti-canonical system.
We let $S_X$ be a general elephant on a germ $P \in X$ of
a three-fold non-Gorenstein terminal singularity.
Its pre-image $S_X^\sharp:=\pi^{-1}(S_X)_\red$
on the index-one cover $P^\sharp \in X^\sharp$ is a $\bQ$-Cartier divisor,
and satisfies $\mathrm{S}_2$ condition by \cite[Corollary 5.25]{KM98}.
Hence the covering $S_X^\sharp \to S_X$ is a crepant morphism
between normal surfaces with Du Val singularities.
We reproduce the list of this covering in Table \ref{tbl:reproduce}.

\begin{table}[ht]
\caption{}\label{tbl:reproduce}
\begin{tabular}{l|rcl}
\hline
\multicolumn{1}{c|}{type} &
\multicolumn{3}{c}{covering} \\
\hline
c$A/n$  & $A_{k-1}$  & $\to$ & $A_{kn-1}$ \\
c$Ax/2$ & $A_{2k-1}$ & $\to$ & $D_{k+2}$  \\
c$Ax/4$ & $A_{2k-2}$ & $\to$ & $D_{2k+1}$ \\
c$D/2$  & $D_{k+1}$  & $\to$ & $D_{2k}$   \\
c$D/3$  & $D_4$      & $\to$ & $E_6$      \\
c$E/2$  & $E_6$      & $\to$ & $E_7$      \\
\hline
\end{tabular}
\end{table}

Reid also noted that
any three-fold terminal singularity $P \in X$ has a small deformation
to a basket of terminal quotient singularities $P_i$,
called \textit{fictitious singularities},
although Mori had already observed it in \cite{Mo85}.
This basket plays a pivotal role
in the singular Riemann--Roch formula \cite[Theorem 10.2]{Re87} of Reid.
If $P$ is of type other than c$Ax/4$,
then any local index at $P_i$ equals that at $P$,
and if $P$ is of type c$Ax/4$,
then one of $P_i$ has local index $4$ and all the others have local index $2$.

We return to the study of three-fold divisorial contractions.
We let
\begin{align*}
f \colon (Y \supset E) \to (X \ni P)
\end{align*}
be a germ of a three-fold divisorial contraction
whose exceptional divisor $E$ contracts to a point $P$ of index $n$.
We let $a/n$ denote its discrepancy
defined by the relation $K_Y=f^*K_X+(a/n)E$.
We shall expand the numerical results
\cite[Section 4]{Ka01} and \cite[Section 2]{Ka03} into our general situation.
Unlike the case where $X$ is Gorenstein,
in general $E$ does not generate
the torsion part of the local class group at a point of $Y$.
By this reason,
we need to consider both the divisors $K_Y$ and $E$ simultaneously.
We set
\begin{align}\label{eqn:definition_D}
D_{i,j}:=iK_Y+jE
\end{align}
for integers $i$, $j$.
There exists a natural exact sequence
\begin{align}\label{eqn:exact_sequence}
0 \to \cO_Y(D_{i,j-1}) \to \cO_Y(D_{i,j}) \to \cQ_{i,j} \to 0.
\end{align}
The quotient sheaf $\cQ_{i,j}$ is an $\cO_E$-module
satisfying $\mathrm{S}_2$ condition by \cite[Proposition 5.26]{KM98}.
Note that $\cQ_{0,0}=\cO_E$ and $\cQ_{1,1}=\omega_E$.
We set
\begin{align*}
d(i,j):= \chi(\cQ_{i,j})=\chi(\cO_Y(D_{i,j}))-\chi(\cO_Y(D_{i,j-1})).
\end{align*}
We need to extend $f$ to a morphism of proper varieties
by the algebraisation theorems \cite{Ar68} and \cite{Ar69} of Artin
when we regard $d(i, j)$ as the difference
$\chi(\cO_Y(D_{i,j}))-\chi(\cO_Y(D_{i,j-1}))$.
Note that the value $d(i,j)$ depends only on
the class of $(i,j)$ in $\bZ^2/\bZ(n, -a)$.

We shall compute $d(i,j)$ by the singular Riemann--Roch formula.
We let $I_0:=\{Q,\ \textrm{of type} \ \frac{1}{r_Q}(1, -1, b_Q)\}$
denote the basket of fictitious singularities from singularities on $Y$,
and $e_Q$ for $Q \in I_0$ the smallest positive integer
such that $E \sim e_QK_Y$ at $Q$.
Such an $e_Q$ exists by \cite[Corollary 5.2]{Km88}.
Note that $n \equiv ae_Q$ modulo $r_Q$ by the global relation $nK_Y \sim aE$.
By replacing $b_Q$ with $r_Q-b_Q$ if necessary,
we may assume that $v_Q:= \overline{e_Qb_Q}\le r_Q / 2$,
where $\bar{\ }$ denotes the residue modulo $r_Q$,
that is, $\overline{k}:=k-\rd{k/r_Q}r_Q$.
We set $I:=\{Q\in I_0 \mid v_Q\neq0\}$.
If a non-Gorenstein point $Q$ on $Y$ is of type other than c$Ax/4$,
then $v_{Q_i}$ takes the same value
for every fictitious singularity $Q_i$ from $Q$.
If $Q$ is of type c$Ax/4$,
then $v_{Q_0}$ is either $\overline{\pm1}$, $2$ or $0$
for the fictitious singularity $Q_0$ of index $4$,
and any other $v_{Q_i}$ is respectively $1$, $0$, $0$.

By this notation, the singular Riemann--Roch formula implies that
\begin{align}\label{eqn:d(i,j)}
d(i,j)&= \frac{1}{12}
\biggl(6\Bigl(\frac{a}{n}i+j\Bigr)^2 \\
\nonumber& \qquad \qquad -6\Bigl(\frac{a}{n}+1\Bigr)\Bigl(\frac{a}{n}i+j\Bigr)
+\Bigl(\frac{a}{n}+1\Bigr)\Bigl(\frac{a}{n}+2\Bigr)\biggr)E^3 \\
\nonumber&\quad\ +\frac{1}{12} E\cdot c_2(Y) +A_{i,j}-A_{i,j-1}.
\end{align}
Here the contribution term $A_{i,j}$ is given by
$A_{i,j}:=\sum_{Q \in I}A_Q(\overline{i+je_Q})$,
where
\begin{align*}
A_Q(k):=-k\frac{r_Q^2-1}{12r_Q}+\sum_{l=1}^{k-1}
\frac{\overline{lb_Q}(r_Q-\overline{lb_Q})}{2r_Q}.
\end{align*}
$A_Q(k)$ is determined by $k$ modulo $r_Q$.
We emphasise that
the formula (\ref{eqn:d(i,j)}) no longer encodes
any information on the fictitious singularities at which $E$ is Cartier.
This is never the case when $X$ is Gorenstein,
but in general the set $I$ does not necessarily coincide with $I_0$.
We say that a non-Gorenstein point $Q$ on $Y$ is
a \textit{hidden non-Gorenstein point} if $E$ is Cartier at $Q$.
Note that the local index at any hidden non-Gorenstein point
is a divisor of $n$ by the relation $nK_Y \sim aE$.
It is difficult to control the value of $A_Q(k)$ itself,
but the following formulae make $A_Q(k)$ more accessible.
\begin{align}
A_Q(k+1)-A_Q(k)&=-\frac{r_Q^2-1}{12r_Q}+B_{r_Q}(kb_Q),
\label{eqn:difference_of_A} \\
A_Q(k)+A_Q(-k)&=-B_{r_Q}(kb_Q),
\end{align}
where $B_r(k):=(\overline{k}\cdot\overline{r-k})/2r$.
By (\ref{eqn:d(i,j)}) and (\ref{eqn:difference_of_A}),
we obtain that
\begin{align}\label{eqn:difference_of_d(i,j)}
d(i+1,j)-d(i,j)=\Bigl(\frac{a}{n}i+j-\frac{1}{2}\Bigr)\frac{a}{n}E^3
+B_{i,j}-B_{i,j-1},
\end{align}
where $B_{i,j}:=\sum_{Q \in I}B_{r_Q}(ib_Q+jv_Q)$.

On the other hand,
the relative vanishing theorem \cite[Theorem 1-2-5]{KMM87} implies that
$R^kf_*\cO_Y(D_{i,j})=0$ for $k\ge1$ whenever $ia/n+j \le a/n$.
In particular by (\ref{eqn:exact_sequence}), if $ia/n+j \le a/n$ then
\begin{align}\label{eqn:H^k(Q)}
H^k(\cQ_{i,j})=
\begin{cases}
f_*\cO_Y(D_{i,j})/f_*\cO_Y(D_{i,j}-E) & \textrm{for} \ k=0,\\
0                                     & \textrm{for} \ k\ge 1.
\end{cases}
\end{align}
Therefore
\begin{align}\label{eqn:d(i,j)=dim}
d(i,j)=\dim f_*\cO_Y(D_{i,j})/f_*\cO_Y(D_{i,j}-E)
\quad \textrm{if} \ ia/n+j \le a/n.
\end{align}
The elements of the $\cO_X$-module $f_*\cO_Y(D_{i,j})/f_*\cO_Y(D_{i,j}-E)$
correspond to the effective divisors on $X$ linearly equivalent to $iK_X$
which have multiplicity $-(ia/n+j)$ along $E$.
In particular,
\begin{align}\label{eqn:d(i,j)_non-negative}
d(i,j)=
\begin{cases}
0 & \textrm{if} \ 0 \le ia/n+j \le a/n \
    \textrm{and} \ (i,j)\not\in\bZ(n, -a),\\
1 & \textrm{if} \ (i,j)\in\bZ(n, -a).
\end{cases}
\end{align}
We explain (\ref{eqn:d(i,j)_non-negative}) in the case where $ia/n+j=0$.
Set $g:=\gcd(a, n)$ and $a=ga'$, $n=gn'$.
Then $D_{n',-a'}=n'f^*K_X$ is numerically trivial
and $g$ is the smallest positive integer such that
$gD_{n',-a'}$ is Cartier.
We can write $(i,j)=k(n',-a')$ by some integer $k$ if $ia/n+j=0$.
Then
\begin{align*}
d(i,j)=\dim \cO_X(kn'K_X)/f_*\cO_Y(kn'f^*K_X-E),
\end{align*}
whence $d(i, j)=1$ if and only if $kn'K_X$ is Cartier at $P$,
or equivalently $(i,j)\in\bZ(n, -a)$.

We classify $f$ according to the numerical data $J:=\{(r_Q,v_Q)\}_{Q\in I}$.

\begin{theorem}\label{thm:classification}
$f$ belongs to one of the types in Table \textup{\ref{tbl:classification}}.

\begin{table}[ht]
\caption{}\label{tbl:classification}
\begin{tabular}{c|l|c}
\hline
\textup{No} & \multicolumn{1}{c|}{$J$}                     & $(a/n)E^3$    \\
\hline
\textup{1}  & $(6,3)$                                      & $1/2$         \\
\textup{2}  & $(7,3)$                                      & $2/7$         \\
\textup{3}  & $(8,3)$                                      & $1/8$         \\
\textup{4}  & \textup{$(4,2)$, $(r,1)$}                    & $1/r$         \\
\textup{5}  & \textup{$(5,2)$, $(2,1)$}                    & $3/10$        \\
\textup{6}  & \textup{$(5,2)$, $(3,1)$}                    & $2/15$        \\
\textup{7}  & \textup{$(5,2)$, $(4,1)$}                    & $1/20$        \\
\textup{8}  & \textup{$(6,2)$, $(2,1)$}                    & $1/6$         \\
\textup{9}  & \textup{$(7,2)$, $(2,1)$}                    & $1/14$        \\
\hline
\end{tabular}
\
\begin{tabular}{c|l|c}
\hline
\textup{No} & \multicolumn{1}{c|}{$J$}                     & $(a/n)E^3$    \\
\hline
\textup{10} & \textup{$(2,1)$, $(2,1)$, $(r,1)$}           & $1/r$         \\
\textup{11} & \textup{$(2,1)$, $(3,1)$, $(3,1)$}           & $1/6$         \\
\textup{12} & \textup{$(2,1)$, $(3,1)$, $(4,1)$}           & $1/12$        \\
\textup{13} & \textup{$(2,1)$, $(3,1)$, $(5,1)$}           & $1/30$        \\
\textup{14} & $(r,2)$                                      & $4/r$         \\
\textup{15} & \textup{$(r_1,1)$, $(r_2,1)$, $r_1 \le r_2$} & $1/r_1+1/r_2$ \\
\textup{16} & $(r,1)$                                      & $1+1/r$       \\
\textup{17} & $\emptyset$                                  & $2$           \\
            &                                              &               \\
\hline
\end{tabular}
\end{table}
\end{theorem}

\begin{proof}
$d(1,0)-d(0,0)=-1$ by (\ref{eqn:d(i,j)_non-negative}).
Applying (\ref{eqn:difference_of_d(i,j)}) for $(i,j)=(0,0)$, we have
\begin{align*}
-1=
-\frac{1}{2}\Bigl(\frac{a}{n}E^3\Bigr)-\sum_{Q \in I}\frac{v_Q(r_Q-v_Q)}{2r_Q}.
\end{align*}
All the possibilities of $J$ satisfying this equality are
listed in Table \ref{tbl:classification}.
\end{proof}

We divide No 15 in Table \ref{tbl:classification} into cases for convenience.
We say that $f$ is of No $\textrm{15}'$, No $\textrm{15}''$
if $J$ comes from one, respectively two non-Gorenstein points on $Y$.
No $\textrm{15}'$ is further divided into two cases:
No $\textrm{15}'$a if $J=\{(r, 1), (r, 1)\}$,
and No $\textrm{15}'$b if $J=\{(2, 1), (4, 1)\}$.
We say that $f$ is \textit{of exceptional type}
if $f$ is of No 1-14, $\textrm{15}'$ or 17,
and \textit{of general type}
if $f$ is of No $\textrm{15}''$ or 16.

We shall present some auxiliary results based on
the singular Riemann--Roch formula and the relative vanishing theorem.
We let $R$ denote the least common multiple of all the $r_Q$ for $Q\in I$,
and $R^*$ that of all the $r_Q/\gcd(r_Q, v_Q)$ for $Q\in I$.

\begin{lemma}\label{lem:integer}
$R^*E^3$ and $(a/n)^2RE^3$ are integers.
\end{lemma}

\begin{proof}
By (\ref{eqn:d(i,j)}),
$d(0,R^*+j)-d(0,j)=(j+(R^*-a/n-1)/2)R^*E^3$.
Thus $R^*E^3$ is an integer.
By (\ref{eqn:difference_of_d(i,j)}),
$(d(R+i+1,j)-d(R+i,j))-(d(i+1,j)-d(i,j))=(a/n)^2RE^3$.
Thus $(a/n)^2RE^3$ is an integer.
\end{proof}

The following is an easy corollary to Lemma \ref{lem:integer}.

\begin{corollary}\label{cor:a=1}
Suppose that $f$ is of exceptional type with index $n \ge 2$.
If $a$ is co-prime to $n$, then $a=1$.
\end{corollary}

\begin{proof}
For instance,
we consider the case No $\textrm{15}'$a.
Then $R^*E^3=2n/a$ and $(a/n)^2RE^3=2a/n$.
Hence $a=1$ by the assumption and Lemma \ref{lem:integer}.
\end{proof}

We shall consider a non-trivial effective divisor $L_{i,j}$ on $Y$
which is linearly equivalent to $D_{i,j}$ and intersects $E$ properly.
Note that such an $L_{i,j}$ exists if and only if
$f_*\cO_Y(D_{i,j}-E) \neq f_*\cO_Y(D_{i,j})$ and $D_{i,j} \not\sim 0$.
In particular,
the assumption that $ia/n+j<0$ is implicitly imposed.
By the following exact sequences,
we construct an $\cO_{L_{i,j} \cap E}$-module $\calR_{i,j}^{i',j'}$
which is naturally isomorphic to
$\cO_Y(D_{i',j'}) \otimes \cO_{L_{i,j} \cap E}$ outside finitely many points.
\begin{gather}
0 \to \cO_Y(D_{i',j'}-E) \to \cO_Y(D_{i',j'}) \to \cQ_{i',j'} \to 0, \\
0 \to \cQ_{i'-i,j'-j} \to \cQ_{i',j'} \to \calR_{i,j}^{i',j'} \to 0,
\label{eqn:exact_R}
\end{gather}
\begin{align}
\begin{array}{rcccl}
& \cO_Y(D_{i',j'}) & & & \\
& \downarrow  & & & \\
0 \to (\textrm{finite length}) \to &
\cO_Y(D_{i',j'}) \otimes \cO_E & \to & \cQ_{i',j'} & \to 0 \\
& \downarrow & & \downarrow & \\
0 \to (\textrm{finite length}) \to &
\cO_Y(D_{i',j'}) \otimes \cO_{L_{i,j} \cap E} & \to &
\calR_{i,j}^{i',j'} & \to 0.
\end{array}
\end{align}
The map $\cO_Y(-L_{i,j}) \otimes \cO_E \to \cQ_{-i,-j}$ is surjective
because it factors the surjective map
$\cO_Y(-L_{i,j}) \twoheadrightarrow \cQ_{-i,-j}$.
Hence $\cQ_{-i,-j}$ is the ideal sheaf of $L_{i,j} \cap E \subset E$,
and thus $\calR_{i,j}^{0,0}=\cO_{L_{i,j} \cap E}$.

\begin{lemma}\label{lem:QandR}
\begin{enumerate}
\item\label{itm:QandR_Q}
$f_*\cO_Y(D_{i,j}) \twoheadrightarrow H^0(\cQ_{i,j})$ for any $i$, $j$.
\begin{align*}
&h^1(\cQ_{i,j})=0 \quad
\textrm{unless} \ a/n<ia/n+j<1.\\
&h^2(\cQ_{i,j})=
\begin{cases}
0 & \textrm{if} \ ia/n+j \le a/n,\\
d(1-i,1-j) & \textrm{if} \ ia/n+j \ge 1.
\end{cases}
\end{align*}
\item\label{itm:QandR_R}
Consider $L_{i,j}$ as above.
Then,
\begin{align*}
H^0(\cQ_{i',j'}) \twoheadrightarrow H^0(\calR_{i,j}^{i',j'}) \
\textrm{unless} \ a/n<(i'-i)a/n+(j'-j)<1.
\end{align*}
\begin{align*}
h^1(\calR_{i,j}^{i',j'})&=
\begin{cases}
0 & \textrm{if} \ (i'-i)a/n+(j'-j) \le a/n+1 \\
  & \quad \textrm{and} \ (i'-i-1,j'-j-1)\not\in\bZ(n, -a),\\
1 & \textrm{if} \ (i'-i-1,j'-j-1)\in\bZ(n, -a)
\end{cases}
\end{align*}
unless $a/n<i'a/n+j'<1$.
\end{enumerate}
\end{lemma}

\begin{proof}
The lemma is derived from
the Grothendieck duality $H^k(\cQ_{i,j})=H^{2-k}(\cQ_{1-i,1-j})$
and (\ref{eqn:H^k(Q)}), (\ref{eqn:d(i,j)=dim}),
(\ref{eqn:d(i,j)_non-negative}), (\ref{eqn:exact_R}).
\end{proof}

\begin{corollary}\label{cor:P1}
Suppose that $-(ia/n+j)\le a/n+1$ with $(i+1,j+1)\not\in\bZ(n, -a)$
and that the estimate $a/n<-(ia/n+j)<1$ does not hold.
Then $h^0(\cO_{L_{i,j} \cap E})=1$ and $h^1(\cO_{L_{i,j} \cap E})=0$.
In particular,
$L_{i,j} \cap E$ has no embedded points and
$(L_{i,j} \cap E)_\red$ is a tree of $\bP^1$.
\end{corollary}

The way of the embedding of $L_{i,j} \cap E$ into $Y$ encodes
the behaviour of the general elephants of $Y$.
According to Corollary \ref{cor:P1},
the fundamental ingredients we should study are smooth curves on $Y$.
We let $C$ be a smooth curve on $Y$,
and $Q$ a quotient singularity $\bC^3_{x_1x_2x_3}/\frac{1}{r_Q}(1, -1, b_Q)$
on $Y$ which $C$ passes through.
We shall review the local study of the germ $Q \in C \subset Y$
in \cite[Section 2]{Ka03} after \cite[Section 2]{Mo88}.

We take the index-one cover $\pi \colon (Y^\sharp \ni Q^\sharp) \to (Y \ni Q)$
and set $C^\sharp:=(C \times_Y Y^\sharp)_\red$.
Let $n_Q^C$ denote the number of the irreducible components of $C^\sharp$,
and $Q^\dag \in C^\dag$ the normalisation of
one of the irreducible components of $C^\sharp$.
We let $t \in \cO_{C, Q}$ and $t^{n_Q^C/r_Q} \in \cO_{C^\dag, Q^\dag}$ be
uniformising parameters of $C$ and $C^\dag$.
Let $a_i$ denote the smallest positive integer such that
there exists a semi-invariant function with weight $\wt x_i$
whose image in $\cO_{C^\dag, Q^\dag}$ has
order $a_i/r_Q$ with respect to $t$.
Then the coordinates of $Y^\sharp$ can be taken so that
$x_i |_{C^\dag}=t^{a_i/r_Q}$ for $i=1,2,3$.
Such a description is called a \textit{normal form} of $Q \in C \subset Y$.
Note that
\begin{align*}
(a_1,a_2,a_3) \in \bZ \cdot n_Q^C(1,-1,b_Q) +r_Q\bZ^3.
\end{align*}
For $l \in \bZ/(r_Q)$,
we let $w_Q^C(l)$ denote the smallest non-negative integer
such that $(w_Q^C(l),l) \in \bZ \times \bZ/(r_Q)$ is contained in
the semi-group
\begin{align}\label{eqn:semigroup}
G_Q^C:=\bZ_{\ge0}(a_1,1)+\bZ_{\ge0}(a_2,-1)+\bZ_{\ge0}(a_3,b_Q)
\subset \bZ \times \bZ/(r_Q).
\end{align}
We note that $(w,l) \in G_Q^C$ means the existence of
a semi-invariant function with weight $l$ whose restriction to $C^\dag$
is $t^{w/r_Q}$.
In particular $w_Q^C(0)=0$ and $(r_Q, 0) \in G_Q^C$.

Let $\cF$ denote the reflexive sheaf on the germ $Q \in Y$
which is isomorphic to the ideal sheaf defined by $x_1=0$ outside $Q$.
Note that
$\cF^{[\otimes v_Q]} \cong \cO_Y(E)$ and
$\cF^{[\otimes b_Q]} \cong \cO_Y(K_Y)$,
where $\cF^{[\otimes i]}$ denotes the double dual of $\cF^{\otimes i}$.
Then for arbitrary integers $j_1, \ldots, j_k$
with $\sum_{1 \le i \le k} j_i=0$,
\begin{multline}\label{eqn:how.to.comp}
\Img \Bigl(\cF^{[\otimes j_1]} \otimes \cdots \otimes
\cF^{[\otimes j_k]} \otimes \cO_C \to \cO_C\Bigr) \\
= \cO_C\Bigl(-\bigl({\sum_{1 \le i \le k} w_Q^C(-j_i)/r_Q}\bigl) \cdot Q\Bigr).
\end{multline}

\section{Co-primeness}\label{sec:co-primeness}
In this section we discuss the co-primeness of $a$ and $n$
and prove Theorem \ref{thm:co-prime}
except for a few cases to be treated in Section \ref{sec:small}.
We start with listing the types of $f$ for which
$a$ and $n$ may have a non-trivial common divisor.

Let $p$ be a common divisor of $a$ and $n$,
and set $a=pa_p$, $n=pn_p$.
We take the index-$n_p$ cover $\alpha_X \colon (X' \ni P') \to (X \ni P)$,
the covering defined by $n_pK_X$,
and lift it over $Y$ by taking the normalisation of $Y$
in the function field of $X'$.
By this procedure,
we obtain the commutative diagram
\begin{align}\label{eqn:covering}
\begin{CD}
(Y' \supset E') @>\alpha_Y>> (Y \supset E) \\
@Vf'VV                       @VVfV \\
(X' \ni P')     @>\alpha_X>> (X \ni P),
\end{CD}
\end{align}
where $E':=\alpha_Y^*E$.
The morphism $\alpha_Y$ is just the covering
defined by $n_pf^*K_X=n_pK_Y-a_pE$.
In particular $\alpha_Y$ is \'etale outside finitely many points,
which implies that $K_Y'=\alpha_Y^*K_Y$ and
that $E'$ is an $f'$-anti-ample reduced divisor.
Hence $Y'$ has also at worst terminal singularities,
and $f'$ satisfies the conditions to be a divisorial contraction
except for the irreducibility of $E'$.

Now we should notice that all the results in Section \ref{sec:numerical} hold
without the condition of the irreducibility of $E$,
in contrast to the use of the reducedness of $E$ to deduce that $d(0,0)=1$.
By applying Theorem \ref{thm:classification} to both $f$ and $f'$,
we can restrict the possibilities of the covering (\ref{eqn:covering}).
It is worth while mentioning that,
although we can construct the covering (\ref{eqn:covering})
by using $(n/q)K_X$ for any divisor $q$ of $n$,
the morphism $\alpha_Y$ is no longer \'etale
at the generic point of $E$ unless $q$ is also a divisor of $a$.
Once $\alpha_Y$ ramifies at the generic point of $E$,
there is no guarantee that $Y'$ has terminal singularities,
and we can not apply the singular Riemann--Roch formula to $Y'$.

We return to our situation.
The covering (\ref{eqn:covering}) is a cyclic $\mu_p$-cover for
a common divisor $p$ of $a$ and $n$.
In order to compare the numerical data of $f$ and $f'$,
we examine the variations of
the basket of fictitious singularities and the value of $(a/n)E^3$
by this covering.
It is obvious that
\begin{align}\label{eqn:variation}
\frac{a_p}{n_p}{E'}^3=p\cdot\frac{a}{n}E^3.
\end{align}
We shall investigate the variation of the basket.
Since the covering (\ref{eqn:covering}) is decomposed into
coverings of prime degrees,
we assume that $p$ is a prime divisor for the time being.
We let $Q$ be a singularity of $Y$.
There exist two cases of the local behaviour of $\alpha_Y$ at $Q$.
If $\alpha_Y$ is \'etale at $Q$,
then each fictitious singularity
$Q_i$ of type $\frac{1}{r_{Q_i}}(1,-1,b_{Q_i})$ from $Q$
is divided into $p$ fictitious singularities
of type $\frac{1}{r_{Q_i}}(1,-1,b_{Q_i})$ on $Y'$.
If $\alpha_Y$ is ramified at $Q$,
then $p \mid r_Q$ and
each $Q_i$ varies into a fictitious singularity
of type $\frac{1}{r_{Q_i}/p}(1,-1,b_{Q_i})$ on $Y'$,
provided that $Q$ is of type other than c$Ax/4$.
If $Q$ is of type c$Ax/4$,
then $p=2$ and the basket
$\{\frac{1}{4}(1,-1,\pm1),s\times\frac{1}{2}(1,1,1)\}$
varies into $\{(2s+1)\times\frac{1}{2}(1,1,1)\}$.
On the other hand,
the value $e_Q$ is preserved modulo index
since $K_Y'=\alpha_Y^*K_Y$ and $E'=\alpha_Y^*E$.
Therefore we have obtained the following lemma on $J=\{(r_Q,v_Q)\}_{Q\in I}$.

\begin{lemma}\label{lem:variation}
Suppose that $p$ is a common prime divisor of $a$ and $n$.
Let $Q \in I_0$ be a fictitious singularity of $Y$
equipped with the data $(r_Q, v_Q)$.
\begin{enumerate}
\item\label{itm:variation_etale}
If $\alpha_Y$ is \'etale at $Q$,
then there exist exactly $p$ fictitious singularities $Q'$ of $Y'$ above $Q$.
Each such $Q'$ is equipped with the data $(r_{Q'},v_{Q'})=(r_Q,v_Q)$.
\item\label{itm:variation_ramified}
If $\alpha_Y$ is ramified at $Q$,
then there exists exactly one fictitious singularity $Q'$ of $Y'$ above $Q$.
$Q'$ is equipped with the data $(r_{Q'},v_{Q'})=(r_Q/p,\overline{\pm v_Q})$,
where $\overline{\pm v_Q}$ is
the minimum of $\overline{v_Q}$ and $\overline{-v_Q}$.
\end{enumerate}
\end{lemma}

We list all the possibilities of $\alpha_Y$ in terms of $J$
by Theorem \ref{thm:classification}, (\ref{eqn:variation})
and Lemma \ref{lem:variation}.

\begin{theorem}\label{thm:ramification}
Suppose that $p$ is a common prime divisor of $a$ and $n$.
Let $J':=\{(r_Q',v_Q')\}_{Q'\in I'}$ be the data of $Y'$ constructed
in the same fashion as $J$.
Then $\alpha_Y$ belongs to one of the cases
in Table \textup{\ref{tbl:ramification}}.

\begin{table}[ht]
\caption{}\label{tbl:ramification}
\begin{tabular}{c|l|c|l}
\hline
\textup{No of $f$}
&\multicolumn{1}{c|}{$J$}    & $p$ & \multicolumn{1}{c}{$J'$}    \\
\hline
\textup{1}  & $(6,3)$                             &
$3$ & $(2,1)$                              \\
\textup{2}  & $(7,3)$                             &
$7$ & $\emptyset$                          \\
\textup{4}  & \textup{$(4,2)$, $(r,1)$}           &
$2$ & \textup{$(r,1)$, $(r,1)$}            \\
\textup{8}  & \textup{$(6,2)$, $(2,1)$}           &
$2$ & \textup{$(3,1)$, $(2,1)$, $(2,1)$}   \\
            &                                     &
$3$ & \textup{$(2,1)$, $(2,1)$, $(2,1)$}   \\
\textup{10} & \textup{$(2,1)$, $(2,1)$, $(r,1)$}  &
$2$ & \textup{$(2,1)$, $(2,1)$, $(r/2,1)$} \\
            &                                     &
$2$ & \textup{$(r,1)$, $(r,1)$}            \\
\textup{11} & \textup{$(2,1)$, $(3,1)$, $(3,1)$}  &
$3$ & \textup{$(2,1)$, $(2,1)$, $(2,1)$}   \\
\textup{12} & \textup{$(2,1)$, $(3,1)$, $(4,1)$}  &
$2$ & \textup{$(3,1)$, $(3,1)$, $(2,1)$}   \\
\textup{14} & \textup{$(r,2)$}                    &
$p$ & \textup{$(r/p,2)$, $r/p \ge 2$}      \\
\textup{$\textrm{15}'$a, $\textrm{15}''$} & \textup{$(r_1,1)$, $(r_2,1)$} &
$p$ & \textup{$(r_1/p,1)$, $(r_2/p,1)$}    \\
\hline
\end{tabular}
\end{table}
\end{theorem}

\begin{proof}
For instance, we consider the case No 1.
Then $(a/n)E^3=1/2$.
Let $Q$ be the fictitious singularity of $Y$ such that $(r_Q,v_Q)=(6,3)$.
Because $J'$ contains at most one $(6,3)$ by Table \ref{tbl:classification},
$\alpha_Y$ must be ramified at $Q$ by Lemma \ref{lem:variation},
whence $p=2$ or $3$.
If $p=2$ then $J'=\emptyset$ and $(a_p/n_p){E'}^3=1$
by (\ref{eqn:variation}) and Lemma \ref{lem:variation},
but it is impossible by the list of the possible types of $f'$
in Table \ref{tbl:classification}.
If $p=3$ then $J'=\{(2,1)\}$ and $(a_p/n_p){E'}^3=3/2$.
In this case $f'$ is of No 16.
\end{proof}

If $a$ is not co-prime to $n$, then in particular $a\neq1$.
The following numerical lemma is derived from
the results \cite{Km92} and \cite{Ma96} on minimal discrepancies by
Kawamata in the non-Gorenstein case
and by Markushevich in the Gorenstein case.

\begin{lemma}\label{lem:minimal_discrepancy}
Suppose that $X$ is singular at $P$ and that $a \neq 1$.
Then there exist a singularity $Q \in Y$ of index $r_Q$
at which $E \sim e_QK_Y$ and
positive integers $n_m$, $n_a$ satisfying the following conditions.
\begin{enumerate}
\item\label{itm:minimal_discrepancy_1}
$an_m+nn_a=r_Q$.
\item\label{itm:minimal_discrepancy_0}
$n_m+e_Qn_a \equiv 0$ modulo $r_Q$.
\end{enumerate}
In particular, $a+n\le r_Q$.
\end{lemma}

\begin{proof}
The results \cite{Km92} and \cite{Ma96} by Kawamata and Markushevich
are the following.
For a three-fold terminal singularity $P \in X$ of index $n$,
there exists an exceptional divisor above $P$
whose discrepancy with respect to $K_X$ is $1/n$
unless $X$ is smooth at $P$.
Hence if $a\neq1$ then there exist a birational morphism $g \colon Z \to Y$
from a smooth variety $Z$ and a $g$-exceptional prime divisor $F$ such that
$K_Z=(f \circ g)^*K_X+(1/n)F+(\others)$.
$g(F)$ has to be a non-Gorenstein point $Q$ of index say $r_Q$.
Let $e_Q$ be an integer such that $E \sim e_QK_Y$ at $Q$.
Set
\begin{align*}
K_Z &=g^*K_Y + (n_a/r_Q)F+(\others), \\
g^*E &=g^{-1}_*E+(n_m/r_Q)F+(\others).
\end{align*}
Then the equalities
(\ref{itm:minimal_discrepancy_1}) and (\ref{itm:minimal_discrepancy_0}) hold.
\end{proof}

\begin{corollary}\label{cor:ramification_again}
\begin{enumerate}
\item
$a$ is co-prime to $n$ unless $f$ is of No 4, 8, 10, 12, 14 $\textrm{15}'$a
or $\textrm{15}''$.\label{itm:ramification_again_ge}
\item
If $f$ is of No 8 or 12,
then $a$ is co-prime to $n$ or $(a,n)=(2,2)$.\label{itm:ramification_again_sp}
\end{enumerate}
\end{corollary}

\begin{proof}
For (\ref{itm:ramification_again_ge}),
it suffices to exclude the cases No 1, 2 and 11
from Table \ref{tbl:ramification}.
For instance, we consider the case No 1.
$a$ and $n$ are divided by $3$,
and $Q$ in Lemma \ref{lem:minimal_discrepancy} must be equipped with
the data $(r_Q, v_Q)=(6, 3)$ and $e_Q=3$.
Thus $an_m/r_Q \ge 3/2$ by
(\ref{itm:minimal_discrepancy_0}) in Lemma \ref{lem:minimal_discrepancy},
but this contradicts
(\ref{itm:minimal_discrepancy_1}) in Lemma \ref{lem:minimal_discrepancy}.
(\ref{itm:ramification_again_sp}) is proved in the same way.
\end{proof}

For a common divisor $p$ of $a$ and $n$,
we set
\begin{align*}
d_Q^p:=\overline{n_p-e_Qa_p}.
\end{align*}
Then $pd_Q^p \equiv 0$ modulo $r_Q$ by $nK_Y \sim aE$,
and $\alpha_Y$ is locally the covering defined by the divisor
$n_pf^*K_X=n_pK_Y-a_pE \sim d_Q^pK_Y$ at $Q$.
In Lemma \ref{lem:variation},
(\ref{itm:variation_etale}) occurs when $d_Q^p \equiv 0$ modulo $r_Q$,
while (\ref{itm:variation_ramified}) occurs when
$pd_Q^p \equiv 0$ but $d_Q^p \not\equiv 0$ modulo $r_Q$.
We obtain Table \ref{tbl:ramification_again} which lists
all the possible coverings of prime degrees,
based on Lemma \ref{lem:variation},
Theorem \ref{thm:ramification} and Corollary \ref{cor:ramification_again}.

\begin{table}[ht]
\caption{}\label{tbl:ramification_again}
\begin{tabular}{c|c|l|l}
\hline
No of $f$ & $p$ & \multicolumn{1}{c|}{$Q_i \in I$ with $(r_{Q_i},v_{Q_i})$} &
\multicolumn{1}{c}{$d_i:=d_{Q_i}^p$} \\
\hline
4  & $2$ & $Q_1$: $(4,2)$, $Q_2$: $(r,1)$    & $(d_1,d_2)=(2,0)$ \\
8  & $2$ & $Q_1$: $(6,2)$, $Q_2$: $(2,1)$    & $(d_1,d_2)=(3,0)$ \\
10 & $2$ & $Q_1$: $(2,1)$, $Q_2$: $(2,1)$, $Q_3$: $(r,1)$ &
$(d_1,d_2,d_3)=(0,1,r/2)$ \\
   & $2$ & & $(d_1,d_2,d_3)=(1,1,0)$ \\
12 & $2$ & $Q_1$: $(2,1)$, $Q_2$: $(3,1)$, $Q_3$: $(4,1)$ &
$(d_1,d_2,d_3)=(1,0,2)$ \\
14 & $p$ & $Q$: $(r,2)$                       & \\
$\textrm{15}'$a, $\textrm{15}''$
   & $p$ & $Q_1$: $(r_1,1)$, $Q_2$: $(r_2,1)$ & \\
\hline
\end{tabular}
\end{table}

We examine No 14 and 15 without the assumption that $p$ is a prime.
Recall the notation that $g=\gcd(a, n)$, $a=ga'$, $n=gn'$.

\begin{lemma}\label{lem:d_Q}
\begin{enumerate}
\item\label{itm:d_Q_sp}
If $f$ is of No 14 or $\textrm{15}'$a with $g\ge2$,
then $g=2$ and $d_Q^2=r_Q/2$ for any $Q \in I$.
\item\label{itm:d_Q_ge}
If $f$ is of No $\textrm{15}''$ with $g\ge2$,
then $\overline{d_{Q_1}^gb_1}/r_1+\overline{d_{Q_2}^gb_2}/r_2=1$.
\end{enumerate}
\end{lemma}

\begin{proof}
By (\ref{eqn:difference_of_d(i,j)}),
we have
\begin{multline*}
(d(n'+1,-a')-d(n',-a'))-(d(1,0)-d(0,0)) \\
= (B_{n',-a'}-B_{n',-a'-1})+B_{0,-1}.
\end{multline*}
Suppose that $g\ge2$.
By (\ref{eqn:d(i,j)_non-negative}) the left-hand side is $1$,
whence
\begin{align}\label{eqn:d_Q}
1=\sum_{Q \in I}
\bigl(B_{r_Q}(d_Q^gb_Q)-B_{r_Q}(d_Q^gb_Q-v_Q)+B_{r_Q}(-v_Q)\bigr).
\end{align}

We apply (\ref{eqn:d_Q}) to No 14 and 15.
If $f$ is of No 14 or 15,
then by Lemma \ref{lem:variation} and Theorem \ref{thm:ramification},
the covering $\alpha_Y$ of degree $g$ constructed in (\ref{eqn:covering})
is totally ramified at every fictitious singularity $Q$ in $I$ of $Y$.
Hence $d_Q^g/(r_Q/g)$ is co-prime to $g$.
In No 14, $I$ consists of $Q$ with $(r_Q,v_Q)=(r,2)$,
and $d_Q^g \ge 2$ by $r/g \ge 2$ in Theorem \ref{thm:ramification}.
In No 15, $I$ consists of $Q_1$ with $(r_1,1)$ and $Q_2$ with $(r_2,1)$,
and $d_{Q_1}^g, d_{Q_2}^g \ge 1$.
Note that $r_1=r_2$ and $d_1=d_2$ in No $\textrm{15}'$a.
When $f$ is of No 14 or 15, the equality (\ref{eqn:d_Q}) becomes
\begin{align*}
1 &= \sum_{Q \in I}
\Biggl( \frac{\overline{d_Q^gb_Q}(r_Q-\overline{d_Q^gb_Q})}{2r_Q} \\
&\qquad \qquad \qquad
-\frac{(\overline{d_Q^gb_Q}-v_Q)(r_Q-\overline{d_Q^gb_Q}+v_Q)}{2r_Q}
+\frac{v_Q(r_Q-v_Q)}{2r_Q} \Biggr) \\
&= \sum_{Q \in I} \frac{v_Q(r_Q-\overline{d_Q^gb_Q})}{r_Q}.
\end{align*}
The lemma is deduced from this equality.
\end{proof}

The idea of the proof of Lemma \ref{lem:d_Q} is to examine
the variation of the value of $B_{i,j}-B_{i,j-1}$
by the replacement of $(i,j)$ with $(i+n',j-a')$.
We shall pay attention to $A_{i,j}-A_{i,j-1}$ instead.
This value encodes more information,
although it is more difficult to handle.
By (\ref{eqn:d(i,j)}) we obtain that
\begin{align*}
d(0,0)-d(n_p,-a_p)=-A_{0,-1}-(A_{n_p,-a_p}-A_{n_p,-a_p-1}).
\end{align*}
If $p\ge2$ then
\begin{align}\label{eqn:d_Q_again}
1=\sum_{Q \in I}
\bigl(-A_Q(-e_Q)-A_Q(d_Q^p)+A_Q(d_Q^p-e_Q)\bigr).
\end{align}

\begin{lemma}\label{lem:co-prime_4,10,12,15-2}
$a$ is co-prime to $n$
if $f$ is of No 4, 10 or 12.
\end{lemma}

\begin{proof}
We first remark that in the case No 10 in Table \ref{tbl:ramification_again},
$e_{Q_1}=e_{Q_2}=1$ and thus $d_1=d_2=\overline{a'-n'}$,
whence $(d_1,d_2,d_3)$ has to be $(1,1,0)$.

Suppose that $a$ and $n$ are not co-prime.
Then they have a common divisor $p=2$.
Starting with Table \ref{tbl:ramification_again} and
Lemma \ref{lem:d_Q}(\ref{itm:d_Q_sp}),
we compute the term $c_i:=-A_{Q_i}(-e_{Q_i})-A_{Q_i}(d_i)+A_{Q_i}(d_i-e_{Q_i})$
which appears in the right-hand side of (\ref{eqn:d_Q_again}),
where $d_i:=d_{Q_i}^2$.
\begin{enumerate}
\item
In No 4, $(c_1,c_2)=(1/2,0)$.
\item
In No 10, $(c_1,c_2,c_3)=(1/4,1/4,0)$.
\item
In No 12, $(c_1,c_2)=(1/4,0)$.
$c_3=0$ if $b_{Q_3}=1$ and $c_3=1/2$ if $b_{Q_3}=3$.
\end{enumerate}
Hence (\ref{eqn:d_Q_again}) never holds,
which is a contradiction.
\end{proof}

We prove the co-primeness of $a$ and $n$ when $f$ is of exceptional type
except for a few cases to be investigated in Section \ref{sec:small},
where the values of $a$ and $n$ are small.

\begin{theorem}\label{thm:co-prime_almost}
If $f$ is of exceptional type,
then $a$ is co-prime to $n$ except for the following possibilities.
\begin{enumerate}
\item
$f$ is of No 8 and $(a,n)=(2,2)$.
\item
$f$ is of No 14 and $(a,n)=(2,2)$ or $(4,2)$.
\item
$f$ is of No $\textrm{15}'$a and
$(a,n)=(2,2)$ or $(4,2)$.
\end{enumerate}
\end{theorem}

\begin{proof}
By Corollary \ref{cor:ramification_again} and
Lemma \ref{lem:co-prime_4,10,12,15-2},
we have only to treat the case where $f$ is of No 14 or $\textrm{15}'$a.
In this case $g=1$ or $2$ by Lemma \ref{lem:d_Q}(\ref{itm:d_Q_sp}).
Suppose that $g=2$.
$a=2a'$, $n=2n'$, and $a'$ is co-prime to $n'$.
We have the following by Table \ref{tbl:ramification}.
\begin{enumerate}
\item
If $f$ is of No 14 then
$R^*E^3=2n'/a'$ and $(a/n)^2RE^3=4a'/n'$.
\item
If $f$ is of No $\textrm{15}'$a then
$R^*E^3=2n'/a'$ and $(a/n)^2RE^3=2a'/n'$.
\end{enumerate}
Hence
$(a,n)=(2,2)$, $(2,4)$, $(2,8)$ or $(4,2)$ in No 14
and $(a,n)=(2,2)$, $(2,4)$ or $(4,2)$ in No $\textrm{15}'$a
by Lemma \ref{lem:integer}.
In order to exclude the case $n \neq 2$,
we compute $d(2,0)-d(1,0)$ by (\ref{eqn:difference_of_d(i,j)})
in terms of $r_Q,b_Q$ for $Q \in I$.
This value is, modulo $\bZ$, equal to
\begin{align*}
B_{r_Q}(b_Q)-B_{r_Q}(b_Q-2) &\equiv\frac{2-2b_Q}{r_Q}
&& \textrm{for No 14, $(a,n)=(2,4)$}, \\
-\frac{1}{r_Q}+B_{r_Q}(b_Q)-B_{r_Q}(b_Q-2) &\equiv\frac{1-2b_Q}{r_Q}
&& \textrm{for No 14, $(a,n)=(2,8)$}, \\
2B_{r_Q}(b_Q)-2B_{r_Q}(b_Q-1) &\equiv\frac{1-2b_Q}{r_Q}
&& \textrm{for No $\textrm{15}'$a, $(a,n)=(2,4)$}.
\end{align*}
Note that $r_Q$ is even by Theorem \ref{thm:ramification}.
The value above becomes an integer only if
$f$ is of No 14 with $(a,n)=(2,4)$ and $b_Q=1$ or $r_Q/2+1$.
In this case $e_Q=r_Q/2+2$ by Lemma \ref{lem:d_Q}(\ref{itm:d_Q_sp}),
and thus $4 \mid r_Q$.
This contradicts that $v_Q=\overline{e_Qb_Q}=2$.
\end{proof}

To treat the cases left in Theorem \ref{thm:co-prime_almost},
we are forced to analyse $f$ geometrically
by making use of Corollary \ref{cor:easy_GE_exceptional}.
This is the reason why
we postpone completing the proof of Theorem \ref{thm:co-prime}.
We here collect the numerical data to be used.

\begin{lemma}\label{lem:data}
The cases left in Theorem \textup{\ref{thm:co-prime_almost}}
have the following data.
\begin{enumerate}
\item
$f$ is of No 8.
$I$ consists of $Q=(6,2)$ and $Q'=(2,1)$.
$(a,n)=(2,2)$.
$e_Q=4$ and $b_Q=5$.
$d(0,-1)=0$ and $d(0,-2)=1$.
\item
$f$ is of No 14.
$I$ consists of $Q=(2r',2)$.
$n=2$.
$(a,b_Q)$ is either $(2,r'+2)$ with $r'\ge3$ or $(4,r'+4)$ with $r'\ge5$.
$d(0,-1)=1$ and $d(-1,1)=0$.
\item
$f$ is of No $\textrm{15}'$a.
$I$ consists of two fictitious singularities with data $(2r',1)$
from one singularity $Q$ of $Y$.
$n=2$.
$(a,b_Q)$ is either $(2,r'+1)$ with $r'\ge2$ or $(4,r'+2)$ with $r'\ge3$.
$d(0,-1)=1$ and $d(-1,1)=0$.
\end{enumerate}
\end{lemma}

\begin{proof}
All the data are deduced from Theorem \ref{thm:classification},
(\ref{eqn:difference_of_d(i,j)}), (\ref{eqn:d(i,j)_non-negative}),
Table \ref{tbl:ramification_again} and
Lemmata  \ref{lem:minimal_discrepancy}, \ref{lem:d_Q}(\ref{itm:d_Q_sp}).
Lemma \ref{lem:minimal_discrepancy} is used to estimate $r'$.
\end{proof}

\section{General elephant theorem}\label{sec:GE}
In this section we prove the general elephant Theorem \ref{thm:GE}
in strong form by either of two different ways
according to the value of the discrepancy $a/n$.

The theorem in the strongest form is easily deduced from its local version
when the discrepancy $a/n$ is small.
In this case we can obtain a general elephant of $Y$ by taking
the birational transform of a general elephant of $X$.

\begin{theorem}\label{thm:easy_GE}
If $d(-1, i)=0$ for every $1 \le i < a/n$,
then the birational transform on $Y$ of a general elephant of $X$
is a general elephant of $Y$ and has at worst Du Val singularities.
This is the case in particular when $a/n\le1$.
\end{theorem}

\begin{proof}
Let $S_X$ be a general elephant of $X$ and
$S$ its birational transform on $Y$.
$S_X$ has at worst a Du Val singularity, or equivalently,
a canonical singularity at $P$.
We set $f^*S_X=S+(a/n+m)E$ by an integer $m>-a/n$.
The assumption on $d(-1, i)$ is by (\ref{eqn:d(i,j)=dim})
equivalent to $m \ge 0$.
On the other hand, by the adjunction formula,
$\omega_S$ is in codimension one equal to $f^*\omega_{S_X} \otimes \cO_Y(-mE)$,
whence $m \le 0$ since $S_X$ is canonical.
Thus $m=0$ and $S \to S_X$ is a crepant morphism.
Now the theorem follows obviously.
\end{proof}

Theorem \ref{thm:easy_GE} is applied to the case where
$f$ is of exceptional type
by Corollary \ref{cor:a=1}, Theorem \ref{thm:co-prime_almost}
and Lemma \ref{lem:data}.

\begin{corollary}\label{cor:easy_GE_exceptional}
If $f$ is of exceptional type,
then the assumption in Theorem \textup{\ref{thm:easy_GE}} holds
and consequently the general elephant theorem holds for $f$.
\end{corollary}

Therefore, with regard to the general elephant theorem,
it suffices to treat the case where $f$ is of general type
with discrepancy $a/n>1$.
We keep the following notation when $f$ is of general type.
\begin{enumerate}
\item
If $f$ is of No $\textrm{15}''$,
then $I$ consists of $Q_1=(r_1,1)$ and $Q_2=(r_2,1)$.
Set $b_1:=b_{Q_1}$ and $b_2:=b_{Q_2}$.
\item
If $f$ is of No 16,
then $I$ consists of $Q=(r,1)$.
Set $b:=b_Q$.
We set $r_1:=1$, $r_2:=r$ and $b_1:=0$, $b_2:=b$
to treat both the cases simultaneously.
\end{enumerate}
We have $(a/n)E^3=1/r_1+1/r_2$ by Theorem \ref{thm:classification}.
We let $S$ be a general elephant of $Y$ and
$S_X$ its birational transform on $X$.
The general elephant theorem guarantees that
$S_X$ has at worst a Du Val singularity at $P$
and the induced morphism $f_S \colon S \to S_X$ is crepant.
In our case we prove the theorem with
the description of this partial resolution.

\begin{theorem}\label{thm:hard_GE}
If $f$ is of general type,
then a general elephant $S$ of $Y$ has at worst Du Val singularities.
Moreover,
except for the case $a=1$,
the dual graph
\textup{(}see Section \textup{\ref{sec:small}} for notation\textup{)}
for the partial resolution $f_S \colon S \to S_X$ is one of the following.
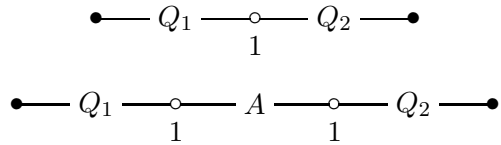
\begin{enumerate}
\item\label{itm:hard_A}
\begin{gather*}
\xymatrix@!0@=30pt{
*=<3.5pt>{\bullet} \ar@{-}[r] &
Q_1 \ar@{-}[r] &
*=<3.5pt>{\circ} \ar@{-}[r]_<*+{\>\!\!\!1} &
Q_2 \ar@{-}[r] &
*=<3.5pt>{\bullet}
}
\\
\xymatrix@!0@=30pt{
*=<3.5pt>{\bullet} \ar@{-}[r] &
Q_1 \ar@{-}[r] &
*=<3.5pt>{\circ} \ar@{-}[r]_<*+{\>\!\!\!1} &
A \ar@{-}[r] &
*=<3.5pt>{\circ} \ar@{-}[r]_<*+{\>\!\!\!1} &
Q_2 \ar@{-}[r] &
*=<3.5pt>{\bullet}
}
\end{gather*}
The attached number at $\circ$
is the coefficient in the fundamental cycle.
In No $\textrm{15}''$,
$Q_1$, $Q_2$ correspond to the points $Q_1$, $Q_2$ in $I$,
and in No 16,
$Q_1$, $Q_2$ correspond to a smooth point and $Q$ in $I$.
$S$ has Du Val singularities of types $A_{r_1-1}$, $A_{r_2-1}$, $A$
respectively at $Q_1$, $Q_2$, $A$
\textup{(}possibly smooth at $A$\textup{)},
and they are all the singularities of $S$.
In particular, $S_X$ has a Du Val singularity of type $A$
and consequently $P$ is of type c$A/n$.
\item
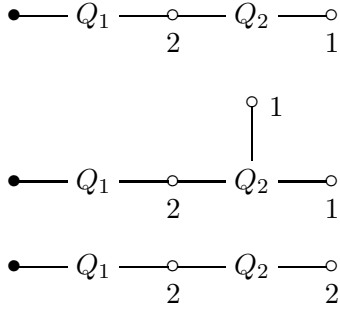
\label{itm:hard_D}
\begin{gather*}
\xymatrix@!0@=30pt{
*=<3.5pt>{\bullet} \ar@{-}[r] &
Q_1 \ar@{-}[r] &
*=<3.5pt>{\circ} \ar@{-}[r]_<*+{\>\!\!\!2} &
Q_2 \ar@{-}[r] \ar@{-}[r]_>*+{\ 1} &
*=<3.5pt>{\circ}
}
\\
\xymatrix@!0@=30pt{
*=<3.5pt>{} &
&
&
*=<3.5pt>{\circ} \ar@{-}[d]^<*+{1} &
\\
*=<3.5pt>{\bullet} \ar@{-}[r] &
Q_1 \ar@{-}[r] &
*=<3.5pt>{\circ} \ar@{-}[r]_<*+{\>\!\!\!2} &
Q_2 \ar@{-}[r] \ar@{-}[r]_>*+{\ 1} &
*=<3.5pt>{\circ}
}
\\
\xymatrix@!0@=30pt{
*=<3.5pt>{\bullet} \ar@{-}[r] &
Q_1 \ar@{-}[r] &
*=<3.5pt>{\circ} \ar@{-}[r]_<*+{\>\!\!\!2} &
Q_2 \ar@{-}[r] \ar@{-}[r]_>*+{\ 2} &
*=<3.5pt>{\circ}
}
\end{gather*}
In this case $f$ is of No $\textrm{15}''$ with $n=1$ or $2$,
and either
\begin{enumerate}
\item
$r_2=r_1+n$ with $a \mid 2r_1+n$ and $a$ is odd\textup{;} or
\item
$r_2=r_1+2n$ with $a \mid r_1+n$.
\end{enumerate}
In each case $a<r_1<r_2$,
$n \mid a+r_i$
and $b_i \equiv (a+r_i)/n$ modulo $r_i$ for $i=1,2$.
$Q_1$, $Q_2$ correspond to the points $Q_1$, $Q_2$ in $I$.
$S$ has Du Val singularities of types $A_{r_1-1}$, $A_{r_2-1}$
respectively at $Q_1$, $Q_2$,
and they are all the singularities of $S$
unless $S$ has the last graph.
In the case of the last graph,
$S$ may have further singularities
on the exceptional curve not passing through $Q_1$.
If $n=2$ then the assumption in Theorem \textup{\ref{thm:easy_GE}} holds.
\end{enumerate}
\end{theorem}

We let $f$ be a divisorial contraction of general type.
We start with the computation of $d(0,-1)$ and $d(-1,0)$.
By (\ref{eqn:difference_of_d(i,j)}),
we obtain that
\begin{align*}
d(1,-1)-d(0,-1) = & -\frac{3}{2}\Bigl(\frac{1}{r_1}+\frac{1}{r_2}\Bigr)
+\sum_{i=1,2}\Bigl(\frac{r_i-1}{2r_i}
-\frac{\overline{2}(r_i-\overline{2})}{2r_i}\Bigr), \\
d(0,0)-d(-1,0) = & -\Bigl(\frac{a}{n}+\frac{1}{2}\Bigr)
\Bigl(\frac{1}{r_1}+\frac{1}{r_2}\Bigr) \\
& + \sum_{i=1,2}\Bigl(\frac{b_i(r_i-b_i)}{2r_i}
-\frac{(b_i+1)(r_i-b_i-1)}{2r_i}\Bigr).
\end{align*}
We note that $d(0,0)=1$ by (\ref{eqn:d(i,j)_non-negative}).
Hence
\begin{align}
d(0,-1)&=d(1,-1)+
\begin{cases}
1 & \textrm{for No 15},\\
2 & \textrm{for No 16},
\end{cases}
\label{eqn:existence_H} \\
d(-1,0)&=2+\sum_{i=1,2}\Bigl(\frac{a}{n}-b_i\Bigr)\frac{1}{r_i}
\ge1.\label{eqn:existence_S}
\end{align}
By (\ref{eqn:d(i,j)=dim}) and (\ref{eqn:existence_H}),
we have $d(1,-1) \ge 0$ and
a general hyperplane section $H_X$ on the germ $P \in X$ has
multiplicity $1$ along $E$, that is,
$f^*H_X=H+E$, where $H$ is the birational transform on $Y$ of $H_X$.
By (\ref{eqn:d(i,j)=dim}) and (\ref{eqn:existence_S}),
there exists an effective divisor $S_X$ on $X$
which is linearly equivalent to $-K_X$ and has
multiplicity $a/n$ along $E$, that is,
$f^*S_X=S+(a/n)E$, where $S$ is the birational transform on $Y$ of $S_X$.
If we take $S_X$ generally
then the surface $S$ is a general elephant of $Y$.
Note that $S$ is not necessarily a prime divisor.
$S$ passes through all the singularities in $I_0$,
whereas $H$ passes through all those in $I$.
By Corollary \ref{cor:P1},
the curves $H \cap E$ and $S \cap E$ have no embedded points,
and their reduced loci are trees of $\bP^1$.
We reduce the case (\ref{itm:hard_A}) of Theorem \ref{thm:hard_GE}
to the following more approachable statement.

\begin{lemma}\label{lem:proper_intersect}
Theorem \textup{\ref{thm:hard_GE}(\ref{itm:hard_A})} follows from
one of the following.
\begin{enumerate}
\item\label{itm:proper_H}
$H \cap E$ has a \textup{(}possibly reducible\textup{)} component
which passes through all the points in $I$
and which intersects $S$ properly at two distinct points.
\item\label{itm:proper_S}
$S \cap E$ has a \textup{(}possibly reducible\textup{)} component
which intersects $H$ properly at two distinct points.
\end{enumerate}
\end{lemma}

\begin{proof}
The lemma is proved in the same fashion on either assumption.
Here we demonstrate it by assuming (\ref{itm:proper_H}).
We let $C$ denote the component of $H \cap E$ in the assumption.
Since $(S \cdot C) \le (S \cdot H \cdot E)=(a/n)E^3=1/r_1+1/r_2$,
the two points on $S \cap C$ are
$Q_1$ and $Q_2$ in No $\textrm{15}''$,
and a smooth point $Q_1$ and $Q_2:=Q$ in No 16.
Moreover $C=H \cap E$ and
$(S \cdot H \cdot E)_{Q_1}=1/r_1$, $(S \cdot H \cdot E)_{Q_2}=1/r_2$.
Hence $S$ is a Du Val singularity of type $A_{r_i-1}$ at $Q_i$ for $i=1,2$,
and $S \cap E$ is irreducible and reduced at each $Q_i$.
By the $f$-ampleness of $H$,
every irreducible component of $S \cap E$ passes through either $Q_1$ or $Q_2$.
Hence $S \cap E$ is isomorphic to $\bP^1$ or the tree of two $\bP^1$
with a node $Q'$ as pictured below.
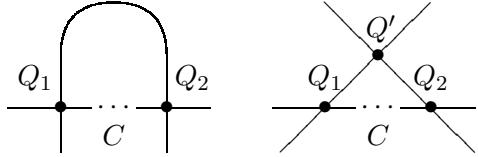
\begin{align*}
\xymatrix@!0@=20pt{
&
&
*=0{} \ar@{-}@/_7pt/[dl] \ar@{-}@/^7pt/[dr] &
&
&
&
*=0{} \ar@{-}[dr]^>>>*-{Q'} &
&
*=0{} \ar@{-}[dl] &
\\
&
*=0{} \ar@{-}[d] &
&
*=0{} \ar@{-}[d] &
&
&
&
*=0{\bullet} \ar@{-}[dr] \ar@{-}[dl] &
&
\\
*=0{} \ar@{-}[r]^>>>>*+{Q_1} &
*=0{\bullet} \ar@{-}[r]_>*+{\quad\;\; C} \ar@{-}[d] &
*=<17pt>{\cdots} \ar@{-}[r] &
*=0{\bullet} \ar@{-}[r]^<<<<*+{Q_2} \ar@{-}[d] &
&
*=0{} \ar@{-}[r]^>*+{Q_1} &
*=0{\bullet} \ar@{-}[r]_>*+{\quad\;\; C} \ar@{-}[dl] &
*=<17pt>{\cdots} \ar@{-}[r] &
*=0{\bullet} \ar@{-}[r]^<*+{Q_2} \ar@{-}[dr] &
\\
&&&&&&&&&
}
\end{align*}

Because $E$ is Cartier outside $Q_1$ and $Q_2$,
$S$ is smooth outside $Q_1$, $Q_2$ and the singularities of $S \cap E$.
The only possible singularity of $S \cap E$ is a node $Q'$,
at which $S$ has embedding dimension $2$ or $3$.
Suppose that $S$ is singular at $Q'$.
In order that $S \cap E$ has a node $Q'$,
the surface $S$ in the tangent space $T_{Q'}S=\bC^3$
has to be given by the equation $x_1x_2+x_3\cdot(\something)=0$
for suitable local coordinates $x_1,x_2,x_3$ at $Q'$,
where $E$ is given by $x_3=0$.
Hence $S$ has a Du Val singularity of type $A$ at $Q'$,
and the lemma is deduced.
\end{proof}

The rest of this section is devoted to the proof of Theorem \ref{thm:hard_GE}
on the assumption that $f$ has discrepancy $a/n>1$.
From now on,
we keep the assumption that
$f$ is of general type with discrepancy $a/n>1$.
The theorem is proved when $n=1$
by \cite[Section 5]{Ka03} and Erratum attached at the end of this paper.
We here assume that $n \ge 2$.
Because $(H \cdot H \cdot E)=E^3=(n/a)(1/r_1+1/r_2)<1/r_1+1/r_2$,
there exists an irreducible reduced component of $H \cap E$
which passes through all the singularities in $I$,
that is, $Q_1$ and $Q_2$ in No $\textrm{15}''$ and $Q$ in No 16.
We let $C \cong \bP^1$ be such an irreducible reduced component of $H \cap E$.
$Q_1, Q_2 \in C$ in No $\textrm{15}''$ and $Q \in C$ in No 16.
We take a normal form explained at the end of Section \ref{sec:numerical}
at each singularity in $I$.
The curve $C$ is determined uniquely in No $\textrm{15}''$
since $(H \cap E)_\red$ is a tree of $\bP^1$,
but in No 16 every irreducible component of $H \cap E$ passes through $Q$.
When $f$ is of No 16,
we choose $C$ wisely as in the following lemma.

\begin{lemma}\label{lem:D_16}
If $f$ is of No 16,
then there exists an effective divisor $D_X$ on $X$
whose multiplicity along $E$ is $1/n$, that is,
$f^*D_X=D+(1/n)E$,
where $D$ is the birational transform on $Y$ of $D_X$.
Moreover, by taking a suitable $H$,
we can choose $C$ so that $D$ intersects $C$ properly
and so that there exists no hidden non-Gorenstein point on $C$.
\end{lemma}

\begin{proof}
$a$ is co-prime to $n$ by Theorem \ref{thm:ramification}.
Thus there exist integers $s$, $t$ such that $as+nt=-1$.
The existence of $D_X$ follows if $d(s,t)\neq 0$ by (\ref{eqn:d(i,j)=dim}).
By (\ref{eqn:difference_of_d(i,j)}),
we obtain that
\begin{align*}
d(s+1,t)-d(s,t)=
&-\Bigl(\frac{1}{n}+\frac{1}{2}\Bigr)\Bigl(1+\frac{1}{r}\Bigr) \\
&+\Bigl(\frac{(u+1)(r-u-1)}{2r}-\frac{u(r-u)}{2r}\Bigr),
\end{align*}
where $u:=\overline{sb+t-1}$.
Since $d(s+1,t)=0$ by (\ref{eqn:d(i,j)_non-negative}),
we have
\begin{align*}
d(s,t)=\frac{1}{n}\Bigl(1+\frac{1}{r}\Bigr)+\frac{u+1}{r}>0,
\end{align*}
whence the existence of $D_X$ is deduced.

We have the surjective map
$f_*\cO_Y(H) \twoheadrightarrow H^0(\cQ_{0,-1})$ by Lemma \ref{lem:QandR}
and the estimate $h^0(\cQ_{0,-1})=d(0,-1)=2$ by
(\ref{eqn:H^k(Q)}) and (\ref{eqn:existence_H}).
Hence the restriction on $E$ of the divisor $H$ moves in its linear system,
although we need to mind the possibility that $E$ is non-normal.
Thus we can take $H$ so that
$H \cap E$ has an irreducible reduced component $C$
which intersects $D$ properly.
Suppose that there exists a hidden non-Gorenstein point $Q'$ on $C$.
Then the index $r'$ at $Q'$ is a divisor of $n$,
whence $D \sim sK_Y+tE$ is not Cartier at $Q'$ and
consequently $D$ passes through $Q'$.
However,
this contradicts that $((D+H) \cap E)_\red$ is a tree of $\bP^1$
in Corollary \ref{cor:P1} by $a/n>1$.
Therefore $Q$ is the only non-Gorenstein point which $C$ passes through.
\end{proof}

We consider $D_{i,j}$ in (\ref{eqn:definition_D})
and write $s_C(i,j)$ for the integer such that
\begin{align*}
(\cO_Y(D_{i,j}) \otimes \cO_C)/(\torsion) \cong \cO_{\bP^1}(s_C(i,j)).
\end{align*}
Note that $s_C(i,j)$ is well-defined for any $C$ on $Y$ isomorphic to $\bP^1$.
We give a condition for Theorem \textup{\ref{thm:hard_GE}}(\ref{itm:hard_A})
to hold in terms of the value of $s_C(-1,0)$.

\begin{lemma}\label{lem:numerical_condition}
Theorem \textup{\ref{thm:hard_GE}(\ref{itm:hard_A})} holds
if and only if $s_C(-1,0)=0$ in No $\textrm{15}''$ or $s_C(-1,0)=1$ in No 16.
In this case,
for any point $Q$ in $I$
\textup{(}$Q_1$ or $Q_2$ in No $\textrm{15}''$
and $Q$ in No 16\textup{)}
we may choose semi-invariant local coordinates
$x_1,x_2,x_3$ with weights $\wt(x_1, x_2, x_3)=(1, -1, b_Q)$
of the index-one cover $Q^\sharp \in Y^\sharp$ so that
\begin{align*}
Q^\sharp \in C^\sharp \subset Y^\sharp
\cong o \in (x_3\axis) \subset \bC^3_{x_1x_2x_3}.
\end{align*}
\end{lemma}

\begin{proof}
Let $R_0$ denote the global Gorenstein index of $Y$.
Compute the image of the map
\begin{align*}
\cO_Y(-K_Y)^{\otimes R_0} \otimes \cO_Y(R_0K_Y) \otimes \cO_C \to \cO_C.
\end{align*}
By (\ref{eqn:how.to.comp}),
its image at any $Q \in I$ is $\cO_C(-R_0w_{Q}^C(b_Q)/r_Q \cdot Q)$.
Hence the image is contained in the ideal sheaf isomorphic to
$\cO_{\bP^1}(-R_0(w_{Q_1}^C(b_1)/r_1+w_{Q_2}^C(b_2)/r_2))$
in No $\textrm{15}''$ and $\cO_{\bP^1}(-R_0w_Q^C(b)/r)$ in No 16.
On the other hand, this image equals that of
\begin{align*}
\bigl((\cO_Y(-K_Y)^{\otimes R_0} \otimes \cO_C)/(\torsion)\bigr)
\otimes \bigl(\cO_Y(R_0K_Y) \otimes \cO_C\bigr) \to \cO_C,
\end{align*}
whence
\begin{align*}
R_0s_C(-1,0)+R_0(K_Y \cdot C)\le
\begin{cases}
\displaystyle
-R_0\Bigl(\frac{w_{Q_1}^C(b_1)}{r_1}+\frac{w_{Q_2}^C(b_2)}{r_2}\Bigr) &
\textrm{for No $\textrm{15}''$}, \\
\displaystyle -R_0 \cdot \frac{w_Q^C(b)}{r} &
\textrm{for No 16}.
\end{cases}
\end{align*}
Note that by Lemma \ref{lem:D_16},
in fact the equality holds if $f$ is of No 16.
If Theorem \ref{thm:hard_GE}(\ref{itm:hard_A}) holds,
then $C=H \cap E$, $w_Q^C(b_Q)=1$ for any $Q \in I$
and the above inequality becomes an equality
by the proof of Lemma \ref{lem:proper_intersect}(\ref{itm:proper_H}),
whence the condition on $s_C(-1,0)$ holds.
On the other hand if the assumption on $s_C(-1,0)$ holds,
then the above inequality becomes
\begin{align*}
(-K_Y \cdot C)\ge
\begin{cases}
\displaystyle \frac{w_{Q_1}^C(b_1)}{r_1}+\frac{w_{Q_2}^C(b_2)}{r_2} &
\textrm{for No $\textrm{15}''$}, \\
\displaystyle 1+\frac{w_Q^C(b)}{r} &
\textrm{for No 16}.
\end{cases}
\end{align*}
We have $(-K_Y \cdot C) \le (-K_Y \cdot H \cdot E)=1/r_1+1/r_2$.
Hence we obtain that $C=H \cap E$, $w_Q^C(b_Q)=1$ for every $Q \in I$,
and there exists no hidden non-Gorenstein point on $C$.
The desired coordinates at $Q$ are constructed by $w_Q^C(b_Q)=1$.
By applying Lemma \ref{lem:QandR} to
$(i,j;i',j')=(0,-1;-1,0)$, $L_{0,-1}=H$ and $D_{-1,0}=-K_Y$,
we obtain the surjective map
\begin{align*}
f_*\cO_Y(-K_Y) \twoheadrightarrow H^0(\calR_{0,-1}^{-1,0})
&\twoheadrightarrow H^0((\cO_Y(-K_Y) \otimes \cO_C)/(\torsion)) \\
&=
\begin{cases}
H^0(\cO_{\bP^1})    & \textrm{for No $\textrm{15}''$}, \\
H^0(\cO_{\bP^1}(1)) & \textrm{for No 16}.
\end{cases}
\end{align*}
This surjective map and the equality $(S \cdot C)=1/r_1+1/r_2$
imply that $S$ intersects $C$ properly at
$Q_1$ and $Q_2$ in No $\textrm{15}''$,
and at $Q$ and another point in No 16.
Hence Theorem \ref{thm:hard_GE}(\ref{itm:hard_A}) holds
by Lemma \ref{lem:proper_intersect}(\ref{itm:proper_H}).
\end{proof}

Unlike the case where $X$ is Gorenstein,
$C$ might pass through hidden non-Gorenstein points in No $\textrm{15}''$.
In spite of this,
the proof of Lemma \ref{lem:numerical_condition} suggests that
there exists no hidden non-Gorenstein point on $C$
even if $f$ is of No $\textrm{15}''$.
We desire this property especially
in using the formula (\ref{eqn:how.to.comp}).
We prove that this is the case after a lemma
on the values of $s_C(i,j)$ which is used frequently.

\begin{lemma}\label{lem:s}
\begin{enumerate}
\item\label{itm:s_15}
If $f$ is of No $\textrm{15}''$, then
\begin{enumerate}
\item
$s_C(i,j)<0$ and $s_C(-i,-j) \ge -1$ if $ia/n+j>0$.
\item
$s_C(i,j)=-1$ and $s_C(-i,-j)=0$ or $-1$ if $0<ia/n+j<a/n$,
where $s_C(-i,-j)=0$ can occur only if $ib_1+j \equiv 0$ modulo $r_1$.
\item
$s_C(i,j)=-1$ if $ia/n+j=0$ with $(i,j)\notin\bZ(n,-a)$.
\item
$(s_C(1,0),s_C(-1,0))=(-2,0)$, $(-2,-1)$ or $(-1,-1)$.
\item
$s_C(1,1)=-1$.
\end{enumerate}
\item\label{itm:s_16}
If $f$ is of No 16, then
\begin{enumerate}
\item
$s_C(i,j)<0$ and $s_C(-i,-j) \ge 0$ if $ia/n+j>0$.
\item
$s_C(i,j)=-1$ and $s_C(-i,-j)=0$ if $0<ia/n+j<a/n$.
\item
$(s_C(1,0),s_C(-1,0))=(-2,1)$, $(-2,0)$ or $(-1,0)$.
\end{enumerate}
\end{enumerate}
\end{lemma}

\begin{proof}
Let $R_0$ be the global Gorenstein index of $Y$.
For (\ref{itm:s_15}),
we have $s_C(i,j)<0$ for $ia/n+j>0$ by the existence of the natural map
$\cO_Y(iK_Y+jE)^{\otimes R_0} \otimes \cO_C
\to \cO_Y(R_0(iK_Y+jE)) \otimes \cO_C$.
By the surjective map
$H^1(\calR_{0,-1}^{i,j}) \twoheadrightarrow H^1(\cO_{\bP^1}(s_C(i,j)))$
and Lemma \ref{lem:QandR},
we obtain that $s_C(i,j)\ge-1$ if $ia/n+j<a/n$ and that $s_C(1,0)=-1$ or $-2$.
We have (a)-(d) by considering that the map
$\cO_Y(iK_Y+jE) \otimes \cO_Y(-iK_Y-jE) \otimes \cO_C \to \cO_C$
is not surjective at $Q$ unless $iK_Y+jE$ is Cartier at $Q$.
Note that $iK_Y+jE$ is not Cartier at $Q_2$ for $0<ia/n+j\le a/n$
by $a<r_2$ in Lemma \ref{lem:minimal_discrepancy}.
By Lemma \ref{lem:QandR}, (\ref{eqn:d(i,j)_non-negative})
and (\ref{eqn:existence_H}),
we have $h^1(\cQ_{1,1})=0$, $h^2(\cQ_{1,1})=d(0,0)=1$
and $h^2(\cQ_{1,2})=d(0,-1)=1$. 
Hence $h^1(\calR_{0,-1}^{1,1})=0$ by (\ref{eqn:exact_R})
and (e) follows.

(\ref{itm:s_16}) is proved in the same way once we obtain
the isomorphism $(\cO_Y(D) \otimes \cO_C)/(\torsion) \cong \cO_{\bP^1}$
for $D=D_{s,t}$ consrtucted in Lemma \ref{lem:D_16}.
Since $(D \cdot C)=(1/n)(H \cdot C)=(1/a)(1+1/r)<1$,
we have $(D \cdot C)_Q=(D \cdot C)$.
Let $D^\sharp$ denote the pre-image of $D$
on the index-one cover $Q^\sharp \in Y^\sharp$.
Then the restriction of the defining function of $D^\sharp$,
whose weight is $-(sb+t)$,
to $C^\dag$ has order ${(D \cdot C)_Q}$ with respect to $t$,
whence $w_Q^C(-sb-t)=r(D \cdot C)_Q=r(D \cdot C)$.
From the map
$\cO_Y(D)^{\otimes r} \otimes \cO_Y(-rD) \otimes \cO_C \to \cO_C$,
we obtain that $rs_C(s,t)-r(D \cdot C)=-w_Q^C(-sb-t)$.
Hence $s_C(s,t)=0$.
\end{proof}

\begin{lemma}\label{lem:no_more_singularity}
There exists no hidden non-Gorenstein point on $C$.
\end{lemma}

\begin{proof}
By Lemma \ref{lem:D_16},
we have only to consider the case No $\textrm{15}''$.
We shall prove the lemma for No $\textrm{15}''$
by using the maps
\begin{align}\label{eqn:no_more_singularity}
\cO_Y(iK_Y+jE) \otimes \cO_Y(-iK_Y-jE) \otimes\cO_C \to \cO_C
\end{align}
for various $i$, $j$.
Note that the map (\ref{eqn:no_more_singularity})
is not surjective at every non-Gorenstein point on $C$
at which $iK_Y+jE$ is not Cartier.
In particular,
the number of those points is at most $-s_C(i,j)-s_C(-i,-j)$.
Suppose that there exists a hidden non-Gorenstein point $Q'$ on $C$.
We apply this principle with Lemma \ref{lem:s}(\ref{itm:s_15}).

For each $i=1,2$,
at least one of the maps (\ref{eqn:no_more_singularity})
for $(i,j)=(1,-1)$ and $(1,-2)$ is not surjective at $Q_i$,
and both of the maps are not surjective at $Q'$.
Since $-s_C(i,j)-s_C(-i,-j)\le2$ for these $(i,j)$
by Lemma \ref{lem:s}(\ref{itm:s_15}),
$Q'$ is the only hidden non-Gorenstein point on $C$,
and that $b_A=1$, $b_B=2$ for $(A,B)=(1,2)$ or $(2,1)$.
If $a/n \ge 2$,
then by the map (\ref{eqn:no_more_singularity}) for $(i,j)=(1,-3)$,
the divisor $K_Y-3E$ has to be Cartier at $Q_1$ or $Q_2$,
whence $(A,B)=(1,2)$ and $r_1=2$.
In particular,
$r_1$ is co-prime to $r_2$ and consequently
$a$ is co-prime to $n$ by Theorem \ref{thm:ramification}.
Thus $a/n>2$ and the map (\ref{eqn:no_more_singularity}) for $(i,j)=(1,-4)$
is not surjective at $Q_1$, $Q_2$ and $Q'$.
This contradicts that $-s_C(1,-4)-s_C(-1,4)=2$
in Lemma \ref{lem:s}(\ref{itm:s_15}).
Therefore $1<a/n<2$.
The map (\ref{eqn:no_more_singularity}) for $(i,j)=(2,-3)$ is not
surjective at $Q_1$, $Q_2$ by $b_A=1$, $b_B=2$.
Hence this map has to be surjective at $Q'$,
or equivalently the index $r_{Q'}$ of $Y$ at $Q'$ is $2$.
However also the map (\ref{eqn:no_more_singularity}) for $(i,j)=(3,-4)$
is not surjective at $Q_1$, $Q_2$, whence $r_{Q'}$ must be $3$.
This is a contradiction.
\end{proof}

We proceed to the delicate analysis of the embedding of $C$
by using the concept of normal form in Section \ref{sec:numerical}.
We need to keep Lemma \ref{lem:no_more_singularity} in mind.
First we prove that the embedding dimension of $C^\sharp$ is at most two.
This dimension equals the number of generators of
the semi-group $G_Q^C$ in (\ref{eqn:semigroup}).

\begin{lemma}\label{lem:emb_dim_2}
Consider any point $Q$ in $I$
\textup{(}$Q_1$ or $Q_2$ in No $\textrm{15}''$
and $Q$ in No 16\textup{)},
and take a normal form of $Q \in C \subset Y$.
Then the semi-group $G_Q^C$ is
generated either by $(a_1,1)$ and $(a_3,b_Q)$
or by $(a_2,-1)$ and $(a_3,b_Q)$.
In particular,
we may choose semi-invariant local coordinates
$x_1,x_2,x_3$ with weights $\wt(x_1, x_2, x_3)=(1, -1, b_Q)$
of the index-one cover $Q^\sharp \in Y^\sharp$ so that
$(x_1,x_2,x_3) |_{C^\dag}=(t^{c_Q/r_Q},0,t^{\overline{b_Qc_Q}/r_Q})$
or $(0,t^{c_Q/r_Q},t^{1-\overline{b_Qc_Q}/r_Q})$
for some positive integer $c_Q<r_Q$.
\end{lemma}

\begin{proof}
Because $C$ is smooth at $Q$,
there exists an invariant monomial in $x_1,x_2,x_3$
whose restriction to $C$ is $t$,
or equivalently $(r_Q,0) \in G_Q^C$.
Hence one of the following holds.
\begin{enumerate}
\item\label{itm:p1}
One of $a_1,a_2,a_3$ is $1$.
\item\label{itm:p2}
$(r_Q, 0)$ is generated either by $(a_1,1)$ and $(a_3,b_Q)$
or by $(a_2,-1)$ and $(a_3,b_Q)$ in $\bZ \times \bZ/(r_Q)$.
\item\label{itm:p3}
$r_Q=a_1+a_2$.
\end{enumerate}
The lemma follows directly if (\ref{itm:p1}) or (\ref{itm:p2}) holds.
From now on we suppose that only the statement (\ref{itm:p3}) holds.
In particular $r_Q \ge 5$.
On this assumption,
the inequality $w_Q^C(k)+w_Q^C(-k) \ge 2r_Q$ holds unless $k=0$ or $\pm1$.

If $a$ is not co-prime to $n$,
then $f$ is of No $\textrm{15}''$ by Theorem \ref{thm:ramification}
and we have a numerically $f$-trivial divisor $D_{n',-a'}$,
which is not Cartier at $Q_1$ and $Q_2$
by Lemma \ref{lem:d_Q}(\ref{itm:d_Q_ge}).
By the map
$\cO_Y(-D_{n',-a'}) \otimes \cO_Y(D_{n',-a'}) \otimes \cO_C \to \cO_C$,
we have
\begin{align*}
s_C(-n',a')+s_C(n',-a')=
-\sum_{i=1,2}
\Bigl(\frac{w_{Q_i}^C(n'b_i-a')}{r_i}+\frac{w_{Q_i}^C(a'-n'b_i)}{r_i}\Bigr).
\end{align*}
Since $s_C(n',-a')=s_C(-n',a')=-1$ by Lemma \ref{lem:s}(\ref{itm:s_15}),
we have $w_Q^C(n'b_Q-a')+w_Q^C(a'-n'b_Q)=r_Q$.
Hence (\ref{itm:p1}) or (\ref{itm:p2}) above must occur.

Now we may assume that $a$ is co-prime to $n$.
We take integers $s$, $t$ so that $as+nt=-1$
and consider the divisor $D_{s,t} \equiv -(1/n)E$.
Then for each $Q \in I$,
$r_Q$ is the smallest positive integer such that $r_QD_{s,t}$ is Cartier.
Consider the map
$\cO_Y(-iD_{s,t}) \otimes \cO_Y(iD_{s,t}) \otimes \cO_C \to \cO_C$.
We have
\begin{align*}
s_C(-is,-it)+s_C(is,it)=
-\sum_{Q \in I}
\Bigl(\frac{w_Q^C(i(sb_Q+t))}{r_Q}+\frac{w_Q^C(-i(sb_Q+t))}{r_Q}\Bigr).
\end{align*}
By this equality and Lemma \ref{lem:s},
\begin{align}\label{eqn:emb_dim_2}
w_Q^C(i(sb_Q+t))+w_Q^C(-i(sb_Q+t))=r_Q
\end{align}
for $0<i<a$ except that $f$ is of No $\textrm{15}''$ with $r_1 \mid i$.
On the other hand,
$w_Q^C(k)+w_Q^C(-k) \ge 2r_Q$ unless $k=0$ or $\pm1$ by the assumption.
First by (\ref{eqn:emb_dim_2}) for $i=1$,
we have $sb_Q+t \equiv \pm1$ modulo $r_Q$.
Then $2(sb_Q+t), 3(sb_Q+t) \not\equiv 0, \pm1$ modulo $r_Q$ by $r_Q\ge5$.
Thus by (\ref{eqn:emb_dim_2}) for $i=2$,
$f$ is of No $\textrm{15}''$, $Q=Q_2$ and $r_1=2$.
Hence by (\ref{eqn:emb_dim_2}) for $i=3$ and $a/n>1$,
we have $(a,n)=(3,2)$.
This contradicts that $n\equiv ae_{Q_1}$ modulo $r_1$.
\end{proof}

The precise generators of $G_Q^C$ are written as follows.

\begin{lemma}\label{lem:further}
\begin{enumerate}
\item\label{itm:further_15}
If $f$ is of No $\textrm{15}''$,
then we may choose the coordinates so that the restriction to $C^\dag$ at
one of $Q_1$ and $Q_2$ is of form $(t^{c_i/r_i},0,t^{\overline{b_ic_i}/r_i})$
and that at the other is of form $(0,t^{c_i/r_i},t^{1-\overline{b_ic_i}/r_i})$.
\item\label{itm:further_16}
If $f$ is of No 16,
then we may choose the coordinates so that the restriction to $C^\dag$ at $Q$
is of form $(t^{c/r},0,t^{\overline{bc}/r})$.
\end{enumerate}
\end{lemma}

\begin{proof}
We note that the lemma holds
on the assumption in Lemma \ref{lem:numerical_condition}.
Starting with Lemma \ref{lem:emb_dim_2},
first we consider the case No $\textrm{15}''$.
If both are of form  $(t^{c_i/r_i},0,t^{\overline{b_ic_i}/r_i})$,
then the image of the map
$\cO_Y(-K_Y) \otimes \cO_Y(-E) \otimes \cO_Y(K_Y+E) \otimes \cO_C \to \cO_C$
is $\cO_C(-Q_1-Q_2)$,
whence $s_C(-1,0)+s_C(0,-1)+s_C(1,1)=-2$.
Thus $s_C(-1,0)=0$ by Lemma \ref{lem:s}(\ref{itm:s_15}).
If both are of form $(0, t^{c_i/r_i}, t^{1-\overline{b_ic_i}/r_i})$,
then the image of the map
$\cO_Y(-K_Y) \otimes \cO_Y(E) \otimes \cO_Y(K_Y-E) \otimes \cO_C \to \cO_C$
is $\cO_C(-Q_1-Q_2)$,
whence $s_C(-1,0)+s_C(0,1)+s_C(1,-1)=-2$.
Thus $s_C(-1,0)=0$ by Lemma \ref{lem:s}(\ref{itm:s_15}).
If $f$ is of No 16 and
the restriction is of form $(0,t^{c/r},t^{1-\overline{bc}/r})$,
then the image of the map
$\cO_Y(-K_Y) \otimes \cO_Y(E) \otimes \cO_Y(K_Y-E) \otimes \cO_C \to \cO_C$
is $\cO_C(-Q)$,
whence $s_C(-1,0)+s_C(0,1)+s_C(1,-1)=-1$.
Thus $s_C(-1,0)=1$ by Lemma \ref{lem:s}(\ref{itm:s_16}).
\end{proof}

It is possible now to derive Theorem \ref{thm:hard_GE}
when $f$ is of No 16 with discrepancy $a/n>1$.
The divisor $D$ constructed in Lemma \ref{lem:D_16}
plays a crucial role in the proof.

\begin{theorem}\label{thm:hard_GE16}
If $f$ is of No 16 with discrepancy $a/n>1$,
then Theorem \textup{\ref{thm:hard_GE}} holds.
\end{theorem}

\begin{proof}
By virtue of Lemma \ref{lem:numerical_condition},
it suffices to prove that $s_C(-1,0)=1$.
First we shall deduce that $\overline{bc}<c$
in Lemma \ref{lem:further}(\ref{itm:further_16})
by using the effective divisor $D \sim sK_Y+tE \equiv -(1/n)E$
constructed in Lemma \ref{lem:D_16}.
From the map
$\cO_Y(D)^{\otimes r} \otimes \cO_Y(-rD) \otimes \cO_C \to \cO_C$,
we have $rs_C(s,t)-r(D \cdot C)=-w_Q^C(-sb-t)$,
whence
\begin{align*}
\frac{1}{n}(H \cdot C)=(D \cdot C)=\frac{w_Q^C(-sb-t)}{r}
\end{align*}
by $s_C(s,t)=0$ in Lemma \ref{lem:s}(\ref{itm:s_16}).
On the other hand, from the map
$\cO_Y(H)^{\otimes r} \otimes \cO_Y(-rH) \otimes \cO_C \to \cO_C$,
we have $rs_C(0,-1)-r(H \cdot C)=-w_Q^C(1)=-c$,
whence
\begin{align}\label{eqn:HC_16}
(H \cdot C)=\frac{c}{r}
\end{align}
by $s_C(0,-1)=0$ in Lemma \ref{lem:s}(\ref{itm:s_16}).
Hence $w_Q^C(-sb-t)=c/n<c$.
Because $G_Q^C$ in (\ref{eqn:semigroup}) is generated by
$(c,1)$ and $(\overline{bc},b)$,
this inequality implies that $\overline{bc}<c$.

From the map
$\cO_Y(-K_Y)^{\otimes r} \otimes \cO_Y(rK_Y) \otimes \cO_C \to \cO_C$,
we have $rs_C(-1,0)+r(K_Y \cdot C)=-w_Q^C(b)=-\overline{bc}$,
whence
\begin{align}\label{eqn:-KC_16}
\frac{a}{n}(H \cdot C)=(-K_Y \cdot C)=\frac{\overline{bc}}{r}+s_C(-1,0).
\end{align}
By (\ref{eqn:HC_16}), (\ref{eqn:-KC_16}) and the inequality $\overline{bc}<c$,
we have $s_C(-1,0)>0$.
Therefore the theorem follows from Lemma \ref{lem:s}(\ref{itm:s_16}).
\end{proof}

By virtue of Lemma \ref{lem:numerical_condition},
we assume that $f$ is of No $\textrm{15}''$ with $s_C(-1,0)\neq0$ from now on.
The goal of the rest of this section is to derive
Theorem \ref{thm:hard_GE}(\ref{itm:hard_D}) on this assumption.
The embedding dimension of $C^\sharp$ is one
in the situation of Lemma \ref{lem:numerical_condition},
but in fact this is always the case.

\begin{lemma}\label{lem:true/false}
Suppose that $f$ is of No $\textrm{15}''$ with $s_C(-1,0)\neq0$.
Then for each $i=1,2$,
we may choose semi-invariant local coordinates
$x_{i1},x_{i2},x_{i3}$ with weights $\wt(x_{i1},x_{i2},x_{i3})=(1,-1,b_i)$
of the index-one cover $Q_i^\sharp \in Y^\sharp$ so that
\begin{align*}
{Q_i}^\sharp \in C^\sharp \subset Y^\sharp
\cong o \in (x_{ii}\axis) \subset \bC^3_{x_{i1}x_{i2}x_{i3}}.
\end{align*}
Moreover $C \subset S$, $3\le r_1<r_2$
and $(H \cdot C)=1/r_1-1/r_2$, $(-K_Y \cdot C)=b_1/r_1-b_2/r_2$.
\end{lemma}

\begin{proof}
We note that $s_C(-1,0)=-1$ by Lemma \ref{lem:s}(\ref{itm:s_15}),
and the assumption $r_1\ge3$ is imposed
by the proof of Lemma \ref{lem:further}(\ref{itm:further_15}).
Moreover,
according to Lemmata \ref{lem:proper_intersect}(\ref{itm:proper_H})
and \ref{lem:numerical_condition},
if $S$ intersects $C$ properly then eventually $s_C(-1,0)=0$.
Hence in our case we have $C \subset S$.
By Lemma \ref{lem:further}(\ref{itm:further_15}),
we can choose $(A, B)=(1, 2)$ or $(2, 1)$ so that
the description at $Q_A$ is of form
$(t^{c_A/r_A},0,t^{\overline{b_Ac_A}/r_A})$
and that at $Q_B$ is of form $(0,t^{c_B/r_B},t^{1-\overline{b_Bc_B}/r_B})$.

\begin{step}\label{stp:true/false-1}
\itshape
$c_A < \overline{b_Ac_A}$.
\end{step}

Indeed,
by the map
$\cO_Y(H)^{\otimes r_Ar_B} \otimes \cO_Y(-r_Ar_BH) \otimes \cO_C \to \cO_C$,
we have $r_Ar_Bs_C(0,-1)-r_Ar_B(H \cdot C)=-r_Bw_{Q_A}^C(1)-r_Aw_{Q_B}^C(1)$.
By $w_{Q_A}^C(1)=c_A$, $w_{Q_B}^C(1)=r_B-c_B$ and
$s_C(0,-1)=-1$ in Lemma \ref{lem:s}(\ref{itm:s_15}),
we obtain that
\begin{align}\label{eqn:1+HC}
1+(H \cdot C)=\frac{c_A}{r_A}+\frac{r_B-c_B}{r_B}.
\end{align}
On the other hand,
by the map
$\cO_Y(-K_Y)^{\otimes r_Ar_B} \otimes \cO_Y(r_Ar_BK_Y) \otimes \cO_C
\to \cO_C$,
we have
$r_Ar_Bs_C(-1,0)+r_Ar_B(K_Y \cdot C)=-r_Bw_{Q_A}^C(b_A)-r_Aw_{Q_B}^C(b_B)$.
By $w_{Q_A}^C(b_A)=\overline{b_Ac_A}$, $w_{Q_B}^C(b_B)=r_B-\overline{b_Bc_B}$
and $s_C(-1,0)=-1$,
we obtain that
\begin{align}\label{eqn:1+a/nHC}
1+\frac{a}{n}(H \cdot C)=
\frac{\overline{b_Ac_A}}{r_A}+\frac{r_B-\overline{b_Bc_B}}{r_B}.
\end{align}
From the equalities (\ref{eqn:1+HC}) and (\ref{eqn:1+a/nHC}),
either $c_A < \overline{b_Ac_A}$ or $c_B > \overline{b_Bc_B}$ holds.
However, $c_B > \overline{b_Bc_B}$ never occurs since
$(r_B, 0)$ is contained in the semi-group $G_{Q_B}^C$ in (\ref{eqn:semigroup}).
Hence Step \ref{stp:true/false-1} is deduced.

\begin{step}\label{stp:true/false-2}
\itshape
$c_A=1$.
\end{step}

For Step \ref{stp:true/false-2},
it suffices to prove the equality $b_Ac_A=\overline{b_Ac_A}$
by $(r_A,0)\in G_{Q_A}^C$.
Suppose that $b_Ac_A>\overline{b_Ac_A}$.
Then $b_A \ge 3$ by $c_A<\overline{b_Ac_A}$.
Moreover,
$(\overline{b_Ac_A}-c_A, b_A-1)$ is not contained in $G_{Q_A}^C$,
and either $(\overline{b_Ac_A-2c_A}, b_A-2)$
or $(\overline{2c_A-b_Ac_A}, 2-b_A)$ is not contained in $G_{Q_A}^C$.
Therefore, for $j=1,2$,
the image of the map
$\cO_Y(-K_Y+jE) \otimes \cO_Y(K_Y-jE) \otimes \cO_C \to \cO_C$
is contained in $\cO_C(-2Q_A)$,
whence $s_C(-1,j)+s_C(1,-j)\le-3$ unless $K_Y-jE$ is Cartier at $Q_B$.
It is however impossible by $a/n>1$ and Lemma \ref{lem:s}(\ref{itm:s_15}).
Hence Step \ref{stp:true/false-2} follows.

Now by substituting $c_A=1$ in (\ref{eqn:1+HC}) and (\ref{eqn:1+a/nHC}),
we have $(A, B)=(1, 2)$, $r_1<r_2$, $(H \cdot C)=1/r_1-c_2/r_2$
and $(-K_Y \cdot C)=(a/n)(H \cdot C)=b_1/r_1-\overline{b_2c_2}/r_2$.
It remains proving that $c_2=1$.

\begin{step}\label{stp:true/false-3}
\itshape
$c_2=1$ except for the case where $b_1=r_1-1$ and $d(0,-2)\ge2$.
\end{step}

We have $1-(H \cdot C)=(r_1-1)/r_1+c_2/r_2$
and $1+(a/n)(H \cdot C)=b_1/r_1+(r_2-\overline{b_2c_2})/r_2$.
These equalities imply that
$(r_1-1)/r_1+c_2/r_2<b_1/r_1+(r_2-\overline{b_2c_2})/r_2$,
whence $c_2<r_2-\overline{b_2c_2}$.
Thus $c_2=1$ follows if
$(r_2-\overline{b_2c_2}-c_2,b_2+1)$ is contained in $G_{Q_2}^C$.
Consider the map
$\cO_Y(K_Y+E) \otimes \cO_Y(-K_Y-E) \otimes \cO_C \to \cO_C$.
Then $(r_2-\overline{b_2c_2}-c_2,b_2+1) \in G_{Q_2}^C$ holds
if the image of this map at $Q_2$ is $\cO_C(-Q_2)$.
Since $s_C(1,1)=-1$ and $s(-1,-1)\ge-1$ by Lemma \ref{lem:s}(\ref{itm:s_15}),
the image at $Q_2$ is $\cO_C(-Q_2)$
unless $K_Y+E$ is Cartier at $Q_1$,
or equivalently $b_1=r_1-1$.

Similarly,
the equality $c_2=1$ follows if
$(\overline{r_2-b_2c_2-2c_2},b_2+2)$ and
$(\overline{b_2c_2+2c_2},-b_2-2)$ are contained in $G_{Q_2}^C$,
and this is the case if the image of the map
$\cO_Y(K_Y+2E) \otimes \cO_Y(-K_Y-2E) \otimes \cO_C \to \cO_C$
at $Q_2$ is either $\cO_C$ or $\cO_C(-Q_2)$.
$s(-1,-2)\ge-1$ by Lemma \ref{lem:s}(\ref{itm:s_15}).
As in the proof of (e) in Lemma \ref{lem:s}(\ref{itm:s_15})
we can compute that $h^1(\calR_{0,-1}^{1,2})=d(0,-2)-1$,
whence $s_C(1,2)=-1$ if $d(0,-2)=1$.
Thus when $b_1=r_1-1$,
the image at $Q_2$ is $\cO_C$ or $\cO_C(-Q_2)$ unless $d(0,-2)\ge2$.
Therefore Step \ref{stp:true/false-3} is deduced.

We shall investigate the case left in Step \ref{stp:true/false-3}.

\begin{step}\label{stp:true/false-4}
\itshape
The assumption in Theorem \textup{\ref{thm:easy_GE}} holds
in the case left in Step \textup{\ref{stp:true/false-3}}.
\end{step}

Recall that $3\le r_1<r_2$.
By (\ref{eqn:difference_of_d(i,j)}) we have
\begin{align*}
d(1,-2)-d(0,-2)
= -\frac{5}{2}\Bigl(\frac{1}{r_1}+\frac{1}{r_2}\Bigr)
+\sum_{i=1,2}\Bigl(\frac{2(r_i-2)}{2r_i}-\frac{3(r_i-3)}{2r_i}\Bigr)
= -1,
\end{align*}
whence $d(1,-2)\ge1$ by the assumption that $d(0,-2)\ge2$.
In particular $1<a/n<2$ by (\ref{eqn:d(i,j)_non-negative}),
and we have only to prove that $d(-1,1)=0$.
Note that $b_2\ge3$ when $1<a/n<2$
by $nb_2 \equiv a$ modulo $r_2$
and $n<a<r_2$ in Lemma \ref{lem:minimal_discrepancy}.
By (\ref{eqn:difference_of_d(i,j)}) we have
\begin{align*}
d(0,1)-d(-1,1)
= & \Bigl(-\frac{a}{n}+\frac{1}{2}\Bigr)\Bigl(\frac{1}{r_1}+\frac{1}{r_2}\Bigr)
+\sum_{i=1,2}\bigl(B_{r_i}(-b_i+1)-B_{r_i}(-b_i)\bigr) \\
= & -1+\sum_{i=1,2}\Bigl(b_i-\frac{a}{n}\Bigr)\frac{1}{r_i}.
\end{align*}
By $b_1=r_1-1$ in the assumption
and $d(0,1)=0$ in (\ref{eqn:difference_of_d(i,j)}),
the above equality means that
\begin{align*}
d(-1,1)=\Bigl(\frac{a}{n}+1\Bigr)\frac{1}{r_1}+
\Bigl(\frac{a}{n}-b_2\Bigr)\frac{1}{r_2}.
\end{align*}
The first term $(a/n+1)/r_1$ is $<1$ by $a/n<2$ and $r_1\ge3$,
and the second term $(a/n-b_2)/r_2$ is $<0$ by $a/n<2$ and $b_2\ge3$.
Therefore $d(-1,1)=0$ and Step \ref{stp:true/false-4} follows.

In the case left in Step \ref{stp:true/false-3},
the surface $S$ has at worst Du Val singularities
by Step \ref{stp:true/false-4}.
In particular,
as in the argument before Table \ref{tbl:reproduce},
its pre-image $S^\sharp$ on the index-one cover $Q_2^\sharp \in Y^\sharp$
also has a Du Val singularity at $Q_2^\sharp$.
Thus $S^\sharp$ has multiplicity $\le 2$ at $Q_2^\sharp$,
or equivalently,
the defining function $h_S$ of $S^\sharp$ in $x_{21},x_{22},x_{23}$,
whose weight is $b_2$,
has order $\le2$.
On the other hand, from $C \subset S$,
the restriction of $h_S$ to $C^\dag$ by
$(x_{21},x_{22},x_{23})|_{C^\dag}=(0,t^{c_2/r_2},t^{1-\overline{b_2c_2}/r_2})$
is identically zero.
By $c_2<r_2-\overline{b_2c_2}$,
such a function $h_S$ exists only if $c_2=1$ or $b_2=1,2$.
By $b_2\ge3$ as seen in the proof of Step \ref{stp:true/false-4},
we obtain that $c_2=1$ also in the case left in Step \ref{stp:true/false-3}.
\setcounter{step}{0}
\end{proof}

In general it is difficult to control the denominator $n$ of the discrepancy,
but on our assumption we can narrow the possibilities down to one value.

\begin{lemma}\label{lem:n=2}
Suppose that $f$ is of No $\textrm{15}''$ with $s_C(-1,0)\neq0$.
Then $n=2$, $a<r_1<r_2$ and $b_i=(a+r_i)/2$ for $i=1,2$.
\end{lemma}

\begin{proof}
Recall the notation that $g=\gcd(a,n)$ and $a=ga'$, $n=gn'$.
We claim that $g\le2$.
Suppose that $g \neq 1$ and consider $D_{n',-a'}=n'K_Y-a'E \equiv 0$.
By the map
$\cO_Y(-D_{n',-a'})^{\otimes r_1r_2} \otimes \cO_Y(r_1r_2D_{n',-a'})
\otimes \cO_C \to \cO_C$,
we have
$s_C(-n',a')=-\sum_{i=1,2}w_{Q_i}^C(n'b_i-a')/r_i$.
On the other hand,
$s_C(-n',a')=-1$ by Lemma \ref{lem:s}(\ref{itm:s_15}),
and $w_{Q_1}^C(n'b_1-a')=\overline{d_{Q_1}^gb_1}$ and
$w_{Q_2}^C(n'b_2-a')=r_2-\overline{d_{Q_2}^gb_2}$
by Lemma \ref{lem:true/false}.
Therefore $\overline{d_{Q_1}^gb_1}/r_1=\overline{d_{Q_2}^gb_2}/r_2$.
This equality and Lemma \ref{lem:d_Q}(\ref{itm:d_Q_ge}) imply that
$\overline{d_{Q_1}^gb_1}/r_1=\overline{d_{Q_2}^gb_2}/r_2=1/2$,
whence $g=2$.

We have $(a/n)(H \cdot C)=(-K_Y \cdot C)$.
Thus Lemma \ref{lem:true/false} implies that
\begin{align}\label{eqn:keyeqn}
\Bigl(\frac{a}{n}-b_1\Bigr)\frac{1}{r_1}=
\Bigl(\frac{a}{n}-b_2\Bigr)\frac{1}{r_2}.
\end{align}
By this equality and (\ref{eqn:existence_S}),
we have $(a/n-b_i)/r_i \in \bZ/2$.
In particular $a'/n'-b_i \in \bZ \cdot (g/2)$
by $g \mid r_i$ in Theorem \ref{thm:ramification}.
Therefore $n' \mid 2$ if $g=1$, and $n'=1$ if $g=2$.
Anyway $n=2$ holds.

We then claim that $b_2=(a+r_2)/2$,
which completes the lemma by (\ref{eqn:keyeqn}).
Recall that $a<r_2$ by Lemma \ref{lem:minimal_discrepancy}.
If $a$ is odd,
then the relation $2b_2 \equiv a$ modulo $r_2$ implies that $2b_2=a+r_2$.
If $a$ is even,
then $r_i$ for $i=1,2$ is also even by Theorem \ref{thm:ramification}
and thus $b_i$ is odd.
We have $d_{Q_i}^2=r_i/2$ by Lemma \ref{lem:d_Q}(\ref{itm:d_Q_ge}).
In particular $1-e_{Q_2}a/2 \equiv r_2/2$ modulo $r_2$,
whence $b_2=(a+r_2)/2$ by $v_{Q_2}=1$ and $a<r_2$.
\end{proof}

The general elephant theorem for $f$ is derived as a corollary.

\begin{corollary}\label{cor:easy_GE_D}
If $f$ is of No $\textrm{15}''$ with $s_C(-1,0)\neq0$,
then the assumption in Theorem \textup{\ref{thm:easy_GE}} holds.
\end{corollary}

\begin{proof}
Note that $b_i=(a+r_i)/2$ for $i=1,2$ in Lemma \ref{lem:n=2}
and thus $d(-1,0)=1$ by (\ref{eqn:existence_S}).
For $0 \le j \le a/2$, by (\ref{eqn:difference_of_d(i,j)}) we have
\begin{align*}
d(0,j)-d(-1,j)=\Bigl(j-\frac{a+1}{2}\Bigr)
\Bigl(\frac{1}{r_1}+\frac{1}{r_2}\Bigr)+B_{-1,j}-B_{-1,j-1}.
\end{align*}
By (\ref{eqn:d(i,j)_non-negative}),
$d(0,0)=1$ and $d(0,j)=0$ for $1\le j\le a/2$. Hence
\begin{align*}
d(-1,j)-d(-1,j+1) =& \Bigl(\frac{1}{r_1}+\frac{1}{r_2}\Bigr)
+B_{-1,j+1}-2B_{-1,j}+B_{-1,j-1} \\
&+
\begin{cases}
1 & \textrm{for} \ j=0, \\
0 & \textrm{for} \ 1\le j\le a/2-1.
\end{cases}
\end{align*}
By the notation $u_i^j:=b_i-j>0$,
the term $B_{-1,j+1}-2B_{-1,j}+B_{-1,j-1}$ equals
\begin{align*}
& \sum_{i=1,2} \Bigl(\frac{(u_i^j-1)(r_i-u_i^j+1)}{2r_i}
-2\frac{u_i^j(r_i-u_i^j)}{2r_i}
+\frac{(u_i^j+1)(r_i-u_i^j-1)}{2r_i}\Bigr) \\
&= -\sum_{i=1,2}\frac{1}{r_i}.
\end{align*}
Therefore we obtain that $d(-1,j)=0$ for $1 \le j \le a/2$.
\end{proof}

The following geometric statement is important
in finding a component of $S \cap E$ intersecting $H$ properly.

\begin{lemma}\label{lem:irreducible}
Suppose that $f$ is of No $\textrm{15}''$ with $s_C(-1,0)\neq0$.
Then $H \cap E$ is irreducible but non-reduced.
\end{lemma}

\begin{proof}
We write $[H \cap E]=s[C]+[F]$ cycle-theoretically,
and set $(n/a)(1/r_1+1/r_2)=(H \cdot [H \cap E])=s(H \cdot C)+t/r_1+u/r_2$.
By $(H \cdot C)=1/r_1-1/r_2$ in Lemma \ref{lem:true/false},
we obtain that
\begin{align}\label{eqn:irreducible_eqn}
\frac{s+t-n/a}{r_1}=\frac{s+n/a-u}{r_2}.
\end{align}
By $r_1<r_2$ in Lemma \ref{lem:true/false},
we have $(t,u)=(0,0)$, $(1,0)$ or $(0,1)$.
By Corollary \ref{cor:P1},
if $[F]\neq0$ then $F\cong\bP^1$ and
$F$ intersects $C$ exactly at one of $Q_1$ and $Q_2$,
and $(H \cdot F)=1/r_1$ or respectively $1/r_2$.
In this case we have $a/n<2$ by (\ref{eqn:irreducible_eqn}),
whence $(a,n)=(3,2)$ by Lemma \ref{lem:n=2}.
Let $Q_i$ where $i=1$ or $2$ be the one which $F$ passes through.
Then $(H \cdot F)=1/r_i$ and $(S \cdot F)=3/2r_i$.
Since any hidden non-Gorenstein point has index $2$,
the estimate $(S \cdot F)=3/2r_i$ implies that $r_i$ is an odd integer.

The equality (\ref{eqn:irreducible_eqn}) becomes
\begin{align*}
(3s+2)r_1&=(3s+1)r_2 \quad \textrm{if} \ (t,u)=(1,0), \\
(3s-1)r_1&=(3s-2)r_2 \quad \textrm{if} \ (t,u)=(0,1),
\end{align*}
whence we can set
\begin{align*}
r_1&=(3s+1)k,\quad r_2=(3s+2)k, \quad \textrm{if} \ (t,u)=(1,0), \\
r_1&=(3s-2)k,\quad r_2=(3s-1)k, \quad \textrm{if} \ (t,u)=(0,1),
\end{align*}
for some integer $k$.
Here $s$ is even and $k$ is odd by $2 \nmid r_i$,
whence $r_1+r_2$ is odd.
On the other hand,
$(a/n)^2r_1r_2E^3=3(r_1+r_2)/2$
has to be an integer by Lemma \ref{lem:integer},
which is a contradiction.
Hence $(H \cap E)_\red$ is irreducible.

Consider the map
$\cO_Y(-K_Y) \otimes \cO_Y(K_Y) \otimes \cO_C \to \cO_C$.
By the description in Lemma \ref{lem:true/false},
we have $w_{Q_i}^C(b_i)+w_{Q_i}^C(-b_i)=r_i$ for $i=1,2$.
Thus $s_C(-1,0)+s_C(1,0)=-2$, whence $s_C(1,0)=-1$
by the assumption and Lemma \ref{lem:s}(\ref{itm:s_15}).
If $H \cap E$ is reduced, or equivalently $C=H \cap E$,
then $s_C(1,0)=-2$ by $h^1(\calR_{0,-1}^{1,0})=1$ in Lemma \ref{lem:QandR}.
Therefore $H \cap E$ is non-reduced.
\end{proof}

\begin{corollary}\label{cor:nmid}
$r_1 \nmid 2r_2$.
\end{corollary}

\begin{proof}
Set $[H \cap E]=s[C]$ with $s\ge2$ as in Lemma \ref{lem:irreducible}.
Then $1/r_1+1/r_2=(S \cdot [H \cap E])=s(S \cdot C)=(sa/n)(1/r_1-1/r_2)$
by Lemma \ref{lem:true/false},
whence $r_2/r_1=(as+2)/(as-2)$.
On the condition that $a>n=2$ and $s\ge2$,
the value of $r_2/r_1$ is contained in $\bZ/2$
only if $(a,s,r_2/r_1)=(3,2,2)$ or $(5,2,3/2)$.
In each case,
either $r_1$ or $r_2$ is even,
which contradicts that $b_i=(a+r_i)/2$ in Lemma \ref{lem:n=2}.
Hence the corollary is deduced.
\end{proof}

We recall that the strategy for the case No 16 is
to choose $C$ which intersects $D$ properly.
In No $\textrm{15}''$
we instead prove that $S \cap E$ has an irreducible component
which intersects $H$ properly.

\begin{lemma}\label{lem:move_S}
Suppose that $f$ is of No $\textrm{15}''$ with $s_C(-1,0)\neq0$.
Then $S \cap E$ has an irreducible component which intersects $H$ properly.
\end{lemma}

\begin{proof}
By Corollary \ref{cor:easy_GE_D},
the birational transform $S$ on $Y$ of a general elephant $S_X$ of $X$
is also a general elephant of $Y$ and has at worst Du Val singularities.
We examine the dual graph for the partial resolution $f_S \colon S \to S_X$.
$(H \cap E)_\red$ is irreducible by Lemma \ref{lem:irreducible}.
Suppose that $C=(H \cap E)_\red=(S \cap E)_\red$.
Then $C$ is the only exceptional curve of $f_S$.
By Lemma \ref{lem:irreducible} and $a/n>1$,
the cycle $[S \cap E]$ has the component $[C]$ with coefficient $\ge3$.
Thus the coefficient of $C$ in the fundamental cycle,
given by that of $C$ in $(f^*H_X)|_S= H|_S+E|_S$, is $\ge4$.
Hence the dual graph of the exceptional curve of $f_S$ is as follows.
\begin{align*}
\xymatrix@!0@=30pt{
&
*=<3.5pt>{\circ} \ar@{}[r]_<*+{\!\!\ge4} &
}
\end{align*}
By Table \ref{tbl:Mo85} and
the list of dual graphs from Du Val singularities
at the beginning of Section \ref{sec:small},
$P$ has to be of type c$E/2$ and $C$ is the curve
labelled $F_3$ in the dual graph for type $E_7$ in Section \ref{sec:small}.
Moreover, $[S \cap E]=3[C]$ and thus
$-(C^2)_S=-(E \cdot E \cdot S)/9=(1/r_1+1/r_2)/9$.
On the other hand, $-(F_3)^2=1/12$ by the dual graph.
Therefore $(1/r_1+1/r_2)/9=1/12$,
which contradicts that $3\le r_1<r_2$ in Lemma \ref{lem:true/false}.
\end{proof}

Lemma \ref{lem:move_S} means that
$S \cap E$ has an irreducible component other than $C$.
The next lemma decides $S \cap E$ cycle-theoretically.

\begin{lemma}\label{lem:cycle_S}
$[S \cap E]=[C]+[V]$ cycle-theoretically,
where $[V]$ consists of one or two $[\bP^1]$.
$V$ intersects $C$ at and only at $Q_2$.
\end{lemma}

\begin{proof}
Write $[S \cap E]=x[C]+[V]$ cycle-theoretically
and set $1/r_1+1/r_2=(H \cdot [S \cap E])=x(H \cdot C)+y/r_1+z/r_2$.
By Lemma \ref{lem:true/false}, we obtain that
\begin{align*}
\frac{x-1+y}{r_1}=\frac{x+1-z}{r_2}.
\end{align*}
We have $x \ge 1$ by $C \subset S$ in Lemma \ref{lem:true/false}
and $(y,z)\neq(0,0)$ by Lemma \ref{lem:move_S}.
Hence by $r_1<r_2$ in Lemma \ref{lem:true/false},
we have $(y,z)=(1,0)$ or $(0,1)$, or $(x,y,z)=(1,0,2)$.
If $(x,y,z)=(1,0,2)$ then $[S \cap E]=[C]+[V]$ with $(H \cdot V)=2/r_2$
and the lemma holds by Corollary \ref{cor:nmid}.
From now on we suppose that $(y,z)=(1,0)$ or $(0,1)$
and derive a contradiction.

By Corollary \ref{cor:nmid},
if $(y,z)=(1,0)$ or $(0,1)$ then
$V \cong \bP^1$ and $V$ intersects $C$ exactly at one of $Q_1$ and $Q_2$,
and $(H \cdot V)=1/r_1$ or respectively $1/r_2$.
Let $Q_i$ where $i=1$ or $2$ be the one which $V$ passes through.
Then $(H \cdot V)_{Q_i}=(H \cdot V)=1/r_i$.
Since $E$ is Cartier outside $Q_1,Q_2$,
the isomorphism $V \cong \bP^1$ implies that
$S$ is smooth along $V$ outside $Q_i$,
whence there exists no hidden non-Gorenstein point on $V$.
Thus by the map
$\cO_Y(S)^{\otimes r_i} \otimes \cO_Y(-r_iS) \otimes \cO_V \to \cO_V$,
we have $r_is_V(-1,0)-r_i(S \cdot V)=-w_{Q_i}^V(b_i)$.
Since $(S \cdot V)=(a/2)(H \cdot V)=a/2r_i$,
we obtain that $w_{Q_i}^V(b_i)\equiv a/2$ modulo $r_i$.
On the other hand $(H \cdot V)_{Q_i}=1/r_i$ implies that $w_{Q_i}^V(1)=1$,
whence $w_{Q_i}^V(b_i)=b_i$.
This contradicts that $b_i=(a+r_i)/2$ in Lemma \ref{lem:n=2}.
\end{proof}

It is the time to complete Theorem \ref{thm:hard_GE}
on the assumption that $a/n>1$.

\begin{theorem}\label{thm:hard_GE15}
If $f$ is of No $\textrm{15}''$ with discrepancy $a/n>1$,
then Theorem \textup{\ref{thm:hard_GE}} holds.
\end{theorem}

\begin{proof}
Theorem \ref{thm:hard_GE}(\ref{itm:hard_A}) holds when $s_C(-1,0)=0$
thanks to Lemma \ref{lem:numerical_condition}.
Here we consider the case $s_C(-1,0)\neq0$
and derive Theorem \ref{thm:hard_GE}(\ref{itm:hard_D}).
First we construct the coordinates in Lemma \ref{lem:true/false}
using the surfaces $H$, $E$ and $S$.

\begin{step}\label{stp:hard_GE15-1}
\itshape
We may choose the coordinates in Lemma \textup{\ref{lem:true/false}}
so that the pre-images $E^\sharp$, $S^\sharp$ on $Q_1^\sharp \in Y^\sharp$
of $E$, $S$ are given by respectively $x_{12}=0$, $x_{13}=0$,
and the pre-images $H^\sharp$, $S^\sharp$ on $Q_2^\sharp \in Y^\sharp$
of $H$, $S$ are given by respectively $x_{21}=0$, $x_{23}=0$.
\end{step}

The curve $C$ is the intersection of $S$ and $E$ near $Q_1$
by Lemma \ref{lem:cycle_S},
whence we obtain the desired coordinates $x_{11},x_{12},x_{13}$
by the description of $Q_1\in C \subset Y$ in Lemma \ref{lem:true/false}.
We consider the coordinates $x_{21},x_{22},x_{23}$.
We have $(H \cdot V)_{Q_2}=(H \cdot V)=2/r_2$ by Lemma \ref{lem:cycle_S},
whence $H^\sharp$ has multiplicity $\le 2$ at $Q_2^\sharp$.
If $H^\sharp$ is not smooth at $Q_2^\sharp$,
then $\wt x_{2i} +\wt x_{2j} \equiv 1$ modulo $r_2$
for some $i,j=1$, $2$ or $3$ with $(i,j)\neq(2,2)$ by $C \subset H$.
This implies that either $b_2\equiv2$ or $2b_2\equiv1$ modulo $r_2$.
However,
this contradicts that
$2\le a<r_1<r_2$ and $b_2=(a+r_2)/2$ in Lemma \ref{lem:n=2}.
Hence $H^\sharp$ is smooth at $Q_2^\sharp$,
and we can choose $x_{21}$ so that
$H^\sharp$ at $Q_2^\sharp$ is given by $x_{21}=0$.
On the other hand,
as in the argument before Table \ref{tbl:reproduce},
$S^\sharp$ on $Q_2^\sharp \in Y^\sharp$
has a Du Val singularity at $Q_2^\sharp$,
whence $S^\sharp$ has multiplicity $\le 2$ at $Q_2^\sharp$.
In fact $S^\sharp$ is smooth at $Q_2^\sharp$
by the same argument as for $H^\sharp$.
Therefore in order to complete Step \ref{stp:hard_GE15-1},
we have only to prove that $\wt x_{21}\not\equiv \wt x_{23}$ modulo $r_2$,
but this also follows from the equality $b_2=(a+r_2)/2$ in Lemma \ref{lem:n=2}.

We have decided the cycle $[S \cap E]$ in Lemma \ref{lem:cycle_S}.
In order to obtain the dual graph for
the partial resolution $f_S \colon S \to S_X$,
we must study the cycle $[H \cap S]$.

\begin{step}\label{stp:hard_GE15-2}
\itshape
$[H \cap S]=[C]+[B]$ cycle-theoretically,
where $B$ is the birational transform of the curve $H_X \cap S_X$.
$B$ intersects $E$ at and only at $Q_1$.
\end{step}

Let $c_{HS}$, $c_{HE}$, $c_{SE}$ denote the coefficients of $[C]$
in the cycles $[H \cap S]$, $[H \cap E]$, $[S \cap E]$.
Consider the relation $(H \cap S) \cap (H \cap E) \subset S \cap E$.
$H$ is smooth at the generic point of $C$
according to Step \ref{stp:hard_GE15-1}.
Thus we obtain that $\min\{c_{HS},c_{HE}\} \le c_{SE}$.
This implies that $c_{HS}=1$ since
$c_{HE}>1$ by Lemma \ref{lem:irreducible} and
$c_{SE}=1$ by Lemma \ref{lem:cycle_S}.
Hence the equality $[H \cap S]=[C]+[B]$ holds.
Then $(E \cdot B)=(E \cdot [H \cap S])-(E \cdot C)=2/r_1$,
whence $B$ intersects $E$ at and only at $Q_1$ by Corollary \ref{cor:nmid}.

The coefficient of $C$ in the fundamental cycle is
given by that of $C$ in $(f^*H_X)|_S= H|_S+E|_S$,
and this is two by Lemma \ref{lem:cycle_S} and Step \ref{stp:hard_GE15-2}.
By Lemma \ref{lem:cycle_S} and
Steps \ref{stp:hard_GE15-1}, \ref{stp:hard_GE15-2},
we obtain the dual graph for $f_S$ described
in Theorem \ref{thm:hard_GE}(\ref{itm:hard_D}).
Thanks to Lemma \ref{lem:n=2} and Corollary \ref{cor:easy_GE_D},
it suffices to prove that
either (a) or (b) in Theorem \ref{thm:hard_GE}(\ref{itm:hard_D}) holds.

\begin{step}\label{stp:hard_GE15-3}
\itshape
Either \textup{(a)} or \textup{(b)}
in Theorem \textup{\ref{thm:hard_GE}(\ref{itm:hard_D})} holds.
\end{step}

$r_2-r_1=2(b_2-b_1)$ by Lemma \ref{lem:n=2}.
Set $r_2-r_1=2k$.
Note that if $a$ is even then each $b_i$ is odd and thus $k$ is even.
We also set $[H \cap E]=s[C]$ as in Lemma \ref{lem:move_S}.
Then $1/r_1+1/r_2=(S \cdot [H \cap E])=s(S \cdot C)=(sa/2)(1/r_1-1/r_2)$
by Lemma \ref{lem:true/false},
or equivalently $r_1+r_2=sak$.
Thus
\begin{align}\label{eqn:hard_GE15-3}
r_1=(as/2-1)k,\qquad r_2=(as/2+1)k.
\end{align}
On the other hand,
by using the coordinates in Step \ref{stp:15/16small-2},
we can express $H^\sharp$ on $Q_1^\sharp \in Y^\sharp$ as
$x_{12}\cdot(\something)+x_{13}^s=0$,
and $E^\sharp$ on $Q_2^\sharp \in Y^\sharp$ as
$x_{21}\cdot(\something)+x_{23}^s=0$.
Considering weights,
we obtain that $b_1s\equiv 1$ modulo $r_1$ and $b_2s\equiv -1$ modulo $r_2$.
By $b_i=(a+r_i)/2$,
this implies that $as\equiv 2$ modulo $r_1$ and $as\equiv -2$ modulo $r_2$.
Therefore $(as,k)=(2r_1+2,1)$ or $(r_1+2,2)$ by (\ref{eqn:hard_GE15-3}).
These cases correspond to respectively (a), (b)
in Theorem \ref{thm:hard_GE}(\ref{itm:hard_D}).
This completes Step \ref{stp:hard_GE15-3},
and thus the theorem.
\setcounter{step}{0}
\end{proof}

\section{Case of small discrepancies}\label{sec:small}
In this section we study the cases left in Sections
\ref{sec:co-primeness} and \ref{sec:GE}
by using the general elephant theorem in the strongest form.
In the former half,
we proceed with the discussion on the co-primeness of $a$ and $n$
after Section \ref{sec:co-primeness}
and complete Theorem \ref{thm:co-prime}.
We also construct examples of divisorial contractions
in the cases in the first, second, third and last lines of Table \ref{tbl:MT}.
In the latter half,
we prove Theorem \ref{thm:hard_GE}
when $f$ is of general type with discrepancy $a/n \le 1$.

According to Theorem \ref{thm:easy_GE} and
Corollary \ref{cor:easy_GE_exceptional},
we can start from taking a general elephant of $S$
as the birational transform of a general elephant $S_X$ of $X$.
It is reasonable to examine the crepant morphism $f_S \colon S \to S_X$
in terms of the dual graph as in Theorem \ref{thm:hard_GE}.
Hence we here recall the dual graph for
the minimal resolution of a Du Val singularity.
$\circ$ denotes an exceptional curve and
$\bullet$ the birational transform of a general hyperplane section,
and each exceptional curve $F_i$ is marked with
its coefficient in the fundamental cycle.
Note that $S_X$ has a Du Val singularity of type at worst $E_7$ at $P$
if $X$ is not Gorenstein by Table \ref{tbl:reproduce}.
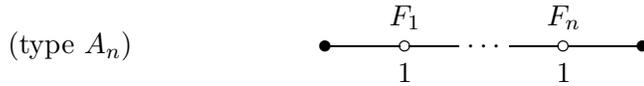
\begin{align}
\xymatrix@!0@=30pt{
*=<3.5pt>{\bullet} \ar@{-}[r] &
*=<3.5pt>{\circ} \ar@{-}[r]^<*+{\>\!\!\!F_1}_<*+{\>\!\!\!1} &
\cdots \ar@{-}[r] & 
*=<3.5pt>{\circ} \ar@{-}[r]^<*+{\>\!\!\!F_n}_<*+{\>\!\!\!1} &
*=<3.5pt>{\bullet}
}\tag{type $A_n$}
\end{align}
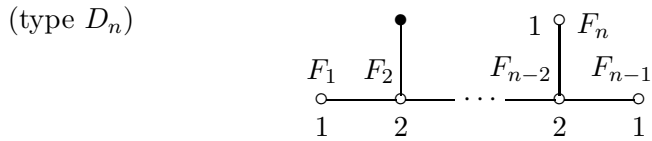
\begin{align}
\xymatrix@!0@=30pt{
&
*=<3.5pt>{\bullet} \ar@{-}[d] &
&
*=<3.5pt>{\circ} \ar@{-}[d]^<*+{F_n}_<*+{1} &
\\
*=<3.5pt>{\circ}
\ar@{-}[r]^<*+{\>\!\!\!F_1}_<*+{\>\!\!\!1}^>>>*+{F_2}_>*+{\ 2} &
*=<3.5pt>{\circ} \ar@{-}[r] &
\cdots \ar@{-}[r]^>>>>>*+{F_{n-2}}_>*+{\ 2} & 
*=<3.5pt>{\circ} \ar@{-}[r]^>>>*+{\ F_{n-1}}_>*+{\ 1} &
*=<3.5pt>{\circ}
}\tag{type $D_n$}
\end{align}
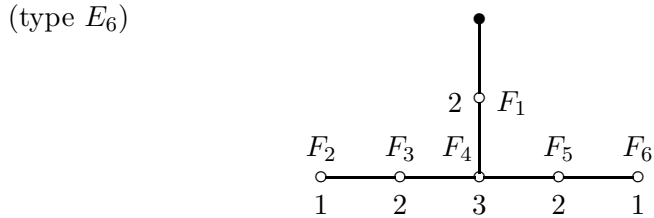
\begin{align}
\xymatrix@!0@=30pt{
&
&
*=<3.5pt>{\bullet} \ar@{-}[d] &
&
\\
&
&
*=<3.5pt>{\circ} \ar@{-}[d]^<*+{F_1}_<*+{2} &
&
\\
*=<3.5pt>{\circ}
\ar@{-}[r]^<*+{\>\!\!\!F_2}_<*+{\>\!\!\!1}^>*+{\ F_3}_>*+{\ 2} &
*=<3.5pt>{\circ} \ar@{-}[r]^>>>*+{F_4}_>*+{\ 3} &
*=<3.5pt>{\circ} \ar@{-}[r]^>*+{\ F_5}_>*+{\ 2} & 
*=<3.5pt>{\circ} \ar@{-}[r]^>*+{\ F_6}_>*+{\ 1} &
*=<3.5pt>{\circ}
}\tag{type $E_6$}
\end{align}
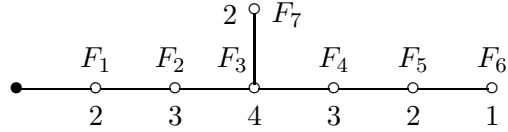
\begin{align}
\xymatrix@!0@=30pt{
&
&
&
*=<3.5pt>{\circ} \ar@{-}[d]^<*+{F_7}_<*+{2} &
&
&
\\
*=<3.5pt>{\bullet} \ar@{-}[r]^>*+{\ F_1}_>*+{\ 2} &
*=<3.5pt>{\circ} \ar@{-}[r]^>*+{F_2}_>*+{\ 3} &
*=<3.5pt>{\circ} \ar@{-}[r]^>>>*+{F_3}_>*+{\ 4} &
*=<3.5pt>{\circ} \ar@{-}[r]^>*+{\ F_4}_>*+{\ 3} & 
*=<3.5pt>{\circ} \ar@{-}[r]^>*+{\ F_5}_>*+{\ 2} &
*=<3.5pt>{\circ} \ar@{-}[r]^>*+{\ F_6}_>*+{\ 1} &
*=<3.5pt>{\circ}
}\tag{type $E_7$}
\end{align}
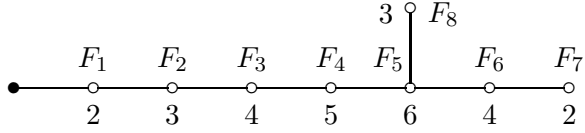
\begin{align}
\xymatrix@!0@=30pt{
&
&
&
&
&
*=<3.5pt>{\circ} \ar@{-}[d]^<*+{F_8}_<*+{3} &
&
\\
*=<3.5pt>{\bullet} \ar@{-}[r]^>*+{\ F_1}_>*+{\ 2} &
*=<3.5pt>{\circ} \ar@{-}[r]^>*+{\ F_2}_>*+{\ 3} &
*=<3.5pt>{\circ} \ar@{-}[r]^>*+{\ F_3}_>*+{\ 4} &
*=<3.5pt>{\circ} \ar@{-}[r]^>*+{\ F_4}_>*+{\ 5} &
*=<3.5pt>{\circ} \ar@{-}[r]^>>>*+{F_5}_>*+{\ 6} &
*=<3.5pt>{\circ} \ar@{-}[r]^>*+{\ F_6}_>*+{\ 4} & 
*=<3.5pt>{\circ} \ar@{-}[r]^>*+{\ F_7}_>*+{\ 2} &
*=<3.5pt>{\circ}
}\tag{type $E_8$}
\end{align}

Throughout this section,
we assume that either
$f$ belongs to one of the cases left in Theorem \ref{thm:co-prime_almost}
or $f$ is of general type with discrepancy $a/n\le1$ for $a\neq1$.
We let $H_X$ be a general hyperplane section on the germ $P \in X$,
and $H$ the birational transform on $Y$ of $H_X$.
We analyse the intersections $H \cap E$ and $S \cap E$
to obtain the dual graph for $f_S \colon S \to S_X$.
In our case we have $f^*S_X=S+(a/n)E$
since $S$ gives a general elephant of $Y$,
whence $S \cap E$ has no embedded points
and $(S \cap E)_\red$ is a tree of $\bP^1$ by Corollary \ref{cor:P1}.
On the other hand,
by Lemma \ref{lem:data} and the equality (\ref{eqn:existence_H}),
we obtain that $f^*H_X=H+E$ except for the case No 8 where $f^*H_X=H+2E$.
Anyway $H \cap E$ has no embedded points
and $(H \cap E)_\red$ is also a tree of $\bP^1$.

We proceed to study the cases left in Theorem \ref{thm:co-prime_almost}.
We use the notation and data in Lemma \ref{lem:data} freely.
By Theorems \ref{thm:co-prime_almost},
\ref{thm:co-prime_8} and \ref{thm:14/15},
we complete Theorem \ref{thm:co-prime}.
First we treat the case No 8.

\begin{theorem}\label{thm:co-prime_8}
If $f$ is of No 8,
then $a$ is co-prime to $n$
except for the case where $P$ is of type c$E/2$ and $(a,n)=(2,2)$.
This exception is placed in the last line of Table \textup{\ref{tbl:MT}}.
\end{theorem}

\begin{proof}
We suppose the case left in Theorem \ref{thm:co-prime_almost}.
Since $(-E \cdot [S \cap E])=E^3=1/6$ and $6E$ is Cartier,
we obtain that $S \cap E \cong \bP^1$ scheme-theoretically.
In particular $S$ is smooth outside $Q$ and $Q'$.
We let $C$ denote this irreducible reduced curve $S \cap E$.
$E$ is smooth at the generic point of $C$,
and the $1$-cycle $2[C]$ defines a member of $\bP H^0(\cQ_{-2,0})$.
By $H \sim 2S$,
the $1$-cycle $[H \cap E]$ also defines a member of $\bP H^0(\cQ_{-2,0})$.
Thus we obtain that $[H \cap E]=2[C]$
since $h^0(\cQ_{-2,0})=d(-2,0)=d(0,-2)=1$ by
(\ref{eqn:H^k(Q)}) and Lemma \ref{lem:data}.
Hence the coefficient of $C$ in the fundamental cycle
of the Du Val singularity $P \in S_X$,
given by that of $C$ in $(f^*H_X)|_S= H|_S+2E|_S$,
is $\ge3$,
and the dual graph of the exceptional curve of $f_S \colon S \to S_X$
is as follows.
\begin{align*}
\xymatrix@!0@=30pt{
&
*=<3.5pt>{\circ} \ar@{}[r]_<*+{\!\!\ge3} &
}
\end{align*}
Therefore $P$ is of type c$E/2$ by Table \ref{tbl:reproduce},
and the theorem is deduced.
\end{proof}

There exists an example of $f$ of No 8 with discrepancy $1$.

\begin{example}\label{exl:cE/2}
Let $P \in X$ be the germ
\begin{align*}
o \in (x_1^2+x_4^3+x_2x_3^3x_4+x_2^4+x_3^8=0) 
\subset \bC^4_{x_1x_2x_3x_4} /\frac{1}{2}(1,1,1,0).
\end{align*}
$P$ is of type c$E/2$.
Let $f$ be the weighted blow-up of $X$
with weights $\wt(x_1,x_2,x_3,x_4)=(4,2,1,3)$.
Then $f$ is a divisorial contraction of No 8 with discrepancy $1$.
\end{example}

We treat the cases No 14 and $\textrm{15}'$a simultaneously.
If $f$ is No 14 then $Q$ is a quotient singularity by $r\ge6$.
If $f$ is of No $\textrm{15}'$a then $Q$ is of type c$A/r$ by $r\ge4$,
and a general elephant of the germ $Q \in Y$ is $A_{2r-1}$ by \cite{Mo85}.

\begin{lemma}\label{lem:D_14/15}
If $a$ is not co-prime to $n$ in No 14 or $\textrm{15}'$a,
then
\begin{enumerate}
\item\label{itm:D_14/15-D}
$S_X$ has a Du Val singularity of type $D_k$ at $P$,
where $k\ge r$ in No 14 and $k \ge 2r$ in No $\textrm{15}'$a.
\item\label{itm:D_14/15-A}
$S$ has a Du Val singularity at $Q$ of type $A_{r-1}$ in No 14,
and of type $A_{2r-1}$ in No $\textrm{15}'$a.
\item\label{itm:D_14/15-proper}
The set-theoretic intersection of $H$, $E$, $S$ is the point $Q$.
\end{enumerate}
\end{lemma}

\begin{proof}
We suppose the case left in Theorem \ref{thm:co-prime_almost}.
If $S_X$ has a Du Val singularity of type $A$ at $P$,
then $H$ intersects $S \cap E$ properly at two points
by the dual graph for $f_S$.
Then $4/r=(H \cdot [S \cap E])\ge1+2/r$ in No 14
and $2/r=(H \cdot [S \cap E])\ge1+1/r$ in No $\textrm{15}'$a,
which is a contradiction.
Hence $S_X$ has a Du Val singularity of type $D$ or $E_7$ at $P$ by $n=2$.
In No 14 $(H \cdot [S \cap E])=a'E^3=2/r'$ and $r'H$ is Cartier,
and in No $\textrm{15}'$a $(H \cdot [S \cap E])=a'E^3=2/r$ and $rH$ is Cartier.
Thus
the $1$-cycle $[S \cap E]$ consists of at most two irreducible reduced cycles.
We use the trichotomy that
$S \cap E$ is either irreducible and reduced,
reduced but reducible, or irreducible but non-reduced.
In the following special case,
we assume that $S \cap E$ is reduced.

\begin{step}\label{stp:14/15-1}
\itshape
We may assume that $S \cap E$ is reduced if
$f$ is of No 14 with $(a,r)=(2,6)$ or $(4,10)$.
\end{step}

We have the surjective map $f_*\cO_Y(S) \twoheadrightarrow H^0(\cQ_{-1,0})$
by Lemma \ref{lem:QandR}.
Thus if $h^0(\cQ_{-1,0}) \ge 2$ then
the restriction on $E$ of the linear system $\system{S}$ moves,
and we can choose $S$ so that $S \cap E$ is reduced.
Since $h^0(\cQ_{-1,0})=d(-1,0)$ by (\ref{eqn:H^k(Q)}),
the estimate $h^0(\cQ_{-1,0}) \ge 2$ holds if and only if
$d(0,0)-d(-1,0) \le -1$ by (\ref{eqn:d(i,j)_non-negative}).
Computing $d(0,0)-d(-1,0)$
by (\ref{eqn:difference_of_d(i,j)}) with Lemma \ref{lem:data},
we obtain Step \ref{stp:14/15-1}.

First we prove (\ref{itm:D_14/15-proper}).

\begin{step}\label{stp:14/15-2}
\itshape
The set-theoretic intersection of $H$, $E$, $S$ is the point $Q$.
\end{step}

By $(H \cdot E \cdot S)<1$,
it suffices to prove that $H$ intersects $S \cap E$ properly.
If $S \cap E$ is non-reduced and if $(S \cap E)_\red \subset H$,
then the dual graph of the exceptional curve of $f_S$ is as follows.
\begin{align*}
\xymatrix@!0@=30pt{
Q \ar@{-}[r]_>*+{\ \ge3} &
*=<3.5pt>{\circ}
}
\end{align*}
$S$ has a Du Val singularity at $Q$ of type $A_{r-1}$ or worse in No 14,
and of type $A_{2r-1}$ or worse in No $\textrm{15}'$a.
By Step \ref{stp:14/15-1} and Lemma \ref{lem:data},
we have $r\ge10$ in No 14 and $r\ge4$ in No $\textrm{15}'$a,
whence $Q \in S$ is $A_7$ or worse.
However,
there exists no possibility for the dual graph to be described as above.
If $S \cap E$ is reduced and if $S \cap E \subset H$,
then $a=2$ and $H \cap E=S \cap E$ scheme-theoretically.
Thus $0$ is linearly equivalent to the restriction of $H-S \sim K_Y-E$ on $E$,
whence $0$ can be regarded as a member of $\bP H^0(\cQ_{1,-1})$.
This contradicts that $h^0(\cQ_{1,-1})=d(1,-1)=0$
in (\ref{eqn:H^k(Q)}) and (\ref{eqn:d(i,j)_non-negative}).

To complete Step \ref{stp:14/15-2},
we need to exclude the case where $S \cap E$ is reducible
and $H$ contains one irreducible component of $S \cap E$.
In this case the dual graph of the exceptional curves of $f_S$ is as follows.
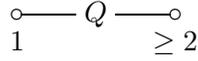
\begin{align*}
\xymatrix@!0@=30pt{
*=<3.5pt>{\circ} \ar@{-}[r]_<*+{\>\!\!\!1} &
Q \ar@{-}[r]_>*+{\ \ge2} &
*=<3.5pt>{\circ}
}
\end{align*}
$S$ has a Du Val singularity at $Q$ of type $A_{r-1}$ or worse in No 14,
and of type $A_{2r-1}$ or worse in No $\textrm{15}'$a.
$S$ is smooth outside $Q$ since $E$ is Cartier outside $Q$.
Hence $S_X$ has a Du Val singularity of type $E_7$ at $P$ and
the exceptional divisor $E|_S$ of $f_S$ has to be $F_1+F_6$ or $F_1+F_7$.
If $E|_S=F_1+F_6$ then $-(E|_S)^2=-(F_1+F_6)^2=3/4$,
which contradicts that
$-(E|_S)^2=4/r$ in No 14 and $2/r$ in No $\textrm{15}'$a.
If $E|_S=F_1+F_7$ then $(E|_S \cdot F_7)=0$,
which contradicts the $f$-ampleness of $-E$.
Therefore Step \ref{stp:14/15-2} follows.

The statement (\ref{itm:D_14/15-A}) follows from (\ref{itm:D_14/15-proper}).

\begin{step}\label{stp:14/15-3}
\itshape
$S$ has a Du Val singularity at $Q$ of type $A_{r-1}$ in No 14,
and of type $A_{2r-1}$ in No $\textrm{15}'$a.
\end{step}

On the index-one cover $Q^\sharp \in Y^\sharp$,
the local intersection number of the pre-images is
$(H^\sharp \cdot E^\sharp \cdot S^\sharp)_{Q^\sharp}=4$ in No 14
and $2$ in No $\textrm{15}'$a.
Hence if $f$ is of No 14 then $S^\sharp$ is smooth at $Q^\sharp$,
because $H^\sharp$ and $E^\sharp$ are not smooth at $Q^\sharp$.
Similarly,
if $f$ is of No $\textrm{15}'$a
then $S^\sharp$ is the restriction of a smooth divisor
on the tangent space $T_{Q^\sharp}Y^\sharp=\bC^4$ of $Y^\sharp$ at $Q^\sharp$.
Now Step \ref{stp:14/15-3} is deduced from the explicit description
of the terminal singularity $Q \in Y$.
Here we should remark the a priori property that
$S$ has a Du Val singularity at $Q$.

By Step \ref{stp:14/15-2},
we can describe the dual graph for $f_S$ as follows
according to the trichotomy.
\begin{enumerate}
\item
$S \cap E$ is irreducible and reduced.
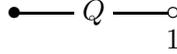
\begin{align*}
\xymatrix@!0@=30pt{
*=<3.5pt>{\bullet} \ar@{-}[r] &
Q \ar@{-}[r]_>*+{\ 1} &
*=<3.5pt>{\circ}
}
\end{align*}
\item
$S \cap E$ is reduced but reducible.
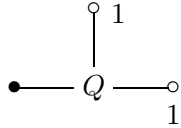
\begin{align*}
\xymatrix@!0@=30pt{
*=<3.5pt>{} &
*=<3.5pt>{\circ} \ar@{-}[d]^<*+{1} &
\\
*=<3.5pt>{\bullet} \ar@{-}[r] &
Q \ar@{-}[r]_>*+{\ 1} &
*=<3.5pt>{\circ}
}
\end{align*}
\item
$S \cap E$ is irreducible but non-reduced.
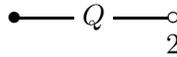
\begin{align*}
\xymatrix@!0@=30pt{
*=<3.5pt>{\bullet} \ar@{-}[r] &
Q \ar@{-}[r]_>*+{\ 2} &
*=<3.5pt>{\circ}
}
\end{align*}
\end{enumerate}
The type of the Du Val singularity $Q \in S$ is described
in Step \ref{stp:14/15-3}.
If $S \cap E$ is reduced then $S$ is smooth outside $Q$,
since $E$ is Cartier outside $Q$.

\begin{step}\label{stp:14/15-4}
\itshape
$S_X$ has a Du Val singularity of type $D$ at $P$.
\end{step}

If otherwise then $S_X$ has a Du Val singularity of type $E_7$ at $P$.
$S$ has a Du Val singularity of type $A_k$ at $Q$
with $k \ge 5$ by the estimate of $r'$ in Lemma \ref{lem:data},
and if $S \cap E$ is non-reduced then $k\ge7$ by Step \ref{stp:14/15-1}.
However,
there exists no possibility for the dual graph to be described as above,
whence Step \ref{stp:14/15-4} is deduced.

According to Step \ref{stp:14/15-3},
the number of the exceptional curves in the minimal resolution of $P \in S_X$
is $\ge r$ in No 14 and $\ge 2r$ in No $\textrm{15}'$a.
Therefore the statement (\ref{itm:D_14/15-D}) follows,
and the lemma is completed.
\setcounter{step}{0}
\end{proof}

In particular
by Lemma \ref{lem:D_14/15}(\ref{itm:D_14/15-D}) and Table \ref{tbl:reproduce},
$P$ is of type c$Ax/2$ or c$D/2$.
We consider the morphism $f' \colon (Y' \supset E') \to (X' \ni P')$,
the double cover of $f$, constructed in (\ref{eqn:covering}).
We set $\alpha_X \colon X' \to X$ and $\alpha_Y \colon Y' \to Y$
as in (\ref{eqn:covering}),
and set $Q':=\alpha_Y^{-1}(Q)$, $H':=\alpha_Y^*H$, $S':=\alpha_Y^*S$.
Let $T'$ be a general member of the linear system
$\system{-K_{Y'}}=\system{S'}=\system{a'H'}$,
and $T_{X'}'$ its birational transform on $X'$.
The surface $T'$ is a general elephant of $Y'$,
and the induced morphism $f_{T'} \colon T' \to T_{X'}'$
is a crepant morphism between normal surfaces with Du Val singularities.
We shall examine the type of the Du Val singularity $P' \in T_{X'}'$.

\begin{lemma}\label{lem:cD/2_14/15}
$T_{X'}'$ has a Du Val singularity of type $D$ at $P'$.
Consequently $P$ is of type c$D/2$.
\end{lemma}

\begin{proof}
We claim that both $Q' \in S'$ and $Q' \in T'$ are
of type $A_{r'-1}$ in No 14,
and of type $A_{2r'-1}$ in No $\textrm{15}'$a.
Indeed by Lemma \ref{lem:D_14/15}(\ref{itm:D_14/15-A}),
the class group of the germ $Q \in S$ is
$\bZ/(r)$ in No 14 and $\bZ/(2r)$ in No $\textrm{15}'$a,
and thus the claim on $Q' \in S'$ follows.
By $S'\sim a'H'\sim -K_{Y'}$,
we can regard $S'$ as a specialisation of $T'$,
whence $P' \in T'$ is at worst $A_{r'-1}$ in No 14,
and at worst $A_{2r'-1}$ in No $\textrm{15}'$a.
On the other hand by Lemma \ref{lem:variation},
the non-Gorenstein point $Q' \in Y'$ generates
one, respectively two fictitious singularities
in No 14, respectively No $\textrm{15}'$a.
Hence the claim on $Q' \in T'$ follows from the classification
of three-fold terminal singularities.

We have known that $P$ is of type c$Ax/2$ or c$D/2$.
Suppose that $T_{X'}'$ has a Du Val singularity of type $A$ at $P'$.
Equivalently,
$P$ is of type c$Ax/2$ and $P' \in S_{X'}':=\alpha_X^*(S_X)$ is of type $A$.
Then $T' \cap E'$ must intersect a general member of $\system{H'}$
properly at two points,
whence $(H' \cdot T' \cdot E')\ge1+1/r'$.
However,
$(H' \cdot T' \cdot E')=2(H \cdot S \cdot E)=4/r'$ in No 14
and $2/r'$ in No $\textrm{15}'$a,
whence $f$ has to be of No 14 with $r'=3$ by Lemma \ref{lem:data}.
In this case by Step \ref{stp:14/15-1} in Lemma \ref{lem:D_14/15},
$S \cap E$ is reduced and $S'$ has the only singularity $Q'$,
which is $A_2$.
The number of the irreducible components of $S' \cap E'$ is $\le 4$.
Hence we obtain that $P' \in S_{X'}'$ is at worst $A_6$
by counting the number of the exceptional curves
in the minimal resolution of $P' \in S_{X'}'$.
On the other hand,
$P \in S_X$ is $D_6$ or worse by Lemma \ref{lem:D_14/15}(\ref{itm:D_14/15-D}).
This contradicts Table \ref{tbl:reproduce}.
Therefore $T_{X'}'$ has a Du Val singularity of type $D$ at $P'$,
and $P$ is of type c$D/2$.
\end{proof}

We can restrict the possibility of the case $(a,n)=(4,2)$ as follows.

\begin{theorem}\label{thm:14/15}
If $f$ is of No 14 or $\textrm{15}'$a,
then $a$ is co-prime to $n$ except for the case
where $P$ is of type c$D/2$ and $(a,n)=(2,2)$ or $(4,2)$.
The case $(a,n)=(4,2)$ happens
only if $f$ is of No 14 with $r'\equiv 1$ or $7$ modulo $8$.
These exceptions are placed in
the first, second and third lines of Table \textup{\ref{tbl:MT}}.
\end{theorem}

\begin{proof}
By Theorem \ref{thm:co-prime_almost} and
Lemmata \ref{lem:data}, \ref{lem:cD/2_14/15},
it suffices to prove the claim on the case $(a,n)=(4,2)$.
Take any irreducible reduced component $C \cong \bP^1$ of $H \cap E$,
and set $C':=\alpha_Y^{-1}(C)_\red$.
If $C'$ is irreducible
then by Lemma \ref{lem:D_14/15}(\ref{itm:D_14/15-proper})
the induced morphism $C' \to C$ of degree $2$ is ramified at and only at $Q$,
which contradicts
the formula $2\chi(\cO_{C'})=2\cdot2\chi(\cO_C)-1$ of Hurwitz.
Hence $C$ splits into two $\bP^1$ by the covering $\alpha_Y$.
Write $[H \cap E]=\sum_ic_i[C_i]$ cycle-theoretically.
For each $C_i\cong\bP^1$
we choose one irreducible component $C_i^1$ of $\alpha_Y^{-1}(C_i)_\red$,
and construct a $1$-cycle $C^1:=\sum_ic_i[C_i^1]$ on $Y'$.
If $f$ is of No $\textrm{15}'$a with $(a,n)=(4,2)$,
then $(H' \cdot C^1)=1/2(H' \cdot H' \cdot E')=1/2r'$.
This contradicts that $r'H'$ is Cartier.
Hence the case $(a,n)=(4,2)$ happens only if $f$ is of No 14.
In this case $(H' \cdot C^1)=1/r'$,
whence $C^1=[\bP^1]$ and thus $H \cap E=\bP^1$.

Suppose that $f$ is of No 14 with $(a,n)=(4,2)$,
and consider a normal form of $Q \in C:=H \cap E \subset Y$.
$Q$ is a quotient singularity of type $\frac{1}{r}(1,-1,r'+4)$.
$S$ intersects $C$ properly at $Q$ with $(S \cdot C)_Q=4/r$
by Lemma \ref{lem:D_14/15}(\ref{itm:D_14/15-proper}).
Hence the restriction of the defining function of its pre-image $S^\sharp$
on $Q^\sharp \in Y^\sharp$,
whose weight is $r'+4$,
to $C^\dag$ has order $4/r$ with respect to $t$.
Thus $w_Q^C(r'+4)=4$.
Since $r'\ge5$ is an odd integer and $(r',0)\in G_Q^C$,
we may choose semi-invariant local coordinates
$x_1,x_2,x_3$ with weights $\wt(x_1,x_2,x_3)=(1,-1,r'+4)$
of the index-one cover $Q^\sharp \in Y^\sharp$ so that
$(x_1,x_2,x_3) |_{C^\dag}=(t^{(r'+1)/r},t^{(r'-1)/r},t^{4/r})$,
and so that $S^\sharp$ is given by $x_3=0$.
The defining functions $h_H$, $h_E$ of $H^\sharp$, $E^\sharp$ have
weights $2$, $-2$ respectively.
By $(H \cdot E \cdot S)_Q=4/r$,
the monomials $x_1^2$, $x_2^2$ appear in respectively $h_H$, $h_E$
with non-zero coefficients.
On the other hand by $C \subset H,E$,
the restrictions $h_H|_{C^\dag}$ and $h_E|_{C^\dag}$ are identically zero.
For the terms $x_1^2|_{C^\dag}=t^{2(r'+1)/r}$ in $h_H|_{C^\dag}$
and $x_2^2|_{C^\dag}=t^{2(r'-1)/r}$ in $h_E|_{C^\dag}$ to be cancelled out,
both of the following must hold.
\begin{enumerate}
\item
One of the monomials
$x_1x_3^{(r'+1)/4}$, $x_2x_3^{(r'+3)/4}$, $x_3^{(r'+1)/2}$ appears in $h_H$
with non-zero coefficient.
Their weights are respectively $2+r'(r'+5)/4$, $2+r'(r'+7)/4$, $2+r'(r'+5)/2$.
Hence $r'\equiv3$ modulo $8$, $r'\equiv1$ modulo $8$,
or $r'\equiv3$ modulo $4$.
\item
One of the monomials
$x_1x_3^{(r'-3)/4}$, $x_2x_3^{(r'-1)/4}$, $x_3^{(r'-1)/2}$ appears in $h_E$
with non-zero coefficient.
Their weights are respectively
$-2+r'(r'+1)/4$, $-2+r'(r'+3)/4$, $-2+r'(r'+3)/2$.
Hence $r'\equiv7$ modulo $8$, $r'\equiv5$ modulo $8$,
or $r'\equiv1$ modulo $4$.
\end{enumerate}
Therefore $r'\equiv1$ or $7$ modulo $8$,
and we complete the theorem.
\end{proof}

There exist examples of $f$ of No 14 or $\textrm{15}'$a with discrepancy $1$.

\begin{example}\label{exl:cD/2_14-1}
Let $P \in X$ be the germ
\begin{align*}
o \in \bigg(
\begin{array}{c}
x_1^2+x_4x_5+x_3^{r'+1}=0 \\
x_2^2+x_3^{r'-1}+x_4^{r'-1}+x_5=0
\end{array}
\bigg)
\subset \bC^4_{x_1x_2x_3x_4x_5} /\frac{1}{2}(1,1,1,0,0),
\end{align*}
where $r'$ is odd with $r'\neq1$.
$P$ is of type c$D/2$.
Let $f$ be the weighted blow-up of $X$
with weights $\wt(x_1,x_2,x_3,x_4,x_5)=(\frac{r'+1}{2},\frac{r'-1}{2},1,1,r')$.
Then $f$ is a divisorial contraction with discrepancy $1$.
$f$ is of No 14 with $J=\{(2r',2)\}$.
\end{example}

\begin{example}\label{exl:cD/2_15}
\begin{enumerate}
\item
Let $P \in X$ be the germ
\begin{align*}
o \in (x_1^2+x_2x_3x_4+x_2^4+x_3^4+x_4^4=0) 
\subset \bC^4_{x_1x_2x_3x_4} /\frac{1}{2}(1,1,1,0).
\end{align*}
$P$ is of type c$D/2$.
Let $f$ be the weighted blow-up of $X$
with weights $\wt(x_1,x_2,x_3,x_4)=(2,2,1,1)$.
Then $f$ is a divisorial contraction with discrepancy $1$.
$f$ is of No $\textrm{15}'$a with $J=\{(4,1),(4,1)\}$.
\item
Let $P \in X$ be the germ
\begin{align*}
o \in (x_1^2+x_2^2x_4+x_3^{2r'}+x_4^{2r'}=0) 
\subset \bC^4_{x_1x_2x_3x_4} /\frac{1}{2}(1,1,1,0),
\end{align*}
where $r'$ is even.
$P$ is of type c$D/2$.
Let $f$ be the weighted blow-up of $X$
with weights $\wt(x_1,x_2,x_3,x_4)=(r',r',1,1)$.
Then $f$ is a divisorial contraction with discrepancy $1$.
$f$ is of No $\textrm{15}'$a with $J=\{(2r',1),(2r',1)\}$.
\end{enumerate}
\end{example}

There also exist examples of $f$ of No 14 with discrepancy $2$.

\begin{example}\label{exl:cD/2_14-2}
\begin{enumerate}
\item
Let $P \in X$ be the germ 
\begin{align*}
o \in \bigg(
\begin{array}{c}
x_1^2+x_4x_5+x_2x_3^{(r'+3)/4}=0 \\
x_2^2+x_3^{(r'-1)/2}+x_4^{r'-1}+x_5=0
\end{array}
\bigg)
\subset \bC^4_{x_1x_2x_3x_4x_5} /\frac{1}{2}(1,1,1,0,0),
\end{align*}
where $r'\equiv1$ modulo $8$ with $r'\neq1$.
$P$ is of type c$D/2$.
Let $f$ be the weighted blow-up of $X$
with weights $\wt(x_1,x_2,x_3,x_4,x_5)=(\frac{r'+1}{2},\frac{r'-1}{2},2,1,r')$.
Then $f$ is a divisorial contraction with discrepancy $2$.
$f$ is of No 14 with $J=\{(2r',2)\}$.
\item
Let $P \in X$ be the germ
\begin{align*}
o \in \bigg(
\begin{array}{c}
x_1^2+x_4x_5+x_3^{(r'+1)/2}=0 \\
x_2^2+x_1x_3^{(r'-3)/4}+x_4^{r'-1}+x_5=0
\end{array}
\bigg)
\subset \bC^4_{x_1x_2x_3x_4x_5} /\frac{1}{2}(1,1,1,0,0),
\end{align*}
where $r'\equiv7$ modulo $8$.
$P$ is of type c$D/2$.
Let $f$ be the weighted blow-up of $X$
with weights $\wt(x_1,x_2,x_3,x_4,x_5)=(\frac{r'+1}{2},\frac{r'-1}{2},2,1,r')$.
Then $f$ is a divisorial contraction with discrepancy $2$.
$f$ is of No 14 with $J=\{(2r',2)\}$.
\end{enumerate}
\end{example}

The rest of this section is devoted to the proof of Theorem \ref{thm:hard_GE}
when $f$ is of general type with discrepancy $a/n\le1$.
We keep the notation given just before Theorem \ref{thm:hard_GE}.
If $S_X$ has a Du Val singularity of type $A$ at $P$,
then $H$ intersects $S \cap E$ properly at two points
by the dual graph for $f_S$,
and consequently Theorem \ref{thm:hard_GE}(\ref{itm:hard_A}) holds
by Lemma \ref{lem:proper_intersect}(\ref{itm:proper_S}).
Hence we assume that $S_X$ has a Du Val singularity of type $D$ or $E$ at $P$.
In particular $n=2$, $3$ or $4$.
We shall describe the partial resolution $f_S \colon S \to S_X$.

\begin{lemma}\label{lem:15/16small}
Suppose that $f$ is of general type with discrepancy $a/n\le1$ for $a\neq1$,
and that $S_X$ has a Du Val singularity of type $D$ or $E$ at $P$.
Then the dual graph of the exceptional curves of $f_S$ is as follows.
\begin{enumerate}
\item
$f$ is of No $\textrm{15}''$.
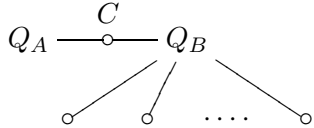
\begin{align*}
\xymatrix@!0@=15pt{
Q_A \ar@{-}[rr]^>*+{\ C} &
&
*=<3.5pt>{\circ}\ar@{-}[rr] &
&
Q_B \ar@{-}[ddlll]\ar@{-}[ddl]\ar@{-}[ddrrr] &
&
&\\
&&&&&&&\\
&
*=<3.5pt>{\circ} &
&
*=<3.5pt>{\circ} &
&
\cdots\cdot &
&
*=<3.5pt>{\circ}
}
\end{align*}
The curve $C$ is contained in $H|_S$.
$(A,B)=(1,2)$ or $(2,1)$
and $S$ has Du Val singularities of types $A_{r_A-1}$, $A_{r_B-1}$ or worse
respectively at $Q_A$, $Q_B$.
$S$ has no singularities other than $Q_1$ and $Q_2$ on
any exceptional curve at the generic point of which $S \cap E$ is smooth.
\item
$f$ is of No 16.
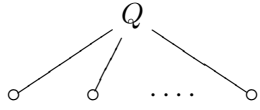
\begin{align*}
\xymatrix@!0@=15pt{
&
&
&
Q \ar@{-}[ddlll]\ar@{-}[ddl]\ar@{-}[ddrrr] &
&
&\\
&&&&&&\\
*=<3.5pt>{\circ} &
&
*=<3.5pt>{\circ} &
&
\cdots\cdot &
&
*=<3.5pt>{\circ}
}
\end{align*}
$S$ has a Du Val singularity of type $A_{r-1}$ or worse at $Q$.
There exists an exceptional curve
at the generic point of which $S \cap E$ is smooth,
and on any such curve,
$S$ has no singularities other than $Q$.
If $\bullet$ intersects an irreducible reduced exceptional curve
at a point other than $Q$,
then this curve is contained in $H|_S$.
\end{enumerate}
\end{lemma}

\begin{proof}
First we consider the case No $\textrm{15}''$.
We claim that $S \cap E$ has an irreducible component
which passes through $Q_1$ and $Q_2$.
If otherwise then  $S \cap E$ is a tree of
$Q_1 \in C_1 \cong \bP^1$ and $Q_2 \in C_2 \cong \bP^1$
by $(H \cdot S \cdot E)=1/r_1+1/r_2$.
Hence $S \cap E \subset H$
by Lemma \ref{lem:proper_intersect}(\ref{itm:proper_S}).
This could happen only if $a/n=1$ and
$H \cap E=S \cap E$ scheme-theoretically.
However in this case,
$0$ is linearly equivalent to the restriction of $H-S \sim K_Y-E$ on $E$,
because $E$ is smooth at the generic points of $C_1$ and $C_2$.
Thus $0$ can be regarded as a member of $\bP H^0(\cQ_{1,-1})$,
which contradicts that $h^0(\cQ_{1,-1})=d(1,-1)=0$
in (\ref{eqn:H^k(Q)}) and (\ref{eqn:d(i,j)_non-negative}).
Hence $S \cap E$ has an irreducible component $C$
which passes through $Q_1$ and $Q_2$.
Such a curve $C$ is uniquely determined
because $(S \cap E)_\red$ is a tree of $\bP^1$.
$C$ is contained in $H$
by Lemma \ref{lem:proper_intersect}(\ref{itm:proper_S}),
whence we obtain the desired dual graph by $(H \cdot S \cdot E)=1/r_1+1/r_2$,
remarking that $H$ and $E$ are Cartier outside $Q_1,Q_2$.

We then consider the case No 16.
In this case $(a,n)=(2,3)$ or $(3,4)$ by
the co-primeness of $a$ and $n$ in Theorem \ref{thm:ramification},
and $r \ge a+n$ by Lemma \ref{lem:minimal_discrepancy}.
By the surjective map
$f_*\cO_Y(S) \twoheadrightarrow H^0(\cQ_{-1,0})$ in Lemma \ref{lem:QandR}
and the estimate $h^0(\cQ_{-1,0})=2+a/n+(a/n-b)/r\ge2$
in (\ref{eqn:existence_S}),
there exists an irreducible component of $S \cap E$
at the generic point of which $S \cap E$ is smooth.
The desired dual graph is obtained
by Lemma \ref{lem:proper_intersect}(\ref{itm:proper_S})
provided that any exceptional curve passes through $Q$.

Suppose that there exists an exceptional curve not passing through $Q$.
Then $S \cap E$ is a tree of
$Q \in C_1 \cong \bP^1$ and $Q \not\in C_2 \cong \bP^1$ with a node $Q'$
by $(H \cdot S \cdot E)=1+1/r$.
In particular $(H \cdot C_1)=1/r$,
and $S$ is smooth outside $Q,Q'$ since $E$ is Cartier outside $Q$.
Hence $Q'$ is a hidden non-Gorenstein point of index $n$
by $(S \cdot C_1)=a/nr$.
We have $C_1 \subset H$
by Lemma \ref{lem:proper_intersect}(\ref{itm:proper_S}).
Thus we obtain the following dual graph of
the exceptional curves of $f_S$.
\begin{align*}
\xymatrix@!0@=30pt{
Q \ar@{-}[r]_>*+{\ \ge2}^>*+{\ C_1} &
*=<3.5pt>{\circ} \ar@{-}[r] &
Q' \ar@{-}[r]_>*+{\ \ge1}^>*+{\ C_2} &
*=<3.5pt>{\circ}
}
\end{align*}
$S$ has Du Val singularities of types $A_{r-1}$, $A_{n-1}$ or worse
respectively at $Q$, $Q'$,
and they are all the singularities of $S$.
If $(a,n)=(2,3)$ then
$P \in S_X$ is of type $E_6$ by Table \ref{tbl:reproduce},
but there exists no possibility for the dual graph to be described as above.
If $(a,n)=(3,4)$ then
$P \in S_X$ is of type $D_k$ by Table \ref{tbl:reproduce},
and $(C_1,C_2)=(F_l,F_1)$ for $5 \le l\le k-3$
or $(F_l,F_k)$ for $r \le l \le k-4$ up to permutation.
If $(C_1,C_2)=(F_l,F_1)$ then $E|_S=F_1+F_l$ and
$(E|_S \cdot C_1)=(F_1+F_l \cdot F_l)=0$,
which is a contradiction.
If $(C_1,C_2)=(F_l,F_k)$
then $E|_S=F_l+F_k$ and $-(E|_S)^2=-(F_l+F_k)^2=1/l+1/(k-l)\le1/r+1/4$,
which also contradicts that $-(E|_S)^2=-(E \cdot E \cdot S)=1+1/r$.
Therefore any exceptional curve passes through $Q$,
and the lemma is deduced.
\end{proof}

We shall narrow down the possibilities of $f_S$
by using Lemma \ref{lem:15/16small}.

\begin{lemma}\label{lem:15/16small_non-E}
If $f$ is of general type with discrepancy $a/n\le1$ for $a\neq1$,
then $S_X$ has a Du Val singularity of type $A$ or $D$ at $P$.
\end{lemma}

\begin{proof}
Indeed if $S_X$ has a Du Val singularity of type $E_7$,
then $(a,n)=(2,2)$ and thus $f$ is of No $\textrm{15}''$
by Theorem \ref{thm:ramification}.
Moreover, $2 \mid r_1,r_2$ by Theorem \ref{thm:ramification} and
$4 \mid r_1,r_2$ by $e_{Q_1}=r_1/2+1$, $e_{Q_2}=r_2/2+1$
in Lemma \ref{lem:d_Q}(\ref{itm:d_Q_ge}).
However, there exists no possibility for the dual graph to be described
as in Lemma \ref{lem:15/16small}.

If $S_X$ has a Du Val singularity of type $E_6$,
then $(a,n)=(2,3)$ or $(3,3)$ and thus
$r_2 \ge 5$ by Lemma \ref{lem:minimal_discrepancy}.
If $f$ is of No $\textrm{15}''$ then $E|_S$ has to be $F_3$ up to permutation
by Lemma \ref{lem:15/16small}.
In this case $r_1=2$, $r_2=5$ and $-(E|_S)^2=1/r_1+1/r_2=7/10$,
which contradicts that $-(F_3)^2=3/10$.
On the other hand,
if $f$ is of No 16 then $E|_S$ has to be
either $F_1$, $F_2$, $F_1+F_2$ or $F_2+F_6$ up to permutation
by Lemma \ref{lem:15/16small}.
Then $-(E|_S)^2$ is respectively
$-(F_1)^2=1/2$, $-(F_2)^2=3/4$, $-(F_1+F_2)^2=4/5$, $-(F_2+F_6)^2=1$,
which contradicts that $-(E|_S)^2=1+1/r$.
\end{proof}

Therefore for Theorem \ref{thm:hard_GE},
the case where $P \in S_X$ is of type $D$ remains.
We list all the possibilities of the exceptional locus of $f_S$
in Tables \ref{tbl:15D} and \ref{tbl:16D} up to permutation
by Lemma \ref{lem:15/16small}.

\begin{table}[ht]
\caption{$f$ is of No $\textrm{15}''$
and $P \in S_X$ is of type $D_k$.}\label{tbl:15D}
\begin{tabular}{l|l|c|c}
\hline
\multicolumn{1}{c|}{$E|_S$} & \multicolumn{1}{c|}{$-(E|_S)^2$} &
$r_A$ & $r_B$ \\
\hline
$F_1+F_l$, $2\le l\le k-3$         & $1$             &         & $\le l-1$   \\
$F_l$, $2\le l\le k-3$             & $1/l$           & $\le l$ &             \\
$F_l+2F_m$, $2\le l< m\le k-2$     & $1/l+1/(m-l)$   & $\le l$ & $\le m-l$   \\
$2F_l+F_m$, $2\le l< m\le k-3$     & $4/l+1/(m-l)$   &         & $\le m-l$   \\
$F_l+F_k$, $2\le l\le k-3$         & $1/l+1/(k-l)$   & $\le l$ & $\le k-l$   \\
$F_l+F_{k-1}+F_k$, $2\le l\le k-3$ & $1/l+1/(k-l-1)$ & $\le l$ & $\le k-l-1$ \\
\hline
\end{tabular}
\end{table}

\begin{table}[ht]
\caption{$f$ is of No 16 and $P$ is of type $D_k$.}\label{tbl:16D}
\begin{tabular}{l|l|c}
\hline
\multicolumn{1}{c|}{$E|_S$} & \multicolumn{1}{c|}{$-(E|_S)^2$} & $r$ \\
\hline
$F_1$                       & $1$         &           \\
$F_1+2F_l$, $2\le l\le k-2$ & $1+1/(l-1)$ & $\le l-1$ \\
$F_1+F_k$                   & $1+1/(k-1)$ & $\le k-1$ \\
$F_1+F_{k-1}+F_k$           & $1+1/(k-2)$ & $\le k-2$ \\
$F_k$                       & $4/k$       & $\le k$   \\
$F_{k-1}+F_k$               & $4/(k-1)$   & $\le k-1$ \\
\hline
\end{tabular}
\end{table}

In Table \ref{tbl:15D},
the case $E|_S=F_1+F_l$ does not occur by $(E|_S \cdot F_l)=0$,
the case $E|_S=F_l$ does not occur by $-(E|_S)^2=1/r_1+1/r_2$,
and the case $E|_S=2F_l+F_m$ does not occur by $(E|_S \cdot F_m)=1/(m-l)>0$.
If $f$ is of No 16 and if $P \in S_X$ is of type $D$,
then $(a,n)=(3,4)$ by Theorem \ref{thm:ramification}
and $r\ge7$ by Lemma \ref{lem:minimal_discrepancy}.
Hence the cases $E|_S=F_1$, $F_k$, $F_{k-1}+F_k$ in Table \ref{tbl:16D}
do not occur by $-(E|_S)^2=1+1/r$.
By this observation and the equality $-(E|_S)^2=1/r_1+1/r_2$
with Tables \ref{tbl:15D} and \ref{tbl:16D},
we obtain the following lemma.
Here the curve $C$ in Lemma \ref{lem:15/16small_D} stands for
$F_l$ in Table \ref{tbl:15D} for No $\textrm{15}''$,
and $F_1$ in Table \ref{tbl:16D} for No 16.

\begin{lemma}\label{lem:15/16small_D}
Theorem \textup{\ref{thm:hard_GE}} holds
except for the case where $f$ satisfies the following conditions.
\begin{enumerate}
\item
$S_X$ has a Du Val singularity of type $D$ at $P$.
In particular $n=2$ or $4$.
\item
$S$ has a Du Val singularity of type $A_{r_Q-1}$ at each $Q \in I$.
\item
$[S \cap E]=[C]+[V]$ cycle-theoretically,
where $[V]$ consists of one or two $[\bP^1]$.
$V$ intersects $C$ at and only at
$Q_2$ in No $\textrm{15}''$ and $Q$ in No 16.
\item
$[S \cap H]=[C]+[B]$ in No $\textrm{15}''$
and $[S \cap H]=[B]$ in No 16 cycle-theoretically,
where $B$ is the birational transform on $S$ of $H_X|_{S_X}$.
$B$ intersects $E$ at and only at
$Q_1$ in No $\textrm{15}''$ and $Q$ in No 16,
except for the case $r_1=2$ in No $\textrm{15}''$.
\item
$(H \cdot C)=1/r_1-1/r_2$.
In particular $r_1<r_2$.
\end{enumerate}
\end{lemma}

It is the time to complete Theorem \ref{thm:hard_GE}.

\begin{theorem}\label{thm:15/16small}
If $f$ is of general type with discrepancy $a/n\le1$,
then Theorem \textup{\ref{thm:hard_GE}} holds.
\end{theorem}

\begin{proof}
It suffices to study the case left in Lemma \ref{lem:15/16small_D}.
We derive Theorem \ref{thm:hard_GE}(\ref{itm:hard_D}) in this case.
Since $S \cap E$ and thus $S$ are smooth along $C$
outside the singularities in $I$,
there exists no hidden non-Gorenstein point on $C$.
We first exclude the case No 16.
If $f$ is of No 16 then $(a,n)=(3,4)$ by Theorem \ref{thm:ramification}.

\begin{step}\label{stp:15/16small-1}
\itshape
$f$ is of No $\textrm{15}''$.
If $(a,n)=(3,4)$ then $V$ is non-reduced.
\end{step}

Indeed,
if $(a,n)=(3,4)$ then we have the divisor $D_{1,-1}=K_Y-E\equiv(1/4)H$
and $(D_{1,-1} \cdot V)=1/2r_2$.
Hence $V$ has to pass through a hidden non-Gorenstein point, say $Q'$.
In particular $V$ is non-reduced
since $V$ is the intersection of $S$ and $E$ outside $C$.
We shall exclude the case No 16.
If $f$ is of No 16 then $(D_{1,-1} \cdot C)=(r-1)/4r$,
whence $(r-1)/4$ is an integer.
On the other hand,
also $(r+1)/4$ must be an integer
by the equality $h^0(\cQ_{-1,0})=2+(3(r+1)/4-b)/r$ in (\ref{eqn:existence_S}).
This is a contradiction,
and Step \ref{stp:15/16small-1} is deduced.

For the case No $\textrm{15}''$,
we shall construct good semi-invariant local coordinates 
$x_{i1},x_{i2},x_{i3}$ with weights $\wt(x_{i1},x_{i2},x_{i3})=(1,-1,b_i)$
of the index-one cover $Q_i^\sharp \in Y^\sharp$
at each non-Gorenstein point $Q_i$ on $C$ for $i=1,2$.

\begin{step}\label{stp:15/16small-2}
\itshape
We may choose the coordinates $x_{i1},x_{i2},x_{i3}$ so that
\begin{align*}
{Q_i}^\sharp \in C^\sharp \subset Y^\sharp
\cong o \in (x_{ii}\axis) \subset \bC^3_{x_{i1}x_{i2}x_{i3}}.
\end{align*}
Moreover, 
the pre-images $E^\sharp$, $S^\sharp$ on $Q_1^\sharp \in Y^\sharp$
of $E$, $S$ are given by respectively $x_{12}=0$, $x_{13}=0$,
and the pre-images $H^\sharp$, $S^\sharp$ on $Q_2^\sharp \in Y^\sharp$
of $H$, $S$ are given by respectively $x_{21}=0$, $x_{23}=0$.
\end{step}

$S$ has a Du Val singularity of type $A_{r_i-1}$ at each $Q_i$.
Thus $S^\sharp$ is smooth at $Q_i^\sharp$
since the class group of the germ $Q_i \in S$ is $\bZ/(r_i)$.
Hence we may choose the coordinates so that
$Q_i^\sharp \in S^\sharp$ is given by $x_{i3}=0$ for each $i$.
The equality $(H \cdot V)_{Q_2}=2/r_2$ implies that
$H^\sharp$ has multiplicity $\le 2$ at $Q_2^\sharp$,
and $H^\sharp$ has multiplicity $2$ only if $V \cong \bP^1$.
If $H^\sharp$ is not smooth at $Q_2^\sharp$,
then $(a,n)\neq(3,4)$ by Step \ref{stp:15/16small-1}.
Hence $a$ and $n$ has a common divisor $2$
and thus $r_2$ is even by Theorem \ref{thm:ramification}.
This contradicts the equality
$\wt x_{2i} +\wt x_{2j} \equiv 1$ modulo $r_2$ for some $i,j=1$, $2$ or $3$.
Thus $H^\sharp$ is smooth at $Q_2^\sharp$.
In order to obtain the coordinates $x_{21},x_{22},x_{23}$
in Step \ref{stp:15/16small-2},
we have only to prove that $b_2\neq1$.
Since $nb_2 \equiv a$ modulo $r_2$ and $r_2>r_1\ge2$,
the equality $b_2=1$ holds only if $(a,n)=(2,2)$ or $(4,4)$.
However if $(a,n)=(2,2)$ or $(4,4)$ and $b_2=1$,
then $e_2=1$ and $d_2^n=0$,
which contradicts Lemma \ref{lem:d_Q}(\ref{itm:d_Q_ge}).
Therefore $b_2\neq1$ and
we obtain the desired coordinates $x_{21},x_{22},x_{23}$.
In particular, $w_{Q_2}^C(1)=r_2-1$.

By the map
$\cO_Y(H)^{\otimes r_1r_2} \otimes \cO_Y(-r_1r_2H) \otimes \cO_C \to \cO_C$,
we have $s_C(0,-1)-(H \cdot C)=-w_{Q_1}^C(1)/r_1-w_{Q_2}^C(1)/r_2$.
By $w_{Q_2}^C(1)=r_2-1$ and $(H \cdot C)=1/r_1-1/r_2$,
we obtain that $w_{Q_1}^C(1) \equiv 1$ modulo $r_1$.
Hence we may choose the coordinates $x_{11},x_{12},x_{13}$
which satisfy all the conditions in Step \ref{stp:15/16small-2}
but the condition that $E^\sharp$ is given by $x_{12}=0$ at $Q_1^\sharp$.
However,
we can obtain this condition
because $C$ is the intersection of $S$ and $E$ near $Q_1$
and its pre-image $C^\sharp$ at $Q_1^\sharp$ is given by the $x_{11}$-axis.

\begin{step}\label{stp:15/16small-3}
\itshape
$(a,n)=(2,2)$.
$4 \mid r_i$ and $b_i=r_i/2+1$ for each $i$.
\end{step}

By the map
$\cO_Y(S)^{\otimes r_1r_2} \otimes \cO_Y(-r_1r_2S) \otimes \cO_C \to \cO_C$,
we have $s_C(-1,0)-(S \cdot C)=-w_{Q_1}^C(b_1)/r_1-w_{Q_2}^C(b_2)/r_2$.
Since $w_{Q_1}^C(b_1)=b_1$ and $w_{Q_2}^C(b_2)=r_2-b_2$
by Step \ref{stp:15/16small-2},
we have $(S \cdot C)=b_1/r_1-b_2/r_2+s_C(-1,0)+1$.
On the other hand $(S \cdot C)=(a/n)(H \cdot C)=(a/n)(1/r_1-1/r_2)$
by Lemma \ref{lem:15/16small_D},
whence $(a/n-b_1)/r_1 \equiv (a/n-b_2)/r_2$ modulo $\bZ$.
Thus we obtain that $(a,n)=(2,2)$ as in the proof of Lemma \ref{lem:n=2}.
The rest of Step \ref{stp:15/16small-3} follow from $e_i=r_i/2+1$
in Lemma \ref{lem:d_Q}(\ref{itm:d_Q_ge}).

By Lemma \ref{lem:15/16small_D} and the above steps,
it remains only proving (b) in Theorem \ref{thm:hard_GE}(\ref{itm:hard_D}).
However the proof of Step \ref{stp:hard_GE15-3} in Theorem \ref{thm:hard_GE15}
is valid also in our situation.
Therefore the theorem is completed.
\setcounter{step}{0}
\end{proof}

\section{Divisorial contractions to points of type c$A/n$}\label{sec:cA/n}
In this section we prove
the ordinary case (\ref{itm:MTord}) of Theorem \ref{thm:MT}.
Our strategy for the classification
is to determine $E$ in the sense of valuation
by applying the next lemma,
which is a reformulation of \cite[Lemma 6.1]{Ka03}.
We recall that for a proper birational morphism $Z \to X$,
the \textit{centre} of $E$ on $Z$ is the locus $\pi_Z(\pi_{Y*}^{-1}E)$
for a common resolution of singularities $W$ of $Y$ and $Z$
with $\pi_Y \colon W \to Y$, $\pi_Z \colon W \to Z$.

\begin{lemma}\label{lem:gen_mtd}
Let $f \colon (Y \supset E) \to (X \ni P)$ be a germ of a
three-fold divisorial contraction to a point $P$ of index $n$
identified with
\begin{align*}
P \in X \cong o \in (\phi=0) \subset
\bar{X}:=\bC^4_{x_1x_2x_3x_4} /\frac{1}{n}(\wt x_1,\wt x_2,\wt x_3,\wt x_4).
\end{align*}
Let $a/n$ denote the discrepancy of $f$,
and $m_i/n$ the multiplicity of $x_i$ along $E$.
Suppose that $(m_1/n,m_2/n,m_3/n,m_4/n)$ is primitive in the group
$N:=\bZ^4+\bZ(\wt x_1/n,\wt x_2/n,\wt x_3/n,\wt x_4/n)$.
Let $d/n$ denote the weighted order of $\phi$ with respect to weights
$\wt(x_1,x_2,x_3,x_4)=(m_1/n,m_2/n,m_3/n,m_4/n)$,
and decompose $\phi$ as
\begin{align*}
\phi=\phi_{d/n}(x_1,x_2,x_3,x_4)+\phi_{>d/n}(x_1,x_2,x_3,x_4),
\end{align*}
where $\phi_{d/n}$ is the weighted homogeneous part of weight $d/n$
and $\phi_{>d/n}$ is the part of weight $>d/n$.
Set $c/n:=m_1/n+m_2/n+m_3/n+m_4/n-1-d/n$.
Let $\bar{g} \colon (\bar{Z} \supset \bar{F}) \to (\bar{X} \ni P)$
be the weighted blow-up of $\bar{X}$ with weights
$\wt(x_1,x_2,x_3,x_4)=(m_1/n,m_2/n,m_3/n,m_4/n)$,
$\bar{F}$ its exceptional divisor,
and $\bar{D}_i$ the birational transform on $\bar{Z}$
of the $\bQ$-Cartier divisor $x_i=0$ on $\bar{X}$.
Let $Z$ denote the birational transform on $\bar{Z}$ of $X$,
and $g \colon Z \to X$ the induced morphism.
Suppose that $\bar{F} \cap Z$ defines an irreducible reduced $2$-cycle $[F]$
on $Z$.
We suppose that $Z$ is smooth at the generic point of $F$,
and that $\dim(\Sing\bar{Z} \cap Z) \le 1$.
Let $L$ denote the centre of $E$ on $Z$.
Then,
\begin{enumerate}
\item\label{itm:non-contained}
$L \not\subseteq \bigcup_{1 \le i \le 4} \bar{D}_i$,
and $F$ has multiplicity $1$ along $E$.
\item\label{itm:L_non-surface}
If $L\neq F$ and $Z$ is canonical at the generic point of $L$,
then $c/n<a/n$.
\item\label{itm:L_surface}
If $L=F$,
then $c/n=a/n$ and $f \cong g$ over $X$.
\end{enumerate}
\end{lemma}

\begin{proof}
$Z$ is normal by the criterion of Serre,
and $F$ is a $\bQ$-Cartier anti-$f$-ample divisor.
Let $M$ and $M_i$ for $1\le i\le 4$ denote
the multiplicities along $E$ of $F$ and $\bar{D}_i|_Z$.
Then $m_i=Mm_i+M_i$ for each $i$.
On the other hand, at least one $M_i$ equals zero
by $\bigcap_{1 \le i \le 4} \bar{D}_i=\emptyset$,
whence $M=1$ and $M_i=0$ for all $i$.
This implies (\ref{itm:non-contained}).
The rest (\ref{itm:L_non-surface}) and (\ref{itm:L_surface}) follow from $M=1$,
the adjunction formula $K_Z=g^*K_X+(c/n)F$ and \cite[Lemma 3.4]{Ka01}.
\end{proof}

The morphisms $\bar{g}$ and $g$ in Lemma \ref{lem:gen_mtd} are
described in terms of toric geometry.
In the situation of Lemma \ref{lem:gen_mtd},
we identify $P \in \bar{X}$ with the germ of the toric variety
corresponding to the fan generated by
$\vec{e}_1:=(1,0,0,0)$, \ldots, $\vec{e}_4:=(0,0,0,1)$ in $N$,
and $\bar{f}$ with the operation of adding the ray generated
by the primitive vector $\vec{e}:=(m_1/n,m_2/n,m_3/n,m_4/n)$.
Then $\bar{Y}$ is covered by four charts $\bar{U}_i$ for $1\le i\le 4$,
where $\bar{U}_i$ corresponds to the fan generated by
$\vec{e}$ and all $\vec{e}_j$ with $j\neq i$.
Note that $\bar{U}_i=\bar{Z} \setminus \bar{D}_i$.

We proceed to the classification in the c$A/n$ case.
From now on,
we assume that $f \colon (Y \supset E) \to (X \ni P)$ is
a divisorial contraction to a point of type c$A/n$ with $n \ge 2$.
The classification of $f$ has been obtained
when $P$ is a terminal quotient singularity by Kawamata in \cite{Km96},
and when the discrepancy $a/n$ is $1/n$
by Hayakawa in \cite[Theorem 6.4]{Ha99}.
We can immediately verify that
Theorem \ref{thm:MT}(\ref{itm:MTord}) holds in these cases.
Therefore in this section
we assume that $a\ge2$ and that $P$ is not a quotient singularity.
Then in particular $f$ is of general type by
Theorem \ref{thm:co-prime} and Corollary \ref{cor:a=1},
and thus we can use Theorem \ref{thm:hard_GE},
the general elephant theorem in strong form.
The dual graph belongs to
the case (\ref{itm:hard_A}) of Theorem \ref{thm:hard_GE}.

We follow the notation in Section \ref{sec:GE}.
We let $H_X$ be a general hyperplane section on the germ $P \in X$
and $H$ its birational transform on $Y$,
let $S$ be a general elephant of $Y$ and
$S_X$ its birational transform on $X$.
Then $f^*H_X=H+E$ by (\ref{eqn:existence_H}),
and $f^*S_X=S+(a/n)E$.
According to Theorem \ref{thm:hard_GE}(\ref{itm:hard_A}),
the surfaces $H$, $E$, $S$ intersect properly at and only at two points
$Q_1$ and $Q_2$ of $Y$.
$Q_1$, $Q_2$ correspond
to the points $Q_1$, $Q_2$ in $I$ in No $\textrm{15}''$,
and to a smooth point and $Q$ in $I$ in No 16.
We let $r_1$, $r_2$ denote the local indices of $Y$ at $Q_1$, $Q_2$.
By $(H\cdot E\cdot S)=(a/n)E^3=1/r_1+1/r_2$ in Table \ref{tbl:classification},
we have the local intersection number $(H\cdot E\cdot S)_{Q_i}=1/r_i$
for each $i$.
In particular we may choose semi-invariant local coordinates
$x_{i1},x_{i2},x_{i3}$ of the index-one cover $Q_i^\sharp \in Y^\sharp$
of the germ $Q_i \in Y$ so that
the pre-images $H^\sharp$, $E^\sharp$, $S^\sharp$ on $Q_i^\sharp \in Y^\sharp$
of $H$, $E$, $S$ are given by
\begin{align}\label{eqn:description}
H^\sharp=(x_{i1}=0), && E^\sharp=(x_{i2}=0), && S^\sharp=(x_{i3}=0).
\end{align}
We let $h_{aH_X}$, $h_{nS_X}$ denote the defining functions on $X$
of $aH_X$, $nS_X$ respectively.
We can regard $h_{aH_X}$ and $h_{nS_X}$ as elements
in the vector space $V_a:=f_*\cO_Y(-aE)/f_*\cO_Y(-(a+1)E)$.
The following lemma plays a pivotal role in applying Lemma \ref{lem:gen_mtd}.

\begin{lemma}\label{lem:pivotal}
$h_{aH_X}$ and $h_{nS_X}$ are linearly independent in $V_a$.
\end{lemma}

\begin{proof}
The statement means that
any function $h:=c_1h_{aH_X}+c_2h_{nS_X}$
for $c_1,c_2\in\bC$ with $(c_1,c_2)\neq(0,0)$ has multiplicity $a$ along $E$.
By the description (\ref{eqn:description}),
on the germ $Q_i \in Y$
the function $h_{aH_X}$, $h_{nS_X}$ are expressed as
$h_{aH_X}=u_1x_{i1}^ax_{i2}^a$ and $h_{nS_X}=u_2x_{i3}^nx_{i2}^a$,
where $u_1$, $u_2$ are units.
Then $h$ is expressed as $h=(c_1u_1x_{i1}^a+c_2u_2x_{i3}^n)x_{i2}^a$,
whence $h$ vanishes with multiplicity $a$ along $x_{i2}=0$,
that is, along $E$.
Thus the lemma is deduced.
\end{proof}

We shall choose a good identification of $P \in X$ in Theorem \ref{thm:Mo85}.

\begin{lemma}\label{lem:id}
We may choose the identification
\begin{align}\label{eqn:id}
P \in X \cong o \in (\phi:=x_1x_2+g(x_3^n,x_4)=0) \subset
\bC^4_{x_1x_2x_3x_4} /\frac{1}{n}(1,-1,b,0)
\end{align}
so that $h_{nS_X}\equiv x_3^n$, $h_{aH_X}\equiv x_4^a$
modulo $f_*\cO_Y(-(a+1)E)$.
\end{lemma}

\begin{proof}
Take the index-one cover $P^\sharp \in X^\sharp$ of $P \in X$,
and let $h_{S_X}$, $h_{H_X}$ denote the defining functions of
the pre-images of $S_X$, $H_X$.
The lemma follows from standard argument on coordinates change.
Since $S_X$ has a Du Val singularity of type $A$ at $P$,
we can take the identification (\ref{eqn:id}) so that
$h_{S_X}=x_3+(\textrm{non-linear terms})$,
see \cite[Corollary 4]{Mo85} for example.
Since $H_X$ is a general hyperplane section on $X$,
we can set $h_{H_X}\equiv x_4+x_1p(x_1,x_3)+q(x_2,x_3)$ modulo $f_*\cO_Y(-2E)$.
By making coordinates change in which $x_4+x_1p \mapsto x_4$,
we eliminate $x_1$ from $h_{H_X}$ modulo $f_*\cO_Y(-2E)$.
The same argument yields the coordinates in which
$h_{H_X}\equiv x_4$ modulo $f_*\cO_Y(-2E)$.

We have preserved the coordinate $x_3$ in the above argument,
whence $h_{S_X}=x_3+(\textrm{non-linear terms})$.
Suppose that we can write
$h_{H_S}\equiv x_3+u_1x_1p(x_1,x_2,x_4)+u_2q(x_2,x_4)$
modulo $f_*\cO_Y(-K_Y-E)$,
where $u_1,u_2$ are units.
Then by the coordinates change in which $x_3+u_1x_1p \mapsto x_3$,
we obtain a new identification (\ref{eqn:id}) in which
$h_{H_S}\equiv v_1x_3+x_1^{n-1}p^nr(x_1,x_2,x_4)+u_2q(x_2,x_4)$
modulo $f_*\cO_Y(-K_Y-E)$,
where $v_1$ is a unit.
Repeating this argument,
we can eventually obtain the desired coordinates.
\end{proof}

We proceed to apply Lemma \ref{lem:gen_mtd} to $f$.
Let $m_i/n$ denote the multiplicity of $x_i$ in Lemma \ref{lem:id} along $E$.
Then $m_3/n=a/n$ and $m_4/n=1$.
We distribute the weights by $\wt(x_1,x_2,x_3,x_4)=(m_1/n,m_2/n,a/n,1)$.
For a semi-invariant function $h$,
we let $\word h$ denote the weighted order of $h$ with respect to
these weights.
Set $D:=\word \phi$
and decompose $\phi$ as $\phi=\phi_D+\phi_{>D}$ as in Lemma \ref{lem:gen_mtd}.

\begin{lemma}\label{lem:phi_d_A}
$\word g=m_1/n+m_2/n$.
Consequently $D=m_1/n+m_2/n$ and $\phi_D=x_1x_2+g_D(x_3^n,x_4)$,
where $g_D$ is the weighted homogeneous part of weighted order $D$ in $g$.
\end{lemma}

\begin{proof}
Set $D_0:=\word g$
and let $g_{D_0}$ denote the weighted homogeneous part
of weighted order $D_0$ in $g$.
It suffices to prove that $D_0=m_1/n+m_2/n$.
If $D_0>m_1/n+m_2/n$ then $x_1x_2=-g_{D_0} \in f_*\cO_Y(-D_0E)$,
which is a contradiction.
Hence $D_0\le m_1/n+m_2/n$,
and it suffices to prove that $g_{D_0} \not\in f_*\cO_Y(-(D_0+1)E)$.
This $g_{D_0}$ always decomposes into a product of
the monomial $x_4^{D_0-a\rd{D_0/a}}$ and
$\rd{D_0/a}$ linear combinations of the monomials $x_3^n$ and $x_4^a$.
The function $x_4$ has multiplicity $1$ along $E$,
and any function defined by a non-zero linear combination of
$x_3^n$ and $x_4^a$ has multiplicity $a$
by Lemmata \ref{lem:pivotal} and \ref{lem:id}.
Hence $g_{D_0}$ has multiplicity $D_0$ along $E$,
and thus the lemma is deduced.
\end{proof}

By using Lemma \ref{lem:phi_d_A},
we shall verify the assumption of Lemma \ref{lem:gen_mtd} for $f$.
Set $N:=\bZ^4+\bZ(1/n,-1/n,b/n,0)$ and $\vec{e}:=(m_1/n,m_2/n,a/n,1)$.
It is obvious that the exceptional divisor $F$ of the weighted blow-up $g$
constructed in Lemma \ref{lem:gen_mtd} is irreducible and reduced.
\begin{enumerate}
\item
$N$ is generated by $\vec{e},\vec{e}_2,\vec{e}_3,\vec{e}_4$ and
\begin{align*}
(-1/n,1/n,-b/n,0)=
\frac{1}{m_1}\Bigl(-\vec{e}+\frac{m_1+m_2}{n}\vec{e}_2
+\frac{-bm_1+a}{n}\vec{e}_3+\vec{e}_4\Bigr).
\end{align*}
We see that
\begin{align*}
&\bar{U}_1=\bC^4_{x_1'x_2'x_3'x_4'}/
\frac{1}{m_1}\big(-1,\frac{m_1+m_2}{n},\frac{-bm_1+a}{n},1\big),\\
&(x_1,x_2,x_3,x_4)=
({x_1'}^{m_1/n},{x_1'}^{m_2/n}x_2',{x_1'}^{a/n}x_3',x_1'x_4').
\end{align*}
In particular $\Sing\bar{Z} \subseteq (x_1'=x_4'=0)$.
$Z|_{\bar{U}_1}$ is given by
\begin{align*}
x_2'+g_D({x_3'}^n,x_4')+x_1'\cdot(\something)=0,
\end{align*}
whence $\Sing Z|_{\bar{U}_1}$ is contained in
\begin{align*}
(x_1'=x_4'=x_2'+g_D({x_3'}^n,0)=0).
\end{align*}
\item
$\bar{U}_2=\bC^4_{x_1'x_2'x_3'x_4'}/
\frac{1}{m_2}(\frac{m_1+m_2}{n},-1,\frac{bm_2+a}{n},1)$.
$\Sing Z|_{\bar{U}_2}$ is contained in
\begin{align*}
(x_2'=x_4'=x_1'+g_D({x_3'}^n,0)=0).
\end{align*}
\item
$N$ is generated by $\vec{e}_1,\vec{e}_2,\vec{e},\vec{e}_4$
and $(-b'/n,b'/n,-1/n,0)$ for an integer $b'$ with $bb'\equiv1$ modulo $n$.
$\bar{U}_3=\bC^4_{x_1'x_2'x_3'x_4'}/
\frac{1}{a}(\frac{-b'a+m_1}{n},\frac{b'a+m_2}{n},-1,1)$.
$\Sing Z|_{\bar{U}_3}$ is contained in
\begin{align*}
\{x_4'\axis\} \cup (x_3'=x_4'=x_1'x_2'+g_D(1,0)=0).
\end{align*}
\item
$\bar{U}_4=
\bC^4_{x_1'x_2'x_3'x_4'}/\frac{1}{n}(1,-1,a,0)$.
$\Sing Z|_{\bar{U}_4}$ is contained in $x_3'$-axis.
\end{enumerate}
Therefore the assumption of Lemma \ref{lem:gen_mtd} holds,
and $\Sing Z \subseteq \bigcup_{1 \le i \le 4} \bar{D}_i$.
Since $c/n:=(m_1/n+m_2/n+a/n+1)-1-D=a/n$ by Lemma \ref{lem:phi_d_A},
we can deduce the following theorem from Lemma \ref{lem:gen_mtd}.

\begin{theorem}\label{thm:wtd_blow-up}
$f \cong g$ over $X$.
Hence $f$ is the weighted blow-up with weights
$\wt(x_1,x_2,x_3,x_4)=(m_1/n,m_2/n,a/n,1)$.
\end{theorem}

Therefore $Z\cong Y$ must be terminal.
$Z$ has a quotient singularity of type $\frac{1}{m_1}(-1,\frac{-bm_1+a}{n},1)$
at the origin of the chart $\bar{U}_1$,
and has a quotient singularity of type $\frac{1}{m_2}(-1,\frac{bm_2+a}{n},1)$
at the origin of the chart $\bar{U}_2$.
Hence $m_1$ is co-prime to $(-bm_1+a)/n$,
and $m_2$ is co-prime to $(bm_2+a)/n$.
On this condition,
$Z$ is terminal on the chart $\bar{U}_3$ only if
$Z$ does not pass through the origin of $\bar{U}_3$,
or equivalently $g_D(1,0)\neq0$.
Now it is obvious that $f$ can be expressed as in
Theorem \ref{thm:MT}(\ref{itm:MTord}).
We can also verify that
any $f$ in Theorem \ref{thm:MT}(\ref{itm:MTord}) is a divisorial contraction
by using Lemma \ref{lem:final}.
Therefore we complete Theorem \ref{thm:MT}(\ref{itm:MTord}).

\begin{lemma}\label{lem:final}
Let $a$, $n$, $b$, $b'$, $r_1$, $r_2$ be integers
which satisfy the following conditions.
\begin{enumerate}
\item
$bb' \equiv 1$, $a \equiv br_1$ modulo $n$,
and $r_1+r_2 \equiv 0$ modulo $an$.
\item
$(a-br_1)/n$ is co-prime to $r_1$.
\end{enumerate}
Then $(a+br_2)/n$ is co-prime to $r_2$,
and $(r_1-b'a)/n$, $(r_2+b'a)/n$ are co-prime to $a$.
\end{lemma}

\begin{proof}
Set $g:=\gcd(a,n)$, $a=ga'$, $n=gn'$.
Then $ga' \equiv br_1$ modulo $gn'$.
Since $b$ is co-prime to $n=gn'$, we have $g\mid r_1$ and thus $g\mid r_2$.
Set $r_1=gr_1'$, $r_2=gr_2'$.
Then $(a'-br_1')/n'$ is co-prime to $gr_1'$,
whence $\gcd(a'-br_1',r_1')\mid n'$.
This implies that $a'$ is co-prime to $r_1'$,
and thus $a'$ is co-prime to $r_2'$ by $r_1'+r_2' \equiv 0$ modulo $a'n'g$.
On the other hand by $g \mid -b(r_1'+r_2')/n'=(a'-br_1')/n'-(a'+br_2')/n'$,
we obtain that $(a'+br_2')/n'$ is also co-prime to $g$.
Therefore $(a+br_2)/n=(a'+br_2')/n'$ is co-prime to $r_2=gr_2'$.
We have seen that $a'$ is co-prime to $r_1'$.
Since $b(r_1'-b'a')/n' \equiv (br_1'-a')/n'$ modulo $g$,
we obtain that $(r_1'-b'a')/n'$ is co-prime to $g$.
Hence $(r_1-b'a)/n=(r_1'-b'a')/n'$ is co-prime to $a=ga'$.
Similarly, $(r_2+b'a)/n$ is co-prime to $a$.
\end{proof}

\section{The exceptional case of type c$D/2$}\label{sec:cD/2}
We have obtained Theorem \ref{thm:co-prime}
in Sections \ref{sec:co-primeness} and \ref{sec:small},
and this implies Theorem \ref{thm:MT}(\ref{itm:MTexc})
in the case where $f$ is of exceptional type by Corollary \ref{cor:a=1}.
We have proved Theorem \ref{thm:MT}(\ref{itm:MTord}) in Section \ref{sec:cA/n}.
In this section
we study the case (\ref{itm:hard_D}) of Theorem \ref{thm:hard_GE}
to complete the theorems stated in Section \ref{sec:introduction}.
We need to prove that in this case $P$ is of type c$D/2$
and $f$ is described as in Theorem \ref{thm:cD/2}.
It should be emphasised that
whichever type in Table \ref{tbl:classification} $f$ belongs to,
we can apply the method of explicit description of $f$ below
by formal modification.
Throughout this section,
we keep assuming that $f$ belongs to the case (\ref{itm:hard_D}) of
Theorem \ref{thm:hard_GE}.
In particular $P$ is of type c$Ax/2$ or c$D/2$ by Table \ref{tbl:reproduce}.
We first provide examples of Theorem \ref{thm:cD/2}.

\begin{example}\label{exl:cD/2_15g}
\begin{enumerate}
\item
Let $P \in X$ be the germ
\begin{align*}
o \in (x_1^2+x_2^2x_4+x_3^{2(r+1)/a}+x_4^{r+1}=0)
\subset \bC^4_{x_1x_2x_3x_4} /\frac{1}{2}(1,1,1,0),
\end{align*}
where $a \mid r+1$ and $a$, $r$ are odd with $r \neq 1$.
$P$ is of type c$D/2$.
Let $f$ be the weighted blow-up of $X$
with weights $\wt(x_1,x_2,x_3,x_4)=((r+2)/2,r/2,a/2,1)$.
Then $f$ is a divisorial contraction with discrepancy $a/2$.
$f$ is of No $\textrm{15}''$ with $J=\{(r,1),(r+2,1)\}$.
\item
Let $P \in X$ be the germ
\begin{align*}
o \in \bigg(
\begin{array}{c}
x_1^2+x_2x_5+x_4^{r+2}=0 \\
x_2x_4+x_3^{(r+2)/a}+x_5=0
\end{array}
\bigg)
\subset \bC^4_{x_1x_2x_3x_4x_5} /\frac{1}{2}(1,1,1,0,1),
\end{align*}
where $a \mid r+2$ and $(r+2)/a$ is odd with $r \neq 1$ and $a \neq r+2$.
$P$ is of type c$D/2$.
Let $f$ be the weighted blow-up of $X$
with weights $\wt(x_1,x_2,x_3,x_4,x_5)=((r+2)/2,r/2,a/2,1,(r+4)/2)$.
Then $f$ is a divisorial contraction with discrepancy $a/2$.
$f$ is of No $\textrm{15}''$ with $J=\{(r,1),(r+4,1)\}$.
\end{enumerate}
\end{example}

We start with computing the dimensions of the vector spaces
\begin{align*}
V^{[0]}(i) := & f_*\cO_Y(-iE)/f_*\cO_Y(-(i+1)E) &&(i \in \bZ),\\
V^{[1]}(i) := & f_*\cO_Y(K_Y-(i+a/2)E)          && \\
& \qquad\qquad /f_*\cO_Y(K_Y-(i+1+a/2)E)        && (i+a/2 \in\bZ).
\end{align*}
Define $V^{[j]}(i)$ for every $i \in \bZ/2$
by setting $V^{[j]}(i):=0$ if $i+ja/2 \not\in \bZ$.
The vector space $V^{[j]}(i)$ is regarded as the space of
semi-invariant functions on $X$
with multiplicity $i$ along $E$ and with weight $j$.
We set
\begin{align*}
V(i):= V^{[0]}(i) \oplus V^{[1]}(i).
\end{align*}
We define the finite set
\begin{align*}
N_k(i):=\Bigl\{\vec{l}=(l_1,l_2,l_3,l_4)\in\bZ_{\ge0}^4 \Bigm|
\frac{r_1+k}{2}l_1+\frac{r_1}{2}l_2+\frac{a}{2}l_3+l_4=i,\ l_1l_2=0\Bigr\}
\end{align*}
for $i \in \bZ/2$ and $k=2,4$,
and decompose it as
\begin{align*}
N_k(i)=N_k^{[0]}(i) \sqcup N_k^{[1]}(i),
\end{align*}
where
$N_k^{[j]}(i) :=\{\vec{l}\in N_k(i) \mid
l_1+l_2+l_3\equiv j \ \textrm{modulo} \ 2\}$.
Note that $N_k^{[j]}(i)=\emptyset$ when $i+ja/2 \not\in \bZ$
by $a\equiv r_1$ modulo $2$.

\begin{lemma}\label{lem:dimension}
$\dim V^{[j]}(i)=\#N_{r_2-r_1}^{[j]}(i)$.
\end{lemma}

\begin{proof}
We have
$\dim V^{[j]}(i)=d(j,-i-ja/2)$ for $i\ge -a/2$ with $i+ja/2 \in \bZ$
by (\ref{eqn:d(i,j)=dim}),
and $jb_{Q_k}-(i+ja/2)v_{Q_k}\equiv -(i+jr_k/2)$ modulo $r_k$ for $k=1,2$
by Theorem \ref{thm:hard_GE}(\ref{itm:hard_D}).
Set $s_k:=i+jr_k/2$.
Thus by (\ref{eqn:difference_of_d(i,j)}),
for any $i\ge0$ with $i+ja/2 \in \bZ$,
\begin{align*}
&\dim V^{[1-j]}(i-a/2)-\dim V^{[j]}(i) \\
&=-\Bigl(i+\frac{1}{2}\Bigr)\Bigl(\frac{1}{r_1}+\frac{1}{r_2}\Bigr)
+\sum_{k=1,2}
\Bigl( \frac{\overline{s_k}(r_k-\overline{s_k})}{2r_k}
- \frac{(\overline{s_k}+1)(r_k-\overline{s_k}-1)}{2r_k} \Bigr) \\
&=j-1-\sum_{k=1,2}\frac{s_k-\overline{s_k}}{r_k}.
\end{align*}
Hence for any $i\ge0$ with $i+ja/2 \in \bZ$,
\begin{align}\label{eqn:dimension_V}
\dim V^{[j]}(i)-\dim V^{[1-j]}(i-a/2)
=1-j+\sum_{k=1,2}\Bigl\lfloor \frac{i}{r_k}+\frac{j}{2} \Bigr\rfloor.
\end{align}
On the other hand,
we have a decomposition
\begin{align*}
N_{r_2-r_1}^{[j]}(i)=
&\{(l_1,l_2,l_3+1,l_4) \mid
\vec{l} \in N_{r_2-r_1}^{[1-j]}(i-a/2)\} \\
&\sqcup \{\vec{l}\in N_{r_2-r_1}^{[j]}(i) \mid l_1,l_3=0\} \\
&\sqcup \{\vec{l}\in N_{r_2-r_1}^{[j]}(i) \mid l_2,l_3=0,\ l_1\ge1\}.
\end{align*}
Hence for any $i\ge0$ with $i+ja/2 \in \bZ$,
\begin{align}\label{eqn:dimension_N}
&\#N_{r_2-r_1}^{[j]}(i)-\#N_{r_2-r_1}^{[1-j]}(i-a/2) \\
\nonumber=&
\#\{\vec{l}\in N_{r_2-r_1}^{[j]}(i) \mid l_1,l_3=0\}
+\#\{\vec{l}\in N_{r_2-r_1}^{[j]}(i) \mid l_2,l_3=0,\ l_1\ge1\} \\
\nonumber=&
\Bigl(1+\Bigl\lfloor \frac{i-jr_1/2}{r_1} \Bigr\rfloor\Bigr)
+\Bigl(1+\Bigl\lfloor \frac{i-(2-j)r_2/2}{r_2} \Bigr\rfloor\Bigr) \\
\nonumber=&
1-j+\sum_{k=1,2}\Bigl\lfloor \frac{i}{r_k}+\frac{j}{2} \Bigr\rfloor.
\end{align}
The lemma holds trivially if $i<0$ or $i+ja/2 \not\in \bZ$.
Therefore the lemma follows from
(\ref{eqn:dimension_V}) and (\ref{eqn:dimension_N}).
\end{proof}

We shall find bases of the vector space $V(i)$ in terms of any identification
\begin{align*}
P \in X \cong o \in (\phi=0) \subset \bC^4_{x_1x_2x_3x_4}/\frac{1}{2}(1,1,1,0).
\end{align*}
For a semi-invariant function $h$,
we let $\mult_E h$ denote the multiplicity of $h$ along $E$.

\begin{lemma}\label{lem:vector}
\begin{enumerate}
\item\label{itm:vector1}
$\mult_E x_4=1$ and $\mult_E x_i\ge a/2$ for any $i=1,2,3$.
There exists some $k=1,2,3$ with $\mult_E x_k=a/2$.
By permutation we may assume that $x_k=x_3$.
\item\label{itm:vector2}
The monomials $x_3^2$ and $x_4^a$ are linearly independent in $V^{[0]}(a)$.
For $i<r_1/2$,
the monomials $x_3^{l_3}x_4^{l_4}$ for
$(0,0,l_3,l_4)\in N_{r_2-r_1}(i)$
form bases of $V(i)$.
In particular for $k=1,2$,
$\mult_E \bar{x}_k \ge r_1/2$ for
$\bar{x}_k:=x_k+\sum_{(0,0,l_3,l_4)\in
\bigcup_{\frac{a}{2}\le i<\frac{r_1}{2}}N_{r_2-r_1}^{[1]}(i)}
c_{kl_3l_4}x_3^{l_3}x_4^{l_4}$
with some $c_{kl_3l_4}\in \bC$.
\item\label{itm:vector3}
There exists some $k=1,2$ with $\mult_E \bar{x}_k=r_1/2$ such that
the monomials $\bar{x}_k$ and $x_3^{l_3}x_4^{l_4}$ for
$(0,0,l_3,l_4)\in N_{r_2-r_1}(r_1/2)$ form bases of $V(r_1/2)$.
By permutation we may assume that $\bar{x}_k=\bar{x}_2$.
\item\label{itm:vector4}
One of the following holds.
\begin{enumerate}
\item
$\bar{x}_2x_4$ and $x_3^{l_3}x_4^{l_4}$ for
$(0,0,l_3,l_4)\in N_{r_2-r_1}^{[1]}(r_1/2+1)$
are linearly independent in $V^{[1]}(r_1/2+1)$.
In this case $\mult_E \hat{x}_1=r_2/2$ for
$\hat{x}_1:=\bar{x}_1+
\sum_{(0,l_2,l_3,l_4) \in
\bigcup_{\frac{r_1}{2}\le i<\frac{r_2}{2}}N_{r_2-r_1}^{[1]}(i)}
c_{l_2l_3l_4}\bar{x}_2^{l_2}x_3^{l_3}x_4^{l_4}$
with some $c_{l_2l_3l_4}\in \bC$.
For $i<(r_1+r_2)/2$,
the monomials $\hat{x}_1^{l_1}\bar{x}_2^{l_2}x_3^{l_3}x_4^{l_4}$ for
$\vec{l}\in N_{r_2-r_1}(i)$ form bases of $V(i)$.
\item
$\bar{x}_2x_4$ and $x_3^{l_3}x_4^{l_4}$ for
$(0,0,l_3,l_4)\in N_{r_2-r_1}^{[1]}(r_1/2+1)$
are linearly dependent in $V^{[1]}(r_1/2+1)$.
They have one non-trivial relation of form $\psi:=\bar{x}_2x_4+h(x_3,x_4)=0$
for which
the monomial $x_3^{(r_1+2)/a}$ appears in $h$ with non-zero coefficient.
In this case $r_2=r_1+4$,
$(r_1+2)/a$ is odd,
$\mult_E\psi=r_2/2$ and $\mult_E \hat{x}_1=r_1/2+1$ for
$\hat{x}_1:=\bar{x}_1+\sum_{(0,l_2,l_3,l_4)\in N_{2}^{[1]}(r_1/2)}
c_{l_2l_3l_4}\bar{x}_2^{l_2}x_3^{l_3}x_4^{l_4}$
with some $c_{l_2l_3l_4}\in \bC$.
For $i<(r_1+r_2)/2$,
there exists a natural exact sequence
\begin{align*}
0 \to \!\!\!\!\!\! &\bigoplus_{\vec{l}\in \tilde{N}_2(i+1-r_2/2)}\!\!\!\!\!\!
\bC \cdot \hat{x}_1^{l_1}\bar{x}_2^{l_2}x_3^{l_3}x_4^{l_4}\psi \to \\
&\quad\,\ \bigoplus_{\vec{l}\in \tilde{N}_2(i)} \
\bC \cdot \hat{x}_1^{l_1}\bar{x}_2^{l_2}x_3^{l_3}x_4^{l_4} \ \oplus
\!\!\!\!\!\!\!\!\!\!\!\! \bigoplus_{(0,0,l_3,l_4) \in \tilde{N}_2(i-r_2/2)}
\!\!\!\!\!\!\!\!\!\!\!\!
\bC \cdot x_3^{l_3}x_4^{l_4}\psi
\to V(i) \to 0,
\end{align*}
where $\tilde{N}_2(i):=\{\vec{l}\in\bZ_{\ge0}^4 \mid
(r_1/2+1)l_1+(r_1/2)l_2+(a/2)l_3+l_4=i\}$.
\end{enumerate}
\end{enumerate}
\end{lemma}

\begin{proof}
$\dim V^{[0]}(1)=1$,
$\dim V^{[1]}(a/2)=1$ and $\dim V^{[1]}(i)=0$ for $i<a/2$
by Lemma \ref{lem:dimension} and $a<r_1<r_2$
in Theorem \ref{thm:hard_GE}(\ref{itm:hard_D}).
This implies (\ref{itm:vector1}).

By permutation we may assume that $\mult_Ex_3=a/2$.
The linear independence of $x_3^2$ and $x_4^a$ is obtained
by the method in Lemma \ref{lem:pivotal}.
By $\dim V^{[0]}(1)=1$,
a general hyperplane section $H_X$ on $X$
has the defining function $x_4+g_H=0$ with $\mult_E g_H \ge 2$.
Similarly by $\dim V^{[1]}(a/2)=1$,
a general elephant $S_X$ of $X$
has the defining function $x_3+g_S=0$ with $\mult_E g_S \ge a/2+1$.
In Step \ref{stp:hard_GE15-1} in Theorem \ref{thm:hard_GE15}
and Step \ref{stp:15/16small-2} in Theorem \ref{thm:15/16small},
we have seen that we may choose the local coordinates
$x_{21},x_{22},x_{23}$ with weights $\wt(x_{21},x_{22},x_{23})=(1,-1,b_2)$
of the index-one cover $Q_2^\sharp \in Y^\sharp$ so that
the pre-images $H^\sharp$, $S^\sharp$
are given by respectively $x_{21}=0$, $x_{23}=0$.
We let $h$ denote a semi-invariant function with weight $-1$
on the germ $Q_2^\sharp \in Y^\sharp$ which defines $E$.
Then on $Q_2 \in Y$ the functions $x_3^2$, $x_4^a$ are expressed as
$x_3^2=u_3x_{23}^2h^a+p_3$ and $x_4^a=u_4x_{21}^ah^a+p_4$,
where $u_3$, $u_4$ are units and $p_3,p_4\in(h^{a+1})$.
Then,
in order to obtain the linear independence of $x_3^2,x_4^a$,
we need to show that $h$ does not divide
any function $c_1x_{23}^2+c_2x_{21}^a$
for units $c_1,c_2$ with $(c_1,c_2)\neq(0,0)$.
This follows from
$\wt x_{23}, \wt x_{23}^2 \not\equiv -1$ modulo $r_2$
by $b_2=(a+r_2)/2$ with $a<r_1\le r_2-2$.

Because of the linear independence of $x_3^2,x_4^a$,
for any $i \in \bZ/2$ the monomials $x_3^{l_3}x_4^{l_4}$
for $(0,0,l_3,l_4)\in N_{r_2-r_1}(i)$ are linearly independent in $V(i)$.
For,
every linear combination of the monomials $x_3^{l_3}x_4^{l_4}$
for $(0,0,l_3,l_4)\in N_{r_2-r_1}^{[j]}(i)$ decomposes
into a product of $x_3^j$, $x_4^{i-(j/2+\rd{i/a-j/2})a}$ and
$\rd{i/a-j/2}$ linear combinations of $x_3^2,x_4^a$.
Set $W(i)$ as the linear subspace of $V(i)$ spanned by these monomials,
and decompose it as $W(i)=W^{[0]}(i) \oplus W^{[1]}(i)$
by setting $W^{[j]}(i):=W(i) \cap V^{[j]}(i)$.
Then $\dim W^{[j]}(i)=\#N_{r_2-r_1}^{[j]}(i)$ for $i<r_1/2$ and
thus we obtain (\ref{itm:vector2}) by Lemma \ref{lem:dimension}.

The vector space $V(r_1/2)$ is spanned by
$\bar{x}_1$, $\bar{x}_2$ and $W(r_1/2)$ by (\ref{itm:vector2}).
Since $\dim V(r_1/2)=\dim W(r_1/2)+1$ by Lemma \ref{lem:dimension},
we obtain (\ref{itm:vector3}).
By permutation we assume that
$\bar{x}_2$ forms a basis of $V(r_1/2)/W(r_1/2)\cong\bC$.
Then $\mult_E x_1' \ge r_1/2+1$ for
$x_1':=\bar{x}_1+\sum_{(0,l_2,l_3,l_4)\in N_{2}^{[1]}(r_1/2)}
c_{l_2l_3l_4}\bar{x}_2^{l_2}x_3^{l_3}x_4^{l_4}$
with some $c_{l_2l_3l_4}\in \bC$,
and $V(r_1/2+1)$ is spanned by the monomials
${x_1'}^{l_1}\bar{x}_2^{l_2}x_3^{l_3}x_4^{l_4}$ with $\vec{l}\in N_2(r_1/2+1)$.
If $r_2=r_1+2$ then
by $\dim V(r_2/2)=\#N_2(r_2/2)$ in Lemma \ref{lem:dimension},
the monomials above form bases of $V(r_2/2)$.
Setting $\hat{x}_1:=x_1'$,
we obtain the claim (\ref{itm:vector4})(a) for $i\le r_2/2$.
The rest of (\ref{itm:vector4})(a) follows
from this and Lemma \ref{lem:dimension} inductively.
If $r_2=r_1+4$ then
$\dim V^{[0]}(r_1/2+1)=\#N_4^{[0]}(r_1/2+1)=\#N_2^{[0]}(r_1/2+1)$ but
$\dim V^{[1]}(r_1/2+1)=\#N_4^{[1]}(r_1/2+1)=\#N_2^{[1]}(r_1/2+1)-1$.
Hence the monomials ${x_1'}$,
$\bar{x}_2x_4$ and $x_3^{l_3}x_4^{l_4}$ in $W^{[1]}(r_1/2+1)$
have one non-trivial relation, say $\psi=0$, in $V^{[1]}(r_1/2+1)$.
If $x_1'$ appears in $\psi$ with non-zero coefficient,
then $\mult_E \hat{x}_1 \ge r_2/2$ for
$\hat{x}_1:=x_1'+\sum_{(0,l_2,l_3,l_4)\in N_{4}^{[1]}(r_1/2+1)}
c_{l_2l_3l_4}\bar{x}_2^{l_2}x_3^{l_3}x_4^{l_4}$
with some $c_{l_2l_3l_4}\in \bC$,
and eventually we obtain the case (\ref{itm:vector4})(a)
as in the case $r_2=r_1+2$.
If $x_1'$ does not appear in $\psi$,
then by (\ref{itm:vector3}) $\psi$ is
of form described in the statement (\ref{itm:vector4})(b) up to constant.
Set $\hat{x}_1:=x_1'$.
The vector space $V^{[1]}(r_2/2)$ is spanned by the function $\psi$ and
the monomials $\hat{x}_1^{l_1}\bar{x}_2^{l_2}x_3^{l_3}x_4^{l_4}$ for
$\vec{l}\in N_2^{[1]}(r_2/2)$
with a relation $x_4\psi=0$ in $V^{[1]}(r_2/2)$.
Since $\dim V^{[1]}(r_2/2)=\#N_4^{[1]}(r_2/2)=\#N_2^{[1]}(r_2/2)$
in Lemma \ref{lem:dimension},
the function $\psi$ does not vanish in $V^{[1]}(r_2/2)$.
This means that $\mult_E \psi=r_2/2$.
Now we obtain the claim (\ref{itm:vector4})(b) inductively,
noticing the estimate $\dim V(i)=\#N_4(i)$ in Lemma \ref{lem:dimension} and
the equality
$\#N_4(i)=\#\tilde{N}_2(i)+\#\tilde{N}_2(i-r_2/2)-\#\tilde{N}_2(i+1-r_2/2)$
for $i<(r_1+r_2)/2$ by direct calculation.
\end{proof}

We derive that $P$ is of type c$D/2$ as an easy corollary.
This completes Theorems \ref{thm:MT}(\ref{itm:MTexc}).

\begin{corollary}\label{cor:cD/2}
$P$ is of type c$D/2$.
\end{corollary}

\begin{proof}
$P$ is of type c$Ax/2$ or c$D/2$ by Table \ref{tbl:reproduce}.
Suppose that $P$ is of type c$Ax/2$.
Then thanks to \cite{Mo85},
we have an identification
\begin{align*}
P \in X \cong o \in (x_1^2+x_4^2+g(x_2,x_3)=0)
\subset \bC^4_{x_1x_2x_3x_4}/\frac{1}{2}(1,1,1,0)
\end{align*}
with $g\in(x_2,x_3)^4$.
$\mult_E x_4=1$ and $\mult_E x_i\ge a/2$ for any $i=1,2,3$
by Lemma \ref{lem:vector}(\ref{itm:vector1}).
Hence $a=2$ and $\mult_E x_1=1$
by the relation $-x_4^2=x_1^2+g$ on $X$ with $\mult_E(x_1^2+g) \ge a$.
Then $\mult_E(x_1^2+x_4^2)=\mult_E g \ge 4$,
which contradicts Lemma \ref{lem:vector}(\ref{itm:vector2}).
\end{proof}

In the rest of this section
we prove Theorem \ref{thm:cD/2} by applying Lemma \ref{lem:gen_mtd},
as in Section \ref{sec:cA/n} and \cite[Section 6]{Ka03}.
According to \cite[Remark 23.1]{Mo85} by Mori,
we can identify $P \in X$ with one of the following
in the space $\bC^4_{x_1x_2x_3x_4} /\frac{1}{2}(1,1,1,0)$.
\begin{align}
&o \in (x_1^2+x_2x_3x_4+x_2^{2\alpha}+x_3^{2\beta}+x_4^\gamma=0),
\label{eqn:phi1} \\
&o \in (x_1^2+x_2^2x_4+\lambda x_2x_3^{2\alpha-1}+g(x_3^2,x_4)=0),
\label{eqn:phi2}
\end{align}
where $\alpha, \beta\ge2$, $\gamma\ge3$,
$\lambda \in \bC$ and $g \in (x_3^4,x_3^2x_4^2,x_4^3)$.
Note that
\begin{align}\label{eqn:ord_estimate}
\frac{r_1+r_2}{2} \le
\begin{cases}
\gamma        & \textrm{in (\ref{eqn:phi1})},\\
\ord g(0,x_4) & \textrm{in (\ref{eqn:phi2})}.
\end{cases}
\end{align}
because a general elephant $S_X$ is of type $D_k$
with $k \ge r_1+r_2$ by the dual graph of $f_S$
in Theorem \ref{thm:hard_GE}(\ref{itm:hard_D}).
The crucial step is to obtain the expression of the germ $P \in X$
described in Theorem \ref{thm:cD/2}.

\begin{lemma}\label{lem:description}
The germ $P \in X$ has the form described in
either \textup{(\ref{itm:cD/2_2})} or \textup{(\ref{itm:cD/2_4})}
in Theorem \textup{\ref{thm:cD/2}}
such that each $\mult_E x_i$ coincides with the weight of $x_i$
distributed in Theorem \textup{\ref{thm:cD/2}}.
\end{lemma}

\begin{proof}
First consider the case (\ref{eqn:phi1}).
We prove that this case does not happen.
By Lemma \ref{lem:vector}(\ref{itm:vector1}),
$\mult_E x_4=1$ and $\mult_E x_i\ge a/2$ for any $i=1,2,3$
with some $\mult_E x_i=a/2$.
$\mult_E x_1 \ge a/2+1$
by the relation $-x_1^2=x_2x_3x_4+x_2^{2\alpha}+x_3^{2\beta}+x_4^\gamma$
and (\ref{eqn:ord_estimate}).
Thus we may assume that $\mult_E x_3=a/2$ without loss of generality.
We construct $\bar{x}_1, \bar{x}_2$
as in Lemma \ref{lem:vector}(\ref{itm:vector2}).

Set $W(i)$ as the linear subspace of $V(i)$ spanned by the monomials
in $x_3$ and $x_4$,
as in the proof of Lemma \ref{lem:vector}.
If $\bar{x}_1 \not\in W(r_1/2)$,
then the triple $(\bar{x}_1,x_3,x_4)$ plays the role
of $(\bar{x}_2,x_3,x_4)$ in Lemma \ref{lem:vector}(\ref{itm:vector3}).
We construct $\hat{x}_2$ as in Lemma \ref{lem:vector}(\ref{itm:vector4})
and obtain the quartuple $(\hat{x}_2,\bar{x}_1,x_3,x_4)$.
Set $\bar{x}_1:=x_1+p_1(x_3,x_4)$
and $\hat{x}_2:=x_2+p_2(\bar{x}_1,x_3,x_4)$.
Then (\ref{eqn:phi1}) is rewritten as
\begin{align*}
o \in ((\bar{x}_1-p_1)^2+(\hat{x}_2-p_2)x_3x_4+
(\hat{x}_2-p_2)^{2\alpha}+x_3^{2\beta}+x_4^\gamma=0).
\end{align*}
Hence there exists a relation $\bar{x}_1^2+(\something)=0$ in $V(r_1)$,
which contradicts Lemma \ref{lem:vector}(\ref{itm:vector4}).
If $\bar{x}_1 \in W(r_1/2)$,
then the triple $(\bar{x}_2,x_3,x_4)$ satisfies the condition
in Lemma \ref{lem:vector}(\ref{itm:vector3}).
We construct $\hat{x}_1$ as in Lemma \ref{lem:vector}(\ref{itm:vector4}),
and set $\bar{x}_2:=x_2+p_2(x_3,x_4)$
and $\hat{x}_1:=x_1+p_1(\bar{x}_2,x_3,x_4)$.
Note that the monomials $\bar{x}_2$, $x_3$ do not appear in $p_1$
by $\mult_E x_1 \ge a/2+1$ and the assumption $\bar{x}_1 \in W(r_1/2)$.
Then (\ref{eqn:phi1}) is rewritten as
\begin{align*}
o \in ((\hat{x}_1-p_1)^2+(\bar{x}_2-p_2)x_3x_4+
(\bar{x}_2-p_2)^{2\alpha}+x_3^{2\beta}+x_4^\gamma=0).
\end{align*}
Hence there exists a relation $\bar{x}_2x_3x_4+(\something)=0$
in $V((r_1+a+2)/2)$.
Thus by Lemma \ref{lem:vector}(\ref{itm:vector4}),
the case (b) in Lemma \ref{lem:vector}(\ref{itm:vector4}) holds
and this relation must equal $(x_3+(\something))\psi=0$.
Then there exists a relation $x_3\psi+(\something)=0$
in $V((r_1+a+4)/2)$,
whence $(r_1+a+4)/2 \ge (r_1+r_2)/2$
by Lemma \ref{lem:vector}(\ref{itm:vector4})(b).
This inequality contradicts that $a<r_1=r_2-4$.

We then consider the case (\ref{eqn:phi2}).
As in the case (\ref{eqn:phi1}),
we have $\mult_E x_4=1$ and $\mult_E x_1 \ge a/2+1$.
Considering the relation $-x_2^2x_4=x_1^2+\lambda x_2x_3^{2\alpha-1}+g$,
we also have $\mult_E x_2 \ge a/2+1$.
Thus $\mult_E x_3=a/2$ by Lemma \ref{lem:vector}(\ref{itm:vector1}).
Construct $\bar{x}_1, \bar{x}_2$ and $W(i)$ as in the case (\ref{eqn:phi1}).
The case $\bar{x}_1 \not\in W(r_1/2)$ does not happen
by the same reason as in the case (\ref{eqn:phi1}),
whence $\bar{x}_1 \in W(r_1/2)$.
Construct $\hat{x}_1$,
and set $\bar{x}_2:=x_2+p_2(x_3,x_4)$
and $\hat{x}_1:=x_1+p_1(\bar{x}_2,x_3,x_4)$ in the same way.
The monomials $\bar{x}_2$, $x_3$ do not appear in $p_1$.
Then (\ref{eqn:phi2}) is rewritten as
\begin{align}\label{eqn:phi2_rewritten}
o \in ((\hat{x}_1-p_1)^2+(\bar{x}_2-p_2)^2x_4+
\lambda(\bar{x}_2-p_2)x_3^{2\alpha-1}+g(x_3^2,x_4)=0).
\end{align}

If $r_2=r_1+2$,
then Lemma \ref{lem:vector}(\ref{itm:vector4})(a) holds
and the defining function in (\ref{eqn:phi2_rewritten}) must consist of
the monomials $\hat{x}_1^{l_1}\bar{x}_2^{l_2}x_3^{l_3}x_4^{l_4}$
for $\vec{l}\in \bigcup_{i\ge (r_1+r_2)/2}\tilde{N}_2(i)$.
In particular $p_2=0$.
$p_1$ is a function of $x_3,x_4$ by $\bar{x}_1 \in W(r_1/2)$.
Hence (\ref{eqn:phi2_rewritten}) is written
as in the case (\ref{itm:cD/2_2}) of Theorem \ref{thm:cD/2}.
The weighted homogeneous part of $\phi$ of weight $(r_1+r_2)/2$
has multiplicity $>(r_1+r_2)/2$ along $E$,
whence it is irreducible by Lemma \ref{lem:vector}(\ref{itm:vector4})(a).
If $r_2=r_1+4$,
then either (a) or (b) in Lemma \ref{lem:vector}(\ref{itm:vector4}) holds.
Since there exists a relation $\bar{x}_2^2x_4+(\something)=0$ in $V(r_1+1)$
by (\ref{eqn:phi2_rewritten}),
the case (a) does not happen
and in the case (b) this relation must equal $(x_2+(\something))\psi=0$.
Note that $p_1$ is a function of $x_3,x_4$ in the case (b).
By replacing $x_3$ with its scalar multiple,
we can set $\psi=\bar{x}_2x_4+x_3^{(r_1+2)/a}+q(x_3^2,x_4)x_3x_4$.
Then (\ref{eqn:phi2_rewritten}) is rewritten as
\begin{multline*}
o \in (\hat{x}_1^2-2\hat{x}_1p_1+(\bar{x}_2-2p_2-qx_3)\psi \\
+\bar{x}_2(\lambda x_3^{2\alpha-1}-x_3^{(r_1+2)/a})+\tilde{g}(x_3^2,x_4)=0).
\end{multline*}
By this expression,
Lemma \ref{lem:vector}(\ref{itm:vector4})(b) implies that
the weighted orders of
$p_1$, $2p_2+qx_3$, $\lambda x_3^{2\alpha-1}-x_3^{(r_1+2)/a}$, $\tilde{g}$
with weights $\wt(x_3,x_4)=(a/2,1)$
are respectively $\ge r_1/2+1$, $\ge r_1/2$, $\ge r_2/2$, $\ge (r_1+r_2)/2$.
Hence $p_1=p_2=0$ and $\lambda x_3^{2\alpha-1}=x_3^{(r_1+2)/a}$.
Replacing the coordinate $\bar{x}_2$ with $\bar{x}_2-qx_3$
and setting $x_5:=-\psi$,
we obtain the expression of $P \in X$
in Theorem \ref{thm:cD/2}(\ref{itm:cD/2_4}).
\end{proof}

It is the time to complete Theorem \ref{thm:cD/2}.
By Lemma \ref{lem:description},
we have the description of the germ $P \in X$ in
either (\ref{itm:cD/2_2}) or (\ref{itm:cD/2_4}) in Theorem \ref{thm:cD/2}.
We take this description and apply Lemma \ref{lem:gen_mtd}.
Of course we need to extend Lemma \ref{lem:gen_mtd}
in treating the case (\ref{itm:cD/2_4}) of Theorem \ref{thm:cD/2},
where $P \in X$ is given by a cyclic quotient
of a complete intersection three-fold in $\bC^5$,
but this extension is straightforward.
We let $g \colon (Z \supset F) \to (X \ni P)$ denote the weighted blow-up
with weights $\wt x_i=\mult_E x_i$,
induced from $\bar{g} \colon (\bar{Z} \supset \bar{F}) \to
(\bar{X}:=\bC_{x_1x_2x_3x_4}^4/\frac{1}{2}(1,1,1,0) \ni P)$
in the case of Theorem \ref{thm:cD/2}(\ref{itm:cD/2_2}),
and from $\bar{g} \colon (\bar{Z} \supset \bar{F}) \to
(\bar{X}:=\bC_{x_1x_2x_3x_4x_5}^5/\frac{1}{2}(1,1,1,0,1) \ni P)$
in the case of Theorem \ref{thm:cD/2}(\ref{itm:cD/2_4}).
By direct calculation as in the paragraph before Theorem \ref{thm:wtd_blow-up},
we verify the assumption of Lemma \ref{lem:gen_mtd} and
the implication $\Sing Z \subset \bigcup_i \bar{D}_i$.
Therefore by Lemma \ref{lem:gen_mtd},
we obtain that $f$ coincides with $g$ over $X$.
This completes Theorem \ref{thm:cD/2}.

\section*{Erratum to the previous paper}
We correct the omission of
the case (\ref{itm:hard_D}) of Theorem \ref{thm:hard_GE} and two minor gaps
in the study of the Gorenstein case in the previous paper \cite{Ka03},
by the methods used in this paper.
First we state the additional case in the classification.

\begin{addendum}
Let $f \colon (Y \supset E) \to (X \ni P)$
be a germ of a three-fold divisorial contraction
whose exceptional divisor $E$ contracts to a Gorenstein point $P$ of type c$D$.
Then $f$ can be of type ${\textrm{IIb}}^{\vee\vee}$ in \cite[Theorem 1.4]{Ka03}
described below,
in addition to the cases stated in \cite[Theorem 1.8, Corollary 1.15]{Ka03}.
\begin{enumerate}
\item\label{itm:D_1}
$f$ belongs to the case \textup{(a)}
in Theorem \textup{\ref{thm:hard_GE}(\ref{itm:hard_D})}.
There exists a suitable identification of $P \in X$ with
\begin{align*}
o \in (\phi:=x_1^2+x_1q(x_3,x_4)+x_2^2x_4+
\lambda x_2x_3^2+\mu x_3^3+p(x_2,x_3,x_4)=0)
\end{align*}
in the space $\bC^4_{x_1x_2x_3x_4}$ with $p \in (x_2,x_3,x_4)^4$,
such that $f$ is the weighted blow-up with weights
$\wt(x_1,x_2,x_3,x_4)=(r+1,r,a,1)$
which satisfies the following conditions.
\begin{enumerate}
\item
$a \mid 2r+1$ with $a\neq1,2r+1$.
\item
$\phi$ has weighted order $2r+1$ with the weights distributed above.
Its weighted homogeneous part of weight $2r+1$ is irreducible.
In fact $x_1q$ has weighted order $2r+1$ unless $q=0$.
\end{enumerate}
\item\label{itm:D_2}
$f$ belongs to the case \textup{(b)}
in Theorem \textup{\ref{thm:hard_GE}(\ref{itm:hard_D})}.
There exists a suitable identification of $P \in X$ with
\begin{align*}
o \in \bigg(
\begin{array}{c}
\phi_1:=x_1^2+x_2x_5+p(x_2,x_3,x_4)=0 \\
\phi_2:=x_2x_4+x_3^{(r+1)/a}+q(x_3,x_4)x_4+x_5=0
\end{array}
\bigg)
\end{align*}
in the space $\bC^5_{x_1x_2x_3x_4x_5}$ with $p \in (x_2,x_3,x_4)^4$,
such that $f$ is the weighted blow-up with weights
$\wt(x_1,x_2,x_3,x_4,x_5)=(r+1,r,a,1,r+2)$
which satisfies the following conditions.
\begin{enumerate}
\item
$a \mid r+1$ with $a\neq1,r+1$.
\item
$\phi_1$ has weighted order $2r+2$ with the weights distributed above.
Similarly $qx_4$ has weighted order $r+1$ unless $q=0$.
\end{enumerate}
\end{enumerate}
\end{addendum}

There exist examples of this addendum.

\begin{example*}
\begin{enumerate}
\item
Let $P \in X$ be the germ
\begin{align*}
o \in (x_1^2+x_2^2x_4+x_3^{(2r+1)/a}+x_4^{2r+1}=0)
\subset \bC^4_{x_1x_2x_3x_4},
\end{align*}
where $a \mid 2r+1$ with $a \neq 1,2r+1$.
$P$ is of type c$D$.
Let $f$ be the weighted blow-up of $X$
with weights $\wt(x_1,x_2,x_3,x_4)=(r+1,r,a,1)$.
Then $f$ is a divisorial contraction with discrepancy $a$.
$f$ is of type ${\textrm{IIb}}^{\vee\vee}$ with $J=\{(r,1),(r+1,1)\}$.
\item
Let $P \in X$ be the germ
\begin{align*}
o \in \bigg(
\begin{array}{c}
x_1^2+x_2x_5+x_4^{2r+2}=0 \\
x_2x_4+x_3^{(r+1)/a}+x_5=0
\end{array}
\bigg)
\subset \bC^4_{x_1x_2x_3x_4x_5},
\end{align*}
where $a \mid r+1$ with $a \neq 1,r+1$.
$P$ is of type c$D$.
Let $f$ be the weighted blow-up of $X$
with weights $\wt(x_1,x_2,x_3,x_4,x_5)=(r+1,r,a,1,r+2)$.
Then $f$ is a divisorial contraction with discrepancy $a$.
$f$ is of type ${\textrm{IIb}}^{\vee\vee}$ with $J=\{(r,1),(r+2,1)\}$.
\end{enumerate}
\end{example*}

The omission is due to an error in the proof of \cite[Lemma 3.14]{Ka03}.
The triple $(x,y,z)=(1,0,2)$ also satisfies the equality
$(x-1+y)/r_1=(x+1-z)/r_2$ appearing in the proof,
and the addendum is derived from this case.
We first deduce Theorem \ref{thm:hard_GE}(\ref{itm:hard_D}),
and then obtain the description of $f$
by using Theorem \ref{thm:hard_GE}(\ref{itm:hard_D}).

Theorem \ref{thm:hard_GE}(\ref{itm:hard_D}) in this case is proved
in the same fashion as Theorem \ref{thm:hard_GE15}.
However we need to obtain the general elephant theorem
in the course of the proof,
because Theorem \ref{cor:easy_GE_D} does not hold in the Gorenstein case.
The proof is modified as follows.
We first see the statement on $Q_1^\sharp$
in Step \ref{stp:hard_GE15-1} in Theorem \ref{thm:hard_GE15},
from which Step \ref{stp:hard_GE15-3} in Theorem \ref{thm:hard_GE15}
is deduced.
We then see that $H^\sharp$ is smooth at $Q_2^\sharp$.
This is the case unless $V\cong\bP^1$ with $a=2$ or $(r_2+1)/2$
by $(H \cdot V)_{Q_2}=2/r_2$.
If $a=2$ with $V\cong\bP^1$,
then there exist new semi-invariant coordinates $x_1,x_2,x_3$ at $Q_2^\sharp$
with $(x_1,x_2,x_3)|_{V^\dag}=(t^{2/r_2},t^{1-2/r_2},0)$,
whence $H^\sharp$ is smooth at $Q_2^\sharp$ by $(H \cdot V)_{Q_2}=2/r_2$.
If $a=(r_2+1)/2$,
then $(a,r_1,r_2)=(3,4,5)$ and $[H \cap E]=3[C]$ by Step \ref{stp:hard_GE15-3}.
$E^\sharp$ is not smooth at $Q_2^\sharp$
by $(x_{21},x_{22},x_{23})|_{C^\dag}=(0,t^{1/r_2},0)$ and $C \subset E$,
whence $H^\sharp$ is smooth at $Q_2^\sharp$
by $(H \cdot E \cdot [x_{22}=0])_{Q_2}=3/5$.
Thus $H^\sharp$ is smooth at $Q_2^\sharp$,
and Step \ref{stp:hard_GE15-2} in Theorem \ref{thm:hard_GE15} is deduced.
This implies Step \ref{stp:hard_GE15-1} in Theorem \ref{thm:hard_GE15}.
We also obtain the general elephant theorem
since $\system{-K_Y}$ is free outside $C$.
Hence Theorem \ref{thm:hard_GE}(\ref{itm:hard_D}) in this case is proved.

$S_X$ has a Du Val singularity of type $D$ at $P$,
whence $P$ is c$A$ or c$D$.
The classification in the c$A$ case in \cite[Section 6]{Ka03} works
even if this extra case is taken into account,
because we do not need the result of \cite[Section 5]{Ka03}
in the argument before \cite[Lemma 6.6]{Ka03} whenever $a<r_1$.
Hence we obtain that $P$ is c$D$
in the case (\ref{itm:hard_D}) of Theorem \ref{thm:hard_GE}.
Now Addendum is proved in the same way as in Section \ref{sec:cD/2}.
First obtain the identification of $P \in X$ with
$o \in ((\bar{x}_1-p_1(x_3,x_4))^2+g(x_2,x_3,x_4)=0)$
with $g\in(x_2,x_3,x_4)^3$ and $\mult_E \bar{x}_1, \mult_E x_2\ge r_1$,
$\mult_E x_3=a$, $\mult_E x_4=1$.
The pivotal remark is that $g$ contains one of the monomials
$x_2^2x_4,x_2x_3x_4,x_3^2x_4,x_2x_4^2,x_3x_4^2,x_4^3$
since $S_X$ has a Du Val singularity of type $D$ at $P$.
Then the rest of the proof is similar.

\medskip
Below we fill two minor gaps in \cite{Ka03}.
Firstly we correct the proof of $c_2=1$ in \cite[Lemma 3.11(i)]{Ka03}.
Its proof implies only that $\overline{ac_2}=ac_2$.
Suppose that $c_2\ge2$.
By the argument in Step \ref{stp:true/false-3} in Lemma \ref{lem:true/false},
we obtain that $a=2$, $r_1=3$ and $r_2$ is an odd integer $>3c_2$.
We then compute
the length of the cokernel of the map \cite[(2.7)]{Ka03} at $Q_2$.
This is $c_2+\min\{0,(r_2+1)/2-2c_2\}$ for $c_2<r_2/3$,
and is bounded by two by \cite[Lemma 2.10]{Ka03}.
Hence $c_2=2$.
By $(H \cdot C)=1/3-2/r_2$ and $(H \cdot H \cdot E)=(1/3+1/r_2)/2$,
we have $[H \cap E]=\alpha[C]+\beta[C']$ with
$(r_2,\alpha,\beta)=(7,5,0)$, $(7,2,1)$, $(9,2,0)$, $(9,1,1)$ or $(15,1,0)$,
where $C'\cong\bP^1$ is a curve intersecting $C$ at $Q_2$.
$H^\sharp$ is not smooth at $Q_1^\sharp$
by $(x_{11},x_{12},x_{13})|_{C^\dag}=(t^{1/3},0,0)$ and $C \subset H$,
whence $\alpha \neq 1$.
If $(r_2,\alpha,\beta)=(7,2,1)$
then $(x_{21},x_{22},x_{23})|_{C^\dag}=(0,t^{2/7},t^{3/7})$
and ${C'}^\dag$ is given by ${x'}_{21}$-axis
for new coordinates ${x'}_{21},{x'}_{22},{x'}_{23}$.
For a suitable surface $T$ on $Q_2 \in Y$ given by
$x_{21}^2-cx_{23}=0$ with $c \in \bC$,
$T$ intersects $C,C'$ properly with
$(T \cdot H \cdot E)_{Q_2}=(T \cdot 2C+C')_{Q_2}=8/7$.
This contradicts that $H^\sharp$, $E^\sharp$ at $Q_2^\sharp$
have multiplicities $\ge 3$  by $C,C' \subset H,E$.
Hence $(r_2,\alpha,\beta)=(7,5,0)$ or $(9,2,0)$.
We obtain the coordinates $x_{i1},x_{i2},x_{i3}$ of $Q_i^\sharp \in Y^\sharp$
such that $Q_1^\sharp \in E^\sharp$ is given by $x_{12}=0$ and 
$Q_2^\sharp \in H^\sharp$ is given by $x_{21}=0$.
Take a general elephant $S$.
By \cite[Lemma 3.4]{Ka03} we have
$[S \cap E]=\gamma[C]+[D_\delta]$ with
$(r_2,\gamma,\delta)=(7,3,1)$, $(7,4,2)$, $(7,1,2)$, $(9,1,1)$ or $(9,1,2)$,
where $D_\delta$ is a union of $\bP^1$ intersecting $C$ at $Q_\delta$.
We can remove the cases $(r_2,\gamma,\delta)=(7,1,2)$ and $(9,1,2)$
using a surface $T$ on $Q_2\in Y$ given by
$x_{21}-cx_{22}x_{23}=0$ with $c\in\bC$.
If $\delta=1$ then $D_1 \cong \bP^1$,
$(H \cdot D_1)_{Q_1}=1/3$ and $D_1 \not\subseteq H$.
This contradicts that $H^\sharp$ is not smooth at $Q_1^\sharp$.
Since $S$ is a general element of $\system{2H}$,
the defining function of $S^\sharp$ at $Q_2^\sharp$
contains the monomial $x_{21}^2$ with non-zero coefficient.
In the case $(r_2,\gamma,\delta)=(7,4,2)$,
this expression contradicts that $C,D \subset S$.

Secondly the smoothness of $E$ at the generic point of $C$ is indispensable
in the argument after the equation (3.6) in \cite[Theorem 3.15]{Ka03}.
Set $[S \cap E]=x[C]+[D]$ as in the proof.
We can express $E^\sharp$ on $Q_2^\sharp \in Y^\sharp$ as
$x_{23}\cdot(\something)+x_{21}^xy=0$,
where $y$ is the semi-invariant function with weight $-1$
for which $D^\sharp$ is defined by $x_{23}=y=0$.
In particular we have $r_2 \mid x$ by considering weights.
Then $r_2=r_1+1$ by (3.6) in \cite[Theorem 3.15]{Ka03},
$a<r_1<r_2$ and $[H \cap S]=[C]+[B]$ with $(E \cdot B)_{Q_1}=2/r_1$.
Thus $E^\sharp$ at $Q_1^\sharp$ is smooth unless
$a=r_1-2$ or $(r_1-1)/2$.
If $a=r_1-2$ or $(r_1-1)/2$ then by $a \mid r_1+r_2$
eventually $(r_1,r_2)=(7,8)$ and
$(a,u)=(5,3)$ or $(3,5)$ with $[H \cap E]=u[C]$.
Hence $E^\sharp$ at $Q_1^\sharp$ is smooth
by $(H \cdot E \cdot [x_{11}=0])_{Q_1}=u/7$,
and thus $E$ is smooth at the generic point of $C$.

\bibliographystyle{amsplain}

\end{document}